\documentclass[onefignum,onetabnum]{siamart190516} % review,

\usepackage{lipsum}
\usepackage{amsfonts}
\usepackage{graphicx}
\usepackage{epstopdf}
\usepackage{algorithmic}
\usepackage{microtype}

\usepackage{amsmath,amsfonts}
\usepackage{url}
\usepackage{comment}

\usepackage{longtable}
\usepackage{caption}
\usepackage{multirow}
\usepackage{xcolor}
\usepackage{xspace}

\usepackage{graphicx}

\ifpdf
  \DeclareGraphicsExtensions{.eps,.pdf,.png,.jpg}
\else
  \DeclareGraphicsExtensions{.eps}
\fi

\newcommand{\bmat}[1]{\begin{bmatrix}#1\end{bmatrix}} %
\newcommand{\inv}{^{-1}}
\newcommand{\norm}[1]{\|#1\|}
\newcommand{\Norm}[1]{\left\|#1\right\|}

\newcommand{\TR}[1]{{\small TR#1}}
\newcommand{\Red}[1]{{\color{black} #1}} % red
\newcommand{\jbb}[1]{{\color{black} #1}} % blue
\newcommand{\jbbo}[1]{{\color{black} #1}} % violet
\newcommand{\jbk}[1]{{\color{black} #1}} % black
\newcommand{\rfm}[1]{{\color{black} #1}} % cyan
\newcommand{\jbr}[1]{{\color{black} #1}} % blue (for revision)

\definecolor{forestgreen}{rgb}{0.13, 0.55, 0.13}
\newcommand{\jbrt}[1]{{\color{black} #1}} %forestgreen, green (for revision two)

\renewcommand{\b}[1]{\ensuremath{#1}} % bold
\newcommand{\bh}[1]{\ensuremath{\hat{#1}}}

\newcommand{\bk}[1]{\ensuremath{{#1}_k}}
\newcommand{\bko}[1]{\ensuremath{{#1}_{k+1}}}

\newcommand{\bik}[2][1]{\ensuremath{{#2}^{-#1}_k}}

\newcommand{\bs}[1]{\ensuremath{#1}}

\newcommand{\biidx}[2]{\ensuremath{#2}_{#1}}
\newcommand{\biidxh}[2]{\widehat{\ensuremath{#2}}_{#1}}

\newcommand{\tp}{^{\top}}   % 09/03/20
\newcommand{\mtp}{^{-\top}} % 11/19/20

\newcommand{\bks}[1]{\ensuremath{{#1}_k}(\sigma)} % 10/07/20

\newsiamremark{remark}{Remark}
\newsiamremark{hypothesis}{Hypothesis}
\crefname{hypothesis}{Hypothesis}{Hypotheses}
\newsiamthm{claim}{Claim}

\headers{RCR: Reduced Compact Representation}
{J. J. Brust, R. F. Marcia, C. G. Petra, and M. A. Saunders}

\title{Large-scale optimization with linear equality
       constraints using reduced compact representation%
\thanks{Dedicated to Dr Oleg Burdakov, 1953--2021.
   Version of \today.  Submitted to SISC 2021.%
   \funding{This work was supported by
   the U.S. Department of Energy, Office of Science,
   Advanced Scientific Computing Research,
   under Contract DE-AC02-06CH11357
   at Argonne National Laboratory.
   This work performed under the auspices of the
   U.S. Department of Energy by
   Lawrence Livermore National Laboratory
   under Contract DE-AC52-07NA27344.}}}

\author{Johannes J. Brust%
   \thanks{Department of Mathematics, University of California San Diego, San Diego, CA (formerly Argonne National Laboratory) %Mathematics and Computer Science Division,
           (\email{jjbrust@ucsd.edu}).} %\url{http://www.imag.com/\string~ddoe/}
\and Roummel F. Marcia%
   \thanks{Department of Applied Mathematics,
           University of California Merced, Merced, CA
           (\email{rmarcia@ucmerced.edu}).} %\email{jesmith@fictional.edu}
\and Cosmin G. Petra%
   \thanks{%\footnotemark[3]
           Center for Applied Scientific Computing,
           Lawrence Livermore National Laboratory,
           Livermore, CA
           (\email{petra1@llnl.gov}).}
\and \hbox{Michael A. Saunders}%
   \thanks{Department of Management Science and Engineering,
           Stanford University, Stanford, CA
           (\email{saunders@stanford.edu}).}}

\usepackage{amsopn}

\ifpdf
\hypersetup{
  pdftitle={Large-Scale Optimization with Linear Equality Constraints},
  pdfauthor={J. J. Brust, R. F. Marcia, C. G. Petra and M. A. Saunders}
}
\fi

\begin{document}

\pagestyle{empty}
  
\vspace{1.75in}

\begin{centering}

ARGONNE NATIONAL LABORATORY

9700 South Cass Avenue

Argonne, Illinois  60439

\vspace{1.5in}

{\large \textbf{Large-scale optimization with linear equality
       constraints using reduced compact representation}}

\vspace{.5in}

\textbf{J. J. Brust, R. F. Marcia, C. G. Petra and M. A. Saunders}

\vspace{.5in}

%% Argonne Leadership Facility and
Mathematics and Computer Science Division

\vspace{.25in}

Preprint ANL/MCS-P9279-0120

\vspace{.5in}

%January 2021
August 2021

\end{centering}

\vspace{2.0in}

\bigskip

\par\noindent
\footnotetext [1]
{
This work was supported by the U.S. Department of Energy, Office of Science,
Advanced Scientific Computing Research, under Contract DE-AC02-06CH11357
at Argonne National Laboratory.
through the Project "Multifaceted Mathematics for Complex Energy Systems."
This work was also performed under the auspices of the U.S. Department of Energy by 
Lawrence Livermore National Laboratory under Contract DE-AC52-07NA27344. 
}

\newpage

\vspace*{\fill}
\begin{center}
\fbox{
\parbox{4in}{
The submitted manuscript has been created by UChicago Argonne, LLC, Operator of Argonne 
National Laboratory (``Argonne''). Argonne, a U.S. Department of Energy Office of Science 
laboratory, is operated under Contract No. DE-AC02-06CH11357. The U.S. Government retains 
for itself, and others acting on its behalf, a paid-up nonexclusive, irrevocable worldwide 
license in said article to reproduce, prepare derivative works, distribute copies to the 
public, and perform publicly and display publicly, by or on behalf of the Government. 
The Department of Energy will provide public access to these results of federally 
sponsored research in accordance with the DOE Public Access 
Plan. \texttt{http://energy.gov/downloads/doe-public-accessplan}
}}
\end{center}
\vfill

\newpage
\pagestyle{plain}
\setcounter{page}{1}

\maketitle

\vspace{-3.5cm}

 \hbox to \textwidth{LLNL Release Number: LLNL-JRNL-818401\hss}
% or 
%\hbox to \textwidth{\hss LLNL Release Number: LLNL-JRNL-xxxxxx}

\vspace{3.0cm}

% REQUIRED
\begin{abstract}
  For %\jbb{minimization}
optimization problems with linear equality constraints,
we \jbr{prove} that the (1,1) block of the inverse KKT matrix remains unchanged
when projected onto the nullspace of the \jbrt{constraint matrix}.
We develop \emph{reduced compact representations} of the limited-memory \jbr{inverse} BFGS Hessian %approximation
to compute search directions efficiently \jbr{when the constraint Jacobian is sparse}.  Orthogonal projections are implemented by \jbrt{a} sparse QR factorization or \jbrt{a} preconditioned LSQR iteration.
In numerical experiments two proposed trust-region algorithms improve in computation times, often significantly, \Red{compared} to previous implementations \jbrt{of related algorithms} and compared to IPOPT.
% In numerical experiments the proposed algorithms improve in computation times, often significantly, \Red{compared} to previous implementations and compared to IPOPT.
\end{abstract}

% REQUIRED
\begin{keywords}
  Large-scale optimization, compact representation, trust-region method, limited memory, LSQR, sparse QR
\end{keywords}

% REQUIRED
\begin{AMS}
  68Q25, 68R10, 68U05
\end{AMS}

%\section{Background}
\jbbo{\section{Introduction}}
%\subsection{Problem formulation}
\label{sec:problem}
Linear equality constrained minimization problems are formulated as % optimization
\begin{equation}
    \label{eq:main}
    \underset{\b{x} \in \mathbb{R}^n}{ \text{ minimize } } f(\b{x}) 
    \quad \text{subject to} \quad \b{A} \b{x} = \b{b},
\end{equation}
where $ f:\mathbb{R}^n \to \mathbb{R} $ and $ \b{A} \in \mathbb{R}^{m \times n} $.
We assume that the number of variables $n$ is large,
$g(x) = \nabla f(x)$ is available,
$ A $ is sparse, % with full row rank (for the purpose of presentation),
and that the initial guess $ x_0 $ is feasible: $Ax_0=b$.
\jbr{If $ \b{A} $ has low rank, one can obtain a full-rank matrix
%$ \underbar{A} $, say,
by deleting rows in $ \b{A} $ that correspond to small
diagonals of the triangular matrix in a sparse QR factorization of $ \b{A}\tp $.
%The methods in this article
Our methods here
use the rank information contained in sparse QR factors, and
thus we assume %without loss of generality
that $A$ has full rank until 
implementation details are described in
section~\ref{sec:appB}.} %()
%(In our implementations, $ \b{A} $
%does not have to have full rank.) %, as will be apparent later.)}
\jbbo{For large problems,} computing the Hessian $ \nabla^2 f(\b{x}) \in \mathbb{R}^{n \times n}  $ is \jbbo{often} not practical, \jbbo{and} we approximate this matrix using
a limited-memory BFGS (Broyden-Fletcher-Goldfarb-Shanno, \cite{Broyden70,Fletcher70,Goldfarb70,Shanno70}) quasi-Newton matrix 
$ \bk{B} \approx \nabla^2 f(\bk{x}) $. Starting from $x_0$, %the feasible initial guess,
we update iterates according to $ \bko{x} = \bk{x} + \bk{s} $. The
step $ \bk{s} $ is computed as the solution of a quadratic trust-region subproblem, in which
the quadratic objective is defined as
$ q(s) \equiv \b{s}\tp \bk{g} + \frac{1}{2} \b{s} \tp \bk{B} \b{s} $ with $\bk{g} \equiv g(\bk{x})$. For a given trust-region radius $\Delta > 0$ and norm $\norm{\cdot}$, the
trust-region subproblem is
\begin{equation}
    \label{eq:trsub}
    \underset{ \| \b{s} \| \le \Delta }{ \text{ minimize } } q(\b{s}) 
    \quad \text{subject to} \quad \b{A} \b{s} = \b{0},
\end{equation}
which ensures that each search direction is in the nullspace of \jbrt{$A$}, and thus each iterate $x_k$ \jbrt{is} feasible.
%\jbb{while at the same time rapid progress towards the solution is made.}
%We refer to \cite[Section 2.2]{BMP19} for more details on 
%these trust-region subproblems.

% optimization 
\subsection{Background}
\label{sec:bkg}
\jbbo{Large problems of the form \eqref{eq:main} are the focus of recent research
because large statistical- and machine-learning problems can be cast in this way. As such,
\eqref{eq:main} constitutes the backbone of the Alternating Direction Method of Multipliers (ADMM)
\cite{BoydEtAl11}, with applications to optimal exchange problems, consensus and sharing problems, support-vector machines, and more. Recent work \cite{FuZhangBoyd20} emphasizes methods that use
gradients %only
of $\jbrt{f}$
and suggest accelerations via quasi-Newton approximations.}
Quasi-Newton methods %are effective techniques that iteratively 
estimate Hessian matrices using low-rank updates at each iteration (typically \hbox{rank-1} or \hbox{rank-2}).
Starting from an initial matrix, the so-called \emph{compact representation} of quasi-Newton matrices \cite{ByrNS94} is a matrix
representation of the recursive low-rank updates. Because the compact representation
enables effective limited memory implementations, which update a small number of previously stored 
vectors, these methods are well suited to large problems.
\jbr{Trust-region and line-search
methods are standard
%are by virtue the standard methods % the de facto In this context t
for determining search directions for smooth problems, 
%each of which has its own merits.
and each approach has its own merits.
Combinations of trust-region methods and quasi-Newton compact representations have 
been developed in \cite{BruBEM19,B18,BruEM15,BurdakovLMTR16}. Widely
used quasi-Newton line-search methods are \cite{Byrd2006,LiuNocedal89,ZhangHager04,ZhuByrdLuNocedal97}. The main ideas
in this article are applicable to both trust-region and line-search methods.}
%Current efforts for large-scale %\jbb{minimization}
% \Red{optimization} have combined quasi-Newton compact % unconstrained optimization 
% representations with trust-region methods, because trust-region methods have desirable
% convergence properties and are regarded as being more robust than alternative
% linesearch methods \cite{BruBEM19,B18,BruEM15,BurdakovLMTR16}.

\subsection{Compact representation}
\label{sec:compact}
A storage-efficient approach to quasi-Newton matrices is
the compact representation of Byrd et al.\ \cite{ByrNS94},
which represents the BFGS matrices in the form
\begin{equation}
    \label{eq:Bk}
   \bk{B} = \gamma_k \b{I} + \bk{J} \bk{M} \bk{J} \tp,
\end{equation}
with scalar $\gamma_k > 0$.
The history of vectors $\jbrt{\{}\bk{s}\jbrt{\}} = \jbrt{\{}\bko{x} - \bk{x}\jbrt{\}}$ and $\jbrt{\{}\bk{y}\jbrt{\}} = \jbrt{\{}\bko{g} - \bk{g}\jbrt{\}}$
is stored in rectangular
    $\bk{S} \equiv \bmat{\b{s}_0, \dots, \b{s}_{k-1}} \in \mathbb{R}^{n \times k}$
and $\bk{Y} \equiv \bmat{\b{y}_0, \dots, \b{y}_{k-1}} \in \mathbb{R}^{n \times k}$.
\jbrt{The} matrices
\begin{align}
             \bk{J} &\equiv \bmat{\bk{S} & \bk{Y} },
\\ \bk{S}\tp \bk{Y} &\equiv \bk{L} + \bk{D} + \bar{\b{T}}_k,
\\           \bk{M} &\equiv
                    - 
                          \bmat{
                             \delta_k \bk{S}\tp \bk{S} & \,\delta_k \bk{L}
                          \\ \delta_k \bk{L}\tp & \,- \bk{D}}^{-1}
                    %  - \bigg[
                    %       \begin{smallmatrix}
                    %          \delta_k \bk{S}\tp \bk{S} & \,\delta_k \bk{L}
                    %       \\ \delta_k \bk{L}\tp & \,- \bk{D}
                    %       \end{smallmatrix}
                    %   \bigg]^{-1}
\end{align}
are defined with $\delta_k = 1/\gamma_k$,
where $\bk{L}$ \jbr{and} $\bar{\b{T}}_k$ are the strictly lower and upper triangular parts
of $\bk{S}\tp \bk{Y}$ and $\bk{D}$ is the diagonal. For large problems, limited-memory versions
store only a small subset of recent pairs $\{ \b{s}_i, \b{y}_i  \}_{i=k-l}^{k-1}  $,
resulting in storage-efficient matrices $ \bk{J} \in \mathbb{R}^{n \times 2l} $
and $ \bk{M} \in \mathbb{R}^{2l \times 2l} $ \jbr{where $ l \ll n $}. % l \in [3,7]
\jbk{%For the L-BFGS matrix (cf.~
Following Byrd et al.\ 
\cite[Theorem 2.2]{ByrNS94},
%$ \bik{B} $ has the form 
the inverse BFGS matrix has the form
\begin{equation}
    \label{eq:iBk}
   %\bik{B} = \delta_k I + \bhk{\Psi} \bhk{M} \bhk{\Psi} % previous notation
   \bik{B} = \delta_k I + \bk{J} \bk{W} \bk{J} \tp,
\end{equation}
where %in this article, $\bk{J}$ and
$ \bk{W} \in \mathbb{R}^{2l \times 2l}  $ is  given by 
%$ \bk{J} = \bmat{ \bk{S} & \bk{Y} }  $ and
\begin{equation}
    \label{eq:Wk}
    \bk{W} =
        \bmat{
             \bk{T}^{-\top} (\bk{D} + \delta_k \bk{Y} \tp \bk{Y}) \bk{T}^{-1} & -\delta_k \bk{T}^{-\top}  \\
            -\delta_k\bik{T} & 0_{l \times l}}.
        % \big[
        % \begin{smallmatrix}
        %      \bk{T}^{-\top} (\bk{D} + \delta_k \bk{Y} \tp \bk{Y}) \bk{T}^{-1} & -\delta_k \bk{T}^{-\top}  \\
        %     -\delta_k\bik{T} & 0_{l \times l}
        % \end{smallmatrix}
        % \big].
\end{equation}
The diagonal matrix $\bk{D}$ (and hence the upper triangular matrix $\bk{T} \equiv \bk{D} + \bar{\b{T}}_k$) are nonsingular as long as $\bk{B}$ is also.}

\subsection{Outline}
\label{sec:outline}
 \jbbo{Section~2 describes our contributions in the context of large problems, while section~3 motivates our proposed representations.
Section~4 develops the reduced compact representation and updating techniques that enable efficient implementations.
Section~5 describes computations of orthogonal
projections, and the trust-region strategy for optimization.
Section~6 gives an efficient method
when an $\ell_2$-norm trust-region subproblem is used.
Sections~7 and~8 develop an effective factorization, and a method that uses a shape-changing norm in the trust-region subproblem.
Numerical experiments are reported in section~9,
and conclusions are drawn in section~10.}

%\section{Article contribution}
\section{Contributions}
\label{sec:contrib}

%The KKT system for \eqref{eq:trsub} without the norm constraint is
\jbr{The first-order necessary conditions for the solution of problem \eqref{eq:trsub}
without the norm constraint
are characterized by the linear system} %subsumed
\begin{equation}\label{eq:KKTsystem}
   \bmat{\b{B}_k & \b{A}\tp
      \\ \b{A}   & \b{0}_{m \times m}}
   \bmat{\jbr{\b{s}_E} \\ \jbr{\bs{\lambda}_E}} % \bs{s}_e, \bs{\lambda}_e
   =
   \bmat{-g_k \\ \b{0}_m},
\end{equation}
where $ \jbr{\bs{\lambda}_E} \in \mathbb{R}^m $ is a vector of Lagrange multipliers 
\jbr{and $ \b{s}_E $ denotes the ``equality'' constrained minimizer of \eqref{eq:trsub}.
Adopting the naming convention of \cite[Sec.~16.1, p.~451]{NocW06}, we refer to
\eqref{eq:KKTsystem} as the KKT system
%We would like to note that this convention is
(a slight misnomer, as use of the system for the equality constrained setting
predates the work of Karush, Kuhn, and Tucker).} % ()
%\jbbo{and $\b{s}=\b{s}_e$}. 
For large $n$, compact representations
of the (1,1) block in the inverse KKT matrix were recently proposed by Brust et al.\ \cite{BMP19}.
% \rfm{[]Johannes: Are these LTRL2-LEC and LTRSC-LEC in \cite{BMP19}?]}
% In BMP19 we ended up calling the algorithms TR1 & TR2
Two limited-memory trust-region algorithms, \rfm{{\small LTRL2-LEC} and {\small LTRSC-LEC} (which we refer to as \TR1 and \TR2 in the numerical experiments in Sec.~\ref{sec:numex})},
use these representations to compute search directions efficiently when $\b{A}$ has relatively few rows.
This article develops efficient algorithms when the number of equality constraints is large and the constraint matrix is sparse. In particular,
by exploiting the \jbr{property} that part of the solution to the KKT system 
\jbbo{is unaltered} when it is projected onto the nullspace of \jbrt{$A$}, 
we develop \emph{reduced compact representations (RCR)}, \jbbo{which need
a small amount of memory and} lead to efficient
methods for \jbrt{solving} problems with many constraints 
\Red{(large $m$ and $n$) and possibly many degrees of freedom (large $n-m$)}.
In numerical experiments \jbrt{when solving} large problems, the proposed methods 
%use less computation time often with significant reductions.}
are often significantly more efficient than both our previous implementations and IPOPT \cite{WaechterBiegler06}.

% extends {\small TR1} and {\small TR2} and 

% \begin{center}
%     \small Statement of problem needed somewhere
% \end{center}

\section{Motivation}
\label{sec:motiv} % equality  constrained and the quasi-Newton matrix $ \bk{B} $ 
%The KKT system, without the norm constraint in \eqref{eq:trsub}, is

The solution $\jbr{\b{s}_E}$ in \eqref{eq:KKTsystem} can be computed from only the (1,1)
block of the inverse KKT matrix, as opposed to both the (1,1) and (1,2) blocks,
because of the zeros in the right-hand side. Let $ \bk{V}$ be the (1,1) block
of the inverse KKT matrix (obtained for example
from a block LDU factorization).  It is given by % , when denoted as $ \bk{V}$, this matrix 
\begin{equation}
    \label{eq:invKKT11}
    \bk{V} \equiv (\bik{B} - \bik{B}\b{A}\tp(\b{A} \bik{B} \b{A}\tp)^{-1}\b{A}\bik{B}),
\end{equation}
%Once $ \bk{V} $ is known, part of the solution in \eqref{eq:KKTsystem} is computed as
and then $ \jbr{\b{s}_E} = -\bk{V}\bk{g} $.
At first sight the expression in \eqref{eq:invKKT11} 
\jbr{appears to be expensive to compute because of the multiple inverse operations
and matrix-vector products.}
%appears involved for practical computation. 
However, as % the purpose of
%$ \bik{B} $ has the form 
$
%\bik{B} = \delta_k I + \bhk{\Psi} \bhk{M} \bhk{\Psi} % previous notation
\bik{B} = \delta_k I + \bk{J} \bk{W} \bk{J} \tp
$,
% where $ \bk{J} \in \mathbb{R}^{n \times 2l} $ and $\bk{W} \in \mathbb{R}^{2l \times 2l}$
% % $ \bk{J} = \bmat{ \bk{S} & \bk{Y} } $
% , (\cite[Theorem 2.2.]{ByrNS94})
we can exploit computationally useful structures. Specifically, with 
$ 
    %\bk{\Omega} \equiv (\b{A} \bik{B} \b{A} \tp)^{-1} % previous notation
    \bk{G} \equiv (\b{A} \bik{B} \b{A} \tp)^{-1}
$
and 
$ 
    \bk{C} \equiv \b{A} \bk{J} \bk{W},
$
\cite[Lemma 1]{BMP19} describes the expression
\begin{equation}
    \label{eq:cmpVk_L1}
    \bk{V} = \delta_k I +
        %\bmat{\b{A}\tp & \bk{S} & \bk{Y} }
        \bmat{\b{A}\tp & \bk{J} }
        \bmat{  -\delta_k^2 \bk{G} & - \delta_k \bk{G} \bk{C} \\
                -\delta_k \bk{C} \tp \bk{G} & \bk{W} - \bk{C} \tp \bk{G} \bk{C}}
        \bmat{ \b{A} \\
                \bk{J}\tp}.        
        % \bmat{ \b{A} \\
        %         \bk{S}\tp \\
        %         \bk{Y} \tp}
\end{equation}
For large $n$, once the components of the middle matrix in \eqref{eq:cmpVk_L1} are available, this compact 
representation of $ \bk{V} $ enables efficient computation
of a matrix-vector product $ \bk{V} \bk{g} $, \jbbo{hence the solution of \eqref{eq:KKTsystem}}, and an economical eigendecomposition $ \bk{V} = \b{U} \b{\Lambda} \b{U}\tp $. % implicit
However, unless $m$ is small (there are few rows in $ \b{A} $), \jbrt{multiplying}
with the $ (m+2l) \times (m+2l) $ middle matrix is not practical.

With large $n$ and $m$ in mind, we
\jbr{note} that the solution $ s_E $ is unchanged if instead of
$\bk{g}$ a projection of this vector onto the nullspace of $\b{A}$ is
used, or if $\jbr{\b{s}_E}$ is projected onto the nullspace of $ \b{A} $. 
This is a consequence of the properties of $ \bk{V} $.
To % s_E
formalize these statements, let the orthogonal projection matrix onto $ \text{null}(\b{A}) $ be 
$ 
	\b{P} =  \b{I}_n - \b{A}\tp (\b{A} \b{A}\tp)^{-1}\b{A}.
$
\jbr{Since the} \jbrt{columns of the} (1,1) \jbr{block of the inverse from} \eqref{eq:KKTsystem}
(namely \jbrt{columns of} $V_k$) \jbr{are in the nullspace of $A$, the orthogonal projection onto $\text{null}(A) $ acts as an identity operator on the vector space spanned by $V_k$:} 
%We observe, using $ \b{P} $ and \eqref{eq:invKKT11}, that
\begin{equation}
\label{eq:projKKT1}	 	
	\bk{V}	= \bk{V}\b{P} % \b{s}_e
	        = \b{P}\tp\bk{V}
			= \b{P}\tp \bk{V} \b{P}.
			%= - \b{P}\tp \bk{V} \b{P}\bk{g}.
\end{equation}
% \label{eq:projKKT1}	 	
% 	\b{s}_e	= - \bk{V}\bk{g} % \b{s}_e
% 			= - \bk{V} \b{P}\bk{g}
% 			= - \b{P}\tp \bk{V} \b{P}\bk{g}.
% \end{equation}
\jbr{Relation \eqref{eq:projKKT1} can equivalently be derived from 
\eqref{eq:invKKT11}, the expression for $P$, and the equality $ V_k A \tp = 0 $.} % ()
The methods in this article are based on representations of projected matrices $ \b{P}\tp \bk{V} \b{P}$ $ \in \mathbb{R}^{n \times n}$, % \in \mathbb{R}^{n \times n} $
whose properties enable desirable numerical advantages for large $ n $ and $m$.
Instead of \jbrt{multiplying}
with the possibly \jbrt{large} $ \bk{G} \in \mathbb{R}^{m \times m} $ and 
$ \bk{C} \in \mathbb{R}^{m \times 2l}  $ in \eqref{eq:cmpVk_L1},
we store the matrices $ \bk{S} \in \mathbb{R}^{n\times l}$ and 
$\bk{Z} \equiv \b{P}\bk{Y} \in \mathbb{R}^{n\times l} $ and 
\jbr{small square matrices that depend on the memory parameter $l$ but
not on $m$. The columns of \bk{Z} are defined as
$
	\bk{z} = \b{P}\bk{y} = \b{P}(\bko{g}-\bk{g}),
$ and they are contained in the nullspace of $A$.}

%a small square matrix $ \bk{N} \in \mathbb{R}^{2l \times 2l} $, which we specify later. 
% The columns of \bk{Z} are defined as
% $
% 	\bk{z} = \b{P}\bk{y} = \b{P}(\bko{g}-\bk{g}).
% $
% Because $\b{A}\bk{S} = 0$, we have $ \b{P}\bk{S} = \bk{S} $, and since $ \b{P} \b{A}\tp = \b{0}  $ we therefore seek a ``reduced" compact representation of \eqref{eq:cmpVk_L1} in the form
% \begin{align*}
%   \b{P}\tp \bk{V} \b{P} &= \delta_k \b{P} + 
%   \bmat{\bk{S} & \bk{Z}} \bk{N} \bmat{\bk{S} & \bk{Z}} \tp,
% \end{align*}
% % \begin{align*}
% %   \b{P}\tp \bk{V} \b{P} &= \delta_k \b{P} + 
% %   \bmat{\bk{S} & \bk{Z}} \bk{N} \bmat{\bk{S} & \bk{Z}} \tp
% % \\    &= \delta_k \b{I}_n +
% %          \bmat{\b{A}\tp & \bk{S} & \bk{Z}}
% %          \bmat{-\delta_k (\b{A} \b{A}\tp)^{-1}
% %             \\      & \bk{N}}
% %          \bmat{\b{A} \\ \bk{S}\tp \\ \bk{Z}\tp},
% % \end{align*}
% where $\delta_k > 0$ is a scalar and $ \bk{N} $ is small.} % n initialization
%Thus far

\jbr{With \eqref{eq:invKKT11} and \eqref{eq:cmpVk_L1}} we motivated the solution of \eqref{eq:trsub} without the norm
constraint (giving the equality-constrained step $s_E$). Computing $s_E$ is important for the implementation of practical algorithms, but it is even more important to solve \eqref{eq:trsub} efficiently with the norm constraint. In Sec.~\ref{sec:l2}, using the $ \ell_2 $ norm, we develop a modified
version of $ \bk{V} $ as a function of a scalar parameter $\sigma>0$, i.e., $ \bks{V} $. In Secs.~\ref{sec:eigen} and~\ref{sec:SCNorm}, %and \ref{sec:SCNorm}, 
we describe how the structure of 
$ \bk{V} $ can be exploited to compute an inexpensive eigendecomposition
that, when combined with a judiciously chosen norm \jbr{(the shape-changing infinity norm
from \cite[Sec. 4.2.1]{BurdakovLMTR16})}, provides a search direction by
an analytic formula.
\jbk{Note that the representation of $ \bk{V} $ is not specific to the L-BFGS matrix,
and other compact quasi-Newton matrices could be used
(Byrd et al.\ \cite{ByrNS94}, DeGuchy et al.\ \cite{DeGuchyEM16}).
}

\section{Reduced compact representation (RCR)}

\jbr{This section describes a computationally effective representation of \eqref{eq:projKKT1}, 
which we call the \emph{reduced compact representation} (RCR). In section~\ref{sec:redHessian},
the RCR is placed into historical context with reduced Hessian methods.  Subsequently, sections \ref{sec:redComp}--\ref{sec:compComplx} develop the specific formulas that enable effective computations.}

\subsection{Reduced Hessian}
\label{sec:redHessian}
\jbk{The name \emph{reduced compact representation} is related to the term
\emph{reduced Hessian} \cite{GillMurray74},
% we bring into context with the relatively known term
%Forsgren and Murray \cite{FM93} or Murtagh and Saunders \cite{MS78}.
%Historically,
where $ \jbr{\bh{Z}} \in \mathbb{R}^{n \times (n-m)} $ denotes a basis for the
nullspace of $\b{A}$ (satisfying $A\jbr{\bh{Z}}=0$).  In turn, $\jbr{\bh{Z}}$ defines the so-called reduced Hessian matrix
as $ \jbr{\bh{Z}}\tp \nabla^2 f_k \jbr{\bh{Z}} $ or $ \jbr{\bh{Z}}\tp \bk{B} \jbr{\bh{Z}} $. % with a 
%quasi-Newton matrix
%Since $ \jbr{\bh{Z}}\tp \bk{B} \jbr{\bh{Z}} \in \mathbb{R}^{(n-m)\times(n-m)} $
%is smaller than $ \bk{B} \in \mathbb{R}^{n \times n} $, it is adequately 
%referred to as a \emph{reduced Hessian}. 
In order to compute an equality-constrained step
$ \jbr{\b{s}_E} $, a reduced Hessian method solves $ (\jbr{\bh{Z}}\tp \bk{B} \jbr{\bh{Z}})\jbr{\bh{s}_E} = -\jbr{\bh{Z}}\tp \bk{g} $ and computes $ \jbr{\b{s}_E} = \jbr{\bh{Z}}\jbr{\bh{s}_E} $. Known computational challenges
with reduced Hessian methods are that a desirable basis $ \jbr{\bh{Z}} $ may be expensive to 
compute, the condition number of the reduced linear system may be larger than the original one, and 
the product $ \jbr{\bh{Z}}\tp \bk{B} \jbr{\bh{Z}} $ is not necessarily sparse even if the 
matrices themselves are. For large-scale problems, these challenges can result in 
significant computational bottlenecks. In the sequel we refer to $ \b{P}\tp \bk{V} \b{P} $
as a \emph{reduced compact representation} because it has a reduced memory footprint compared to
$ \bk{V} $ in \eqref{eq:cmpVk_L1} 
(although the matrices have the same dimensions). We also note that $ \bk{V} $
and $ \b{P}\tp \bk{V} \b{P} $ have the same condition}, \jbbo{and $ \b{P}\tp \bk{V} \b{P} $ 
has structure that enables efficient implementations.}
%\jbb{\subsection{Reduced Compact Representation}}
% represents an
% orthogonal projection, which is typically less expensive computationally than
% forming an orthonormal basis.

\subsection{Reduced compact representation}
\label{sec:redComp}
To simplify \eqref{eq:cmpVk_L1}, we note
that $ \bk{V} = \b{P}\tp \bk{V} \b{P} $, \jbr{that} $ \b{P}\tp\! \b{A} \tp = 0 $,
and $ \b{P} \tp\! \bk{J} = \bmat{\bk{S} & \bk{Z}} $ \jbr{(where $\b{P}\tp \bk{Y} \equiv \bk{Z}$ by definition)}, so that
% \bk{V}  = 
\begin{equation*}
    \b{P}\tp \bk{V} \b{P} 
            = \delta_k \b{P} +
            \bmat{\bk{S} & \bk{Z} }
            ( \bk{W} - \bk{C} \tp \bk{G} \bk{C} ) \bmat{\bk{S} & \bk{Z} }\tp\!.
\end{equation*}
In Appendix A we show that $ \bk{C} \tp \bk{G} \bk{C} $
simplifies to 
$  
    \bk{C} \tp \bk{G} \bk{C} =
    \big[
        \begin{smallmatrix}
            (\bk{C} \tp \bk{G} \bk{C})_{11}                     & 0 \\
            0                                                   & 0
        \end{smallmatrix}
    \big]
$ with $(\bk{C} \tp \bk{G} \bk{C})_{11} = \delta_k  \bk{T}^{-\top} \bk{Y} \tp \b{A} \tp (\b{A} \b{A} \tp)^{-1} \b{A} \bk{Y} \bik{T}  $. Based on this, we derive a \emph{reduced compact representation} of $ \bk{V} $.

\smallskip

{Lemma 1:} The \emph{RCR of $\bk{V}$
in \eqref{eq:cmpVk_L1} for the L-BFGS matrix is given by}
\begin{equation}
    \label{eq:redComp}
    \bk{V} = \delta_k I + 
    \bmat{\b{A}\tp & \bk{S} & \bk{Z} }
    \bmat{ - \delta_k (\b{A} \b{A} \tp)^{-1} &  \\
            & \bk{N}
    }
    \bmat{\b{A} \\ \bk{S} \tp \\ \bk{\jbr{Z}} \tp },
\end{equation}
\emph{where}
\begin{equation*}
    \bk{N} = \bmat{ \bk{T}^{-\top} (\bk{D} + \delta_k \bk{Z} \tp \bk{Z}) \bk{T}^{-1} & -\delta_k \bk{T}^{-\top}  \\
    -\delta_k\bik{T} & 0_{k \times k} }.
\end{equation*}

\smallskip

%{Proof.} 
\begin{proof} Multiplying $ \bk{V} $ in \eqref{eq:cmpVk_L1} from the left and right 
by $ \b{P} \tp $ and $ \b{P} $ yields
$ \bk{V} = \delta_k \b{P} + \bmat{\bk{S} & \bk{Z} } ( \bk{W} - \bk{C} \tp \bk{G} \bk{C} ) \bmat{\bk{S} & \bk{Z} }\tp $. Since only the (1,1) block in $ \bk{C} \tp \bk{G} \bk{C} $
is nonzero, we consider only the (1,1) blocks, namely
\begin{equation*}
 (\bk{W})_{11} - (\bk{C}\tp \bk{G} \bk{C})_{11} =
 \bk{T}^{-\top} (\bk{D} + \delta_k( \bk{Y} \tp \bk{Y} - 
        \bk{Y} \tp \b{A} \tp (\b{A} \b{A} \tp)^{-1} \b{A} \bk{Y} )) \bk{T}^{-1}.
\end{equation*}
\jbr{Since} 
$ \bk{Y} \tp \b{P}\tp \bk{Y} = \bk{Y} \tp \b{P}\tp \b{P} \bk{Y} = \bk{Z} \tp \bk{Z} $, we
obtain the (1,1) block in $ \bk{N} $. Subsequently, by factoring 
\jbr{$P$ as
\begin{equation*}
    P = \b{I} - \bmat{\b{A}\tp & \bk{S} & \bk{Z} }
    \bmat{ - \delta_k (\b{A} \b{A} \tp)^{-1} &  \\
            & \b{0}_{2k \times 2k}
    }
    \bmat{\b{A} \\ \bk{S} \tp \\ \bk{\jbr{Z}} \tp },
    %\bmat{\b{A}\tp & \bk{S} & \bk{Z} }\tp,
\end{equation*}
we see that 
\begin{equation*}
    \b{P}\tp \bk{V} \b{P} 
            = \delta_k \b{I} + \bmat{\b{A}\tp & \bk{S} & \bk{Z} }
            \bmat{ - \delta_k (\b{A} \b{A} \tp)^{-1} &  \\
            & \bk{W} - \bk{C}\tp \bk{G} \bk{C}
    }
    \bmat{\b{A} \\ \bk{S} \tp \\ \bk{\jbr{Z}} \tp }.
\end{equation*}
Because all blocks of $ \bk{W} - \bk{C}\tp \bk{G} \bk{C} $
except for the (1,1) block are equal to those in $\bk{W}$,
all blocks in $ \bk{N} $ are fully specified
%, completing 
and representation \eqref{eq:redComp}
is complete.}
% $ P = I - \b{A} \tp (\b{A} \b{A} \tp )^{-1} \b{A} $ we deduce
% \eqref{eq:redComp}. 
\end{proof} 
%\qed

\smallskip

Note that 
$ \bk{S} \tp \bk{Y} = \bk{D} + \bk{L} + \bar{\b{T}}_k = \bk{S} \tp \bk{Z} $,
which means \jbr{that} $ \bk{D} $ and $ \bk{T} = \bk{D} + \bar{\b{T}}_k $ can be
computed from $ \bk{S} $ and $ \bk{Z} $ alone, and \jbr{that} \bk{G} \jbr{and} \bk{C} \jbrt{need not be explicitly computed}. Therefore, for the RCR, 
only $ \bk{S} $, $ \bk{Z} $, $ \bk{T} $ \jbr{and} $ \bk{D} $ are stored.
% only $ \bk{S} $ and $ \bk{Z} $ are stored
% instead of $ \bk{S} $ and $ \bk{Y} $. 
An \jbbo{addition} is the scalar $ \delta_k $, % exception 
which is typically set to be 
$ \delta_k = \bk{s} \tp \bk{y} \big / \bk{y}\tp \bk{y} = \bk{s} \tp \bk{z} \big / \bk{y}\tp \bk{y} $ and may depend on the most recent $ \bk{y} $.
\jbbo{As $ \b{P} \bk{J} = \bmat{ \bk{S} & \bk{Z} }  $,
we \jbr{note} a key advantage of the RCR: that
\eqref{eq:redComp} can be written as 
\begin{equation}
    \label{eq:redCompConcise}
    \bk{V} = \delta_k \b{P} +  \b{P} \bk{J} \bk{N} \bk{J}\tp \b{P} \tp =
    \delta_k \b{P} +  \bmat{ \bk{S} & \bk{Z} } \bk{N} 
    \bmat{\bk{S} \tp \\[4pt] \bk{Z}\tp}.
    %\bk{V} = \delta_k \b{P} + \bhk{J} \bk{N} \bhk{J}\tp.
\end{equation}
By storing a few columns of $ \bmat{ \bk{S} & \bk{Z} } \in \mathbb{R}^{n \times 2l} $
(as described in section \ref{sec:compComplx}), which in turn define
a small matrix $ \bk{N} \in \mathbb{R}^{2l \times 2l} $ (cf.\ Lemma 1), we can separate the solves with $ \b{A} \b{A}\tp $ from other
calculations. Concretely, \jbr{note} that solves with $ \b{A} \b{A} \tp $
only occur \jbrt{as part of} the orthogonal projection $ \b{P} $, which can be represented
as a linear operator and does not need to be explicitly formed. Also \jbr{note}
that \eqref{eq:iBk} and \eqref{eq:redCompConcise} are related, with the
difference being that $ \bk{Y} $ and $ \delta_k \b{I} $ in \eqref{eq:iBk}
are replaced by $ \bk{Z} $ and $ \delta_k \b{P} $ in \eqref{eq:redCompConcise}.
Hence for large $ n $ and $ m $, computation with 
\eqref{eq:redCompConcise} is efficient and requires little memory,
provided orthogonal projections with $ \b{P} $ are handled effectively (as described in section \ref{sec:Proj}).
On the other hand, the compact representation in \eqref{eq:cmpVk_L1}
does not neatly decouple solves with $ \b{A} \b{A} \tp $, and results
in perhaps prohibitively expensive computations for large $m$. In particular,
$ \bk{G} $ in the middle matrix of \eqref{eq:cmpVk_L1} is defined
by $ \bk{G} \equiv ( \b{A} \bik{B} \b{A} \tp )^{-1} \in \mathbb{R}^{m \times m} $, which interleaves
solves with $ \b{A} \b{A} \tp $ and other terms. 
Therefore, the RCR in \eqref{eq:redComp}--\eqref{eq:redCompConcise} is recognizably more practical for large $ n $
and $m$ than \eqref{eq:cmpVk_L1}.}
We apply $ \bk{V} $ \jbbo{from \eqref{eq:redCompConcise}} to a vector $ \b{g} $ as
 \begin{equation}
    \label{eq:Vg}
     h = \bmat{ \bk{S} \tp \\[4pt] \bk{Z} \tp } \b{g}, \qquad
     \bk{V} \b{g} = 
     \bmat{ \bk{S} & \bk{Z} } \bk{N}h + \delta_k \b{P} \b{g}.
 \end{equation}

\subsection{Computational complexity}
\label{sec:compComplxNew}
\jbbo{With adequate \jbrt{precomputation and storage}, the cost of the matrix-vector product \eqref{eq:Vg} is often inexpensive. If the columns of $ \bk{Z} $ are stored, updating the small $ 2l \times 2l $ matrix $ \bk{N} $ does not depend on
solves with $ \b{A} \b{A} \tp $. Moreover, factors of $ \b{P} $ can be precomputed
once at $k=0$ and reused. In particular, suppose that a (sparse) QR factorization
$ 
\b{A}\tp = 
    \big[
        \begin{smallmatrix}
            \biidx{1}{Q} & \biidx{2}{Q}
        \end{smallmatrix}
    \big]
    \big[
        \begin{smallmatrix}
            \b{R} \\ 
            \b{0}
        \end{smallmatrix}
    \big]
%\b{A}\tp = \bmat{ \b{Q}_1 & \b{Q}_2 } \bmat{ \b{R}\\\b{0} } 
$ 
is obtained once, with $\b{Q} = \big[
        \begin{smallmatrix}
            \biidx{1}{Q} & \biidx{2}{Q}
        \end{smallmatrix}
    \big] $ \jbrt{being sparse}, such that the product $ \b{Q}\tp \!\b{g} $ takes
    $ \mathcal{O}(rn) $ multiplications, where $r$ is constant. Subsequently, the projection
$ \b{P}\b{g} = \b{g} - \biidx{1}{Q}\biidx{1}{Q}\tp \b{g} $ can be computed in
$ \mathcal{O}(n + 2rn) $ multiplications 
(or $ \b{P}\b{g} =  \biidx{2}{Q} \biidx{2}{Q}\tp \b{g} $ in 
$ \mathcal{O}(2rn) $ multiplications). Thus, we summarize
the multiplications in \eqref{eq:Vg} as: $ \b{h} $ with
$ 2nl $, $ \bk{N}\b{h} $ with negligible $(2l)^2$,
$ \bmat{\bk{S} & \bk{Z}} \bk{N}\b{h}  $ with $2nl$, and
$ \b{P} \b{g} $ with, say, $ 2nr $. The total, without negligible terms,
is $ \mathcal{O}( 2n(2l + r)) $. The multiplications scale linearly with $n$,
are related to the sparsity in $ \b{A} $, and are thus suited for large problems.}

\subsection{Updating}
\label{sec:compComplx}
We store and update the columns of $ \bk{Z} = \bmat{ \biidx{k-l}{z} & \cdots & \biidx{k-1}{z} } $
one at a time and \jbr{recall} that $ \biidx{k}{z} = \b{P} \biidx{k+1}{g} - \b{P} \bk{g} $.
Based on this, no additional solves with $ \b{A} \b{A} \tp $ are required to \jbrt{represent}
the matrix $ \biidx{k+1}{V} $. Specifically, suppose that we computed and stored
$ \b{P} \bk{g} $ at the end of the previous iteration, and that we compute 
$ \b{P} \biidx{k+1}{g} $ at the end of the current iteration. We can use this vector
in two places: first to \jbrt{represent} $ \biidx{k+1}{Z} $ with $ \bk{z} = \b{P} \biidx{k+1}{g} - \b{P} \bk{g}$ and hence $ \bko{V} $, and secondly in the 
computation of $ \bko{V} \bko{g} $. Thus only one solve with $ \b{A} \b{A} \tp $
per iteration is necessary to update $ \bko{V} $ and to compute a step of the
form $ \b{s} = - \bko{V} \bko{g} $.

For large problems, the limited-memory representation in \eqref{eq:redComp} is obtained
by storing only the last $l$ columns of $ \bk{S} $ and $\bk{Z} $. With $1 \le l \ll n$, limited-memory strategies enable computational efficiencies and lower storage requirements \cite{Noc80}. Updating $ \bk{S}$ and $ \bk{Z} $ requires replacing or inserting one column at each iteration. Let an underline below a matrix represent the 
matrix with its first column removed. That is, $ \underline{\b{Z}}_k $ represents $ \bk{Z} $ without its first column.
With this notation, a column update of a matrix $ \bk{Z} $ by a vector $ \bk{z} $ is defined as
\begin{equation*}
	\text{colUpdate}\left(\bk{Z},\bk{z} \right) \equiv
	\begin{cases}
		[\: \bk{Z} \: \bk{z}\:  ] 						& \text{ if } k < l, \\
		[\: \underline{\b{Z}}_k \: \bk{z}\:  ] 			& \text{ if } k \ge l. \\
	\end{cases}
\end{equation*}
Such a column update either directly appends a column to a matrix or first removes a column and then appends one. 
This column update will be used, for instance, to obtain $ \bko{Z} $ from $ \bk{Z} $ and $ \bk{z} $, i.e., $\bko{Z}= \text{colUpdate}( \bk{Z}, \bk{z} )$. Next, let an overline above a matrix represent
the matrix with its first row removed. That is, $ \overline{\b{S}\tp_k \b{Z}}_k $ represents $ \b{S}\tp_k \bk{Z} $ without its first row.
With this notation, a product update of $ \bk{S}\tp\bk{Z} $ by matrices $ \bk{S} $ \jbrt{and} \bk{Z} 
and vectors $ \bk{s} $ \jbrt{and} $\bk{z}$ is defined as
\begin{equation*}
	\text{prodUpdate} \left( \bk{S}\tp\bk{Z}, \bk{S}, \bk{Z}, \bk{s}, \bk{z} \right) \equiv 
	\begin{cases}
		\left[
			\begin{array}{ c c }
				\bk{S}\tp\bk{Z} 		& \bk{S}\tp\bk{z} \\
				\bk{s}\tp\bk{Z}	& \bk{s}\tp \bk{z} 
			\end{array}
		\right] & \text{ if } k < l, \vspace{0.1cm} \\		
		\left[
			\begin{array}{ c c }
				\left(\underline{\overline{\b{S}\tp_k \b{Z}_k}}\right) 			& 	\underline{\b{S}}_k\tp\bk{z} \\
				\bk{s}\tp \underline{\b{Z}}_k					&	 \bk{s}\tp \bk{z} 
			\end{array}
		\right] & \text{ if } k \ge l.\\
	\end{cases}
\end{equation*}
This product update is used to compute matrix products such as $ \bko{S}\tp \bko{Z} $ with 
$ \mathcal{O}(2ln) $ multiplications, instead of $ \mathcal{O}(l^2n) $ when the product 
$ \bk{S}\tp \bk{Z} $  is stored and the vectors $ \bk{s} $ \jbrt{and} $ \bk{z} $ have been computed. Note that a diagonal matrix can be updated in this way by setting the rectangular matrices $\bk{S} $ \jbrt{and} $ \bk{Z}$ to zero
and $ \bko{D} = \text{prodUpdate}(\bk{D},\b{0},\b{0},\bk{s},\bk{z})$. An upper triangular matrix
can be updated in a similar way, e.g., $ \bko{T} = \text{prodUpdate}(\bk{T},\bk{S},\b{0},\bk{s},\bk{z}) $. 
To save computation, products with zero matrices are never formed explicitly.
% This product update is used to compute matrix products, such as, $ \bko{Z}\tp \bko{Z} $, with 
% $ \mathcal{O}(2ln) $ multiplications, instead of $ \mathcal{O}(l^2n) $ when the product 
% $ \bk{S}\tp \bk{Z} $  had previously been stored and the vector $ \bk{z} $ been computed. Note that a diagonal matrix can be updated in this way by setting the rectangular matrices (e.g., $\bk{S}, \bk{Z}$) to zero, such that
% $ \bk{D} = \text{prodUpdate}(\bk{D},\b{0},\b{0},\bk{s},\bk{z})$. An upper triangular matrix
% can be updated in a similar way, e.g., $ \bk{T} = \text{prodUpdate}(\bk{T},\bk{S},\b{0},\bk{s},\bk{z}) $. 
% To save computations, products with zeros matrices are never formed explicitly.

% \Red{Should the preceding paragraph be as follows?}

\noindent
% This product update is used to compute matrix products such as $ \bko{Z}\tp \bko{Z} $ with 
% $ \mathcal{O}(2ln) $ multiplications, instead of $ \mathcal{O}(l^2n) $ when the product 
% $ \bk{S}\tp \bk{Z} $  is stored and the vector $ \bko{z} $ has been computed. Note that a diagonal matrix can be updated in this way by setting the rectangular matrices $\bk{S}, \bk{Z}$ to zero
% and $ \bko{D} = \text{prodUpdate}(\bk{D},\b{0},\b{0},\bk{s},\bk{z})$. An upper triangular matrix
% can be updated in a similar way, e.g., $ \bko{T} = \text{prodUpdate}(\bk{T},\bk{S},\b{0},\bk{s},\bk{z}) $. 
% To save computation, products with zero matrices are never formed explicitly.
% \Red{(?)}

% Limited memory discussion

% \begin{equation*}
%     \bk{W} = 
%         \bmat{ \bik{T}(\bk{D} + \delta_k \bk{Y} \tp \bk{Y}) \bk{T}^{-\top} & -\delta_k \bk{T}^{-\top} \\
%         -\bik{T} & 0_{l \times l}}
% \end{equation*}

% Previous notation
% \begin{equation}
%     \label{eq:cmpVk_L1}
%     \bk{V} = \delta_k I +
%         \bmat{\b{A}\tp & \bk{S} & \bk{Y} }
%         \bmat{  -\delta_k^2 \bk{\Omega} & - \delta_k \bk{\Omega} \bk{C} \\
%                 -\delta_k \bk{C} \tp \bk{\Omega} & \bk{W} - \bk{C} \tp \bk{\Omega} \bk{C}}.
% \end{equation}

\section{Computing projections}
% $z = Py$}{z=Py}}
\label{sec:Proj}

With $P = I_{\jbrt{n}} - A\tp (AA\tp)\inv A$, projections $z=Py$ can be computed by direct or
iterative methods.  Their efficiency depends on the sparsity of $A$.

\subsection{QR factorization}
\label{sec:QR}
When $A$ has full row-rank and the QR factorization
\begin{equation}
   \label{eq:QR}
   A\tp = Q\bmat{R \\ 0} = \bmat{Q_1 & Q_2} \bmat{R \\ 0} = Q_1 R
\end{equation}
is available, the projection operator becomes
$P = I - Q_1 Q_1\tp = Q_2 Q_2\tp$.
Thus, $z = Py$ can be computed stably as $z = Q_2(Q_2\tp y)$.
With $m<n$, the QR factors are best obtained using a product of Householder
transformations \cite{GV4}:
\begin{equation}
   \label{eq:HQR}
   Q\tp A\tp = H_m \dots H_3 H_2 H_1 A\tp = \bmat{R \\ 0} = \bmat{Q_1\tp \\ Q_2\tp} A\tp.
\end{equation}
Thus $Q = H_1 H_2 H_3 \dots H_m$ and the operators $Q_1$ and $Q_2$ are available from
\begin{align}
%   Q_1\tp = \bmat{I & 0} H_m \dots H_3 H_2 H_1,  \label{eq:Q1}
%\\ Q_2\tp = \bmat{0 & I} H_m \dots H_3 H_2 H_1.  \label{eq:Q2}
   Q_1 = Q \bmat{I \\ 0} \qquad \jbrt{\textnormal{and}} &\qquad Q_2 = Q \bmat{0 \\ I}.
\end{align}
When $A$ is sparse, the SuiteSparseQR software \cite{SuiteSparseQR} 
permutes the columns of $A\tp$ in \eqref{eq:HQR} to retain sparsity in $H_k$ and $R$.
The projection $z = Py = Q_2(Q_2\tp y)$ can then be computed efficiently.

\jbrt{One can avoid storage of $Q_1$ by noting} that $Q_1 = A\tp R\inv$.
The projection can be computed as $z = (I - Q_1Q_1\tp)y = y - A\tp R\inv R\mtp Ay$, though with lower precision than $z = Q_2(Q_2\tp y)$.

\subsection{Iterative computation of \texorpdfstring{$z$}{z}}
\label{sec:LSQR}

\jbr{%In certain situations it can occur that
Computing QR factors is sometimes not practical 
because $A$ contains one or more relatively dense columns. %(This occurred rarely in the
(In the numerical experiments of section \ref{sec:numex}, this occurred with only 2 out of 142 sparse constrained problems.)
%If $A$ has dense columns then
The multifrontal QR solver SuiteSparseQR \cite{SuiteSparseQR} then has to handle dense factors,
slowing computing times. For problems with thousands of constraints % , say $m > 1000$,
we regard column $j$ as relatively dense if
%it fulfills the condition
$ \text{nnz}(A_{:j}) / m > 0.1 $. When one expects \jbrt{the} QR \jbrt{factorization to be slow} because of dense columns, %in $A$,
an alternative is to solve the least-squares problem} % to form $z$ 
%If QR factors are not practical (e.g., if $A$ contains one or more relatively dense columns), an alternative is to solve the least-squares problem
\begin{equation}
   \label{eq:LS}
   \min_w \norm{A\tp w - y}
\end{equation}
and compute the residual $z = Py = y - A\tp w$.
Suitable iterative solvers for \eqref{eq:LS} are CGLS \cite{HS1952}, LSQR \cite{PS82a}, and LSMR \cite{FonS2011}.
If $\tilde A$ is the same as $A$ with any relatively dense columns deleted,
the factor $\tilde R$ from a sparse QR factorization of ${\tilde A}\tp$ (again with suitable
column permutation) should be a good right-preconditioner to accelerate the
iterative solvers.  If $\tilde A$ does not have full row-rank, the zero or small
diagonals of $\tilde R$ can be changed to 1 before $\tilde R$ is used as a preconditioner.

\subsection{Implementation}

Appendix B describes the implementation of the
two preceding projections. We refer to these operations through the definition
\begin{equation*}
    \b{z} \equiv
    \text{compProj}(\b{A},\b{y},\texttt{P}) \equiv
    \begin{cases}
        \text{Householder QR} & \text{if } \texttt{P} = 1, \\ % Sec. \ref{sec:QR} (App. B.1)\& 
        \text{Preconditioned LSQR} & \text{if } \texttt{P} = 2.
        % Sec. \ref{sec:LSQR} (App. B.2)\& 
    \end{cases}
\end{equation*}
Note that the implementations do not require $ \b{A} $ to have full \jbrt{row} rank.

\subsection{\texorpdfstring{\jbbo{Trust-region algorithm}}{TR algorithm}}
\label{sec:TR}
\jbbo{To solve \eqref{eq:main} we use the trust-region strategy},
which is regarded as a robust minimization method \cite{ConGT00a}.
At each iteration, the method measures progress using the ratio of
actual over predicted reductions:
\begin{equation*}
    \rho_{\jbrt{k}} = \frac{f(\bk{x}) - f(\bk{x} + \jbrt{\bk{s}})}{q(\b{0})-q(\jbrt{\bk{s}})},
\end{equation*}
where \jbrt{$ \bk{s} $} is an intermediate search direction\jbr{, %$s$ is intermediate
in the sense that $ \bk{s} $ will ultimately be used as an update only
if \jbrt{$\rho_k$ is greater than a threshold.}}  %$ \rho > \textnormal{``constant threshold"} $}. %, obtained from % () solving the trust-region subproblem \eqref{eq:trsub}. 
By accepting steps that fulfill the so-called sufficient decrease condition
$ \rho > c_1 $ (\jbrt{suppressing the subscript $k$ on $ \rho_k $}) for a constant $ c_1 > 0 $, the method successively
moves towards \jbrt{a local minimizer (though there is no guarantee that a minimizer
will be reached)}. The trust-region radius $ \Delta > 0 $ % a \jbrt{minimizer (though there is no guarantee )}
controls the norm of the search direction by means
of the constraint $ \| \b{s} \|_2 \le \Delta $. There are two possible
cases for the solution of the TR subproblem: either the search direction
is in the interior of the constraint ($ \| \b{s} \| < \Delta $)
or it is on the boundary ($ \| \b{s} \| = \Delta $). Since the L-BFGS matrix $ \bk{B} $ is positive definite, the solution of 
\eqref{eq:trsub} is given by the unconstrained minimizer $ \b{s} = \jbr{\b{s}_E} $
from \eqref{eq:KKTsystem} if $ \| \jbr{\b{s}_E} \| \le \Delta $.
Otherwise, if $ \| \jbr{\b{s}_E} \| > 0 $, then \eqref{eq:trsub}
is solved with the active norm constraint $ \| \b{s} \| = \Delta $.
Note that even if $ \| \jbr{\b{s}_E} \| \le \Delta $, the condition $ \rho > c_1 $
\jbr{might} not hold. In this situation, or in any case when
$ \rho \le c_1 $, the radius $\Delta$ is reduced and a new problem \eqref{eq:trsub}
(with smaller $ \Delta $) and constraint $ \| \b{s} \| = \Delta $ is solved.
The overall trust-region strategy for one iteration is given next, with radius
$ \Delta > 0 $ and $ c_1 > 0 $ and iteration counter suppressed.

\begin{center}
  \begin{tabular}{ll}
     \multicolumn{2}{l}{Trust-Region Strategy:}
  \\ 1. & Compute the unconstrained step $ \b{s} \leftarrow \jbr{\b{s}_E}         $ from \eqref{eq:KKTsystem} (using \eqref{eq:Vg})
          \color{white}{XXXX}
  %\\ 2. & While ($\norm{s}_2 \ge \Delta$ and $ \rho \le c_1 $)
   \\ 2. & \jbbo{While ($\norm{s}_2 > \Delta$ or $ \rho \le c_1 $)}
  \\    & \quad 2.1. Solve \eqref{eq:trsub} with $ \norm{\b{s}} = \Delta $
  \\    & \quad 2.2. Reduce $\Delta$
  \\    & end
  \\ 3. & \jbr{Increase (or at least do not decrease) $\Delta$} %Adjust (possibly increase) $ \Delta $
  \\ 4. & Update iterate $ \b{x} \leftarrow \b{x} + \b{s} $
    % \color{white}{1.} & \color{white}{Backtracking line-search: $ \bko{x} = \bk{x} - \alpha \bk{g}^P / \|\bk{g}^P \|_2 $ (cf. \cite[Alg. 3.1]{NocW06})}
  \end{tabular}
\end{center}
Practical aspects of an implementation include the setting of constants
and starting the method. Detailed procedures are described in 
sections \ref{sec:l2}, \ref{sec:eigen}, \ref{sec:SCNorm} \jbrt{and \ref{sec:numex}}.

\section{\texorpdfstring{\jbb{$\ell_2$-norm trust-region constraint}}{L2-norm}} %  \jbbo{Trust-region method with the} $\ell_2$-norm}{L2-norm}
% \jbbo{Trust-region method with the} $\ell_2$-norm
\label{sec:l2}
With an $\ell_2$-norm trust-region constraint in \eqref{eq:trsub}, the search direction
is given by
\begin{equation*}
    \biidx{L2}{s} = 
    \underset{ \| \b{s} \|_{2} \le \Delta_k }{\text{ arg min }}
    q(\b{s})
    \quad \text{subject to} \quad \b{A} \b{s} = \b{0}.
    % \b{s} \tp \bk{g} + \frac{1}{2} \b{s} \tp \bk{B} \b{s}
\end{equation*}
\jbrt{With} $\sigma \ge 0$ denoting a scalar Lagrange multiplier, the search
direction is a feasible solution to a shifted KKT system \jbrt{including the norm constraint}:
\begin{equation}
    \label{eq:KKTshifted}
    \bmat{\bk{B} + \sigma \b{I}   & \b{A} \tp
       \\ \b{A}                   & \b{0}}
    \bmat{\biidx{L2}{s} \\ \biidx{L2}{\lambda}}
    =
    \bmat{ -\bk{g} \\ \b{0}},
    \qquad
  \| \biidx{L2}{s} \|_2 \le \Delta_k.
\end{equation}
%(cf.~\cite[Section 4]{BMP19}).
By computing the (1,1) block of the shifted inverse KKT matrix,
we \jbr{note} that a necessary condition for the solution is $ \biidx{L2}{s}(\sigma) = - \bks{V} \bk{g} $, where 
\[ 
    \bks{V} = (\bk{B} + \sigma \b{I})^{-1} -
                (\bk{B} + \sigma \b{I})^{-1} \b{A} \tp
                (\b{A} (\bk{B} + \sigma \b{I})^{-1} \b{A} \tp )^{-1}
                \b{A} (\bk{B} + \sigma \b{I})^{-1}.
\]
For the L-BFGS matrix, with $ \tau_k = \tau_k(\sigma) = (1/ \delta_k + \sigma) $
we have
$ 
    (\bk{B} + \sigma \b{I})^{-1} =  
    \tau_k^{-1} \b{I} + \bk{J} \bks{W} \bk{J} \tp,
$
where the small $ 2l \times 2l $ matrix is
\[
    \bks{W} =
    - \bmat{
            \theta_k \bk{S} \tp \bk{S}          
            & \theta_k \bk{L} + \tau_k \bk{T} \\
            \theta_k \bk{L} \tp + \tau_k \bk{T} \tp     
            & \tau_k (\tau_k \bk{D} + \bk{Y} \tp \bk{Y})     
        }^{-1}
\]
with $ \theta_k = \tau_k(1- \delta_k \tau_k) $.
In terms of $ \bks{C} \equiv \b{A} \bk{J} \bks{W} $ and 
$ \bks{G} \equiv (\b{A} (\bk{B}+\sigma \b{I})^{-1} \b{A} \tp)^{-1} $,
the compact representation of $ \bks{V} $ \cite[Corollary 1]{BMP19} is
\begin{align}
    \label{eq:cmpVks_L1}
   &\bks{V} =
\\ &\frac{1}{\tau_k} I +
    %\bmat{\b{A}\tp & \bk{S} & \bk{Y} }
    \bmat{\b{A}\tp & \bk{J}}
    \bmat{-\frac{1}{\tau_k^2}\bks{G} & -   \frac{1}{\tau_k}\bks{G} \bks{C}
       \\ -\frac{1}{\tau_k}\bks{C} \tp \bks{G} & \bks{W} - \bks{C} \tp \bks{G} \bks{C}}
    \bmat{ \b{A} \\[4pt] \bk{J}\tp}.        \nonumber
        % \bmat{ \b{A} \\
        %         \bk{S}\tp \\
        %         \bk{Y} \tp}
\end{align}
% \begin{equation}
%     \bks{V} = 
%     \tau_k^{-1} \b{I} + 
% \end{equation}
Once the middle matrix in \eqref{eq:cmpVks_L1} is formed, the compact representation
can be used to compute matrix-vector products efficiently. However, when
$ m $ is large (many equality constraints), computing terms such
as $ \bks{G} $ become expensive. Therefore, we describe a reduced representation
similar to \eqref{eq:redComp}, based on \jbr{the property} that 
$ \b{P} \tp \bks{V} \b{P} = \bks{V}  $ and by storing $ \bk{S} $ and $ \bk{Z} $.
Lemma 2 summarizes the outcome.

\smallskip

{Lemma 2:} \emph{The RCR of $ \bks{V} $
in \eqref{eq:cmpVks_L1} for the L-BFGS matrix is given by}
\begin{equation}
    \label{eq:redComps}
    \bks{V} = \frac{1}{\tau_k} I + 
    \bmat{\b{A}\tp & \bk{S} & \bk{Z} }
    \bmat{ - \frac{1}{\tau_k} (\b{A} \b{A} \tp)^{-1} &  \\
            & \bks{N}
    }
    \bmat{\b{A} \\[2pt] \bk{S} \tp \\[2pt] \bk{\jbr{Z}} \tp },
\end{equation}
\emph{where} $ \jbbo{\tau_k = \tau_k(\sigma) = (1/\delta_k + \sigma)} $, $ \jbbo{\theta_k = \theta_k(\sigma) = \tau_k(\sigma)(1-\delta_k \tau_k(\sigma))} $,
\emph{and} %\Red{\ $\bks{N}$ doesn't seem to depend on $\sigma$}
\begin{equation*}
    \bks{N} = -\bmat{ \jbbo{\theta_k(\sigma)} \bk{S} \tp \bk{S} 
    & \jbbo{\theta_k(\sigma)} \bk{L} + \jbbo{\tau_k(\sigma)} \bk{T}  \\
    \jbbo{\theta_k(\sigma)} \bk{L} \tp + \jbbo{\tau_k(\sigma)} \bk{T} \tp 
    & \jbbo{\tau_k(\sigma)} (\jbbo{\tau_k(\sigma)} \bk{D} + \bk{Z} \tp \bk{Z})}^{-1}.
\end{equation*}
% \begin{equation*}
%     \bks{N} = -\bmat{ \jbb{\theta_k(\sigma)} \bk{S} \tp \bk{S} 
%     & \jbb{\theta_k(\sigma)} \bk{L} + \jbb{\tau_k(\sigma)} \bk{T}  \\
%     \jbb{\theta_k(\sigma)} \bk{L} \tp + \jbb{\tau_k(\sigma)} \bk{T} \tp 
%     & \jbb{\tau_k(\sigma)} (\jbb{\tau_k(\sigma)} \bk{D} + \bk{Z} \tp \bk{Z})}^{-1}.
% \end{equation*}

% \theta_k \bk{S} \tp \bk{S}          
%             & \theta_k \bk{L} + \tau_k \bk{T} \\
%             \theta_k \bk{L} \tp + \tau_k \bk{T} \tp     
%             & \tau_k (\tau_k \bk{D} + \bk{Y} \tp \bk{Y})

\smallskip

%{Proof.} 
\begin{proof} To simplify notation, we suppress the explicit dependence on $ \sigma $
in this proof, so that 
$ \bk{V} \equiv \bks{V}$, $\bk{C} \equiv \bks{C}$, and $\bk{W} \equiv \bks{W}$. Multiplying $ \bk{V} $ in \eqref{eq:cmpVks_L1} from the left and right
by $ \b{P} \tp $ and $ \b{P} $ yields
$$ 
\bk{V} = \frac{1}{\tau_k} \b{P} + \bmat{\bk{S} & \bk{Z} } ( \bk{W} - \bk{C} \tp \bk{G} \bk{C} ) \bmat{\bk{S} & \bk{Z} }\tp.
$$ 
Observe that 
$ \bk{C} = \b{A} \bk{J} \bk{W} = \bmat{ 0 & \b{A} \bk{Y}} \bk{W} $ is block-rectangular and that 
$
    \bk{G} = (\b{A} (\frac{1}{\tau_k} \b{I} + \bk{J} \bk{W} \bk{J} \tp) \b{A} \tp)^{-1}
$
depends on $ \bk{W} $.
Defining $ \bk{F} \equiv \tau_k (\b{A} \b{A} \tp)^{-1} $, we \jbr{show that}
the Sherman-Morrison-Woodbury (SMW) inverse gives the simplification
\begin{align*}
    &\bk{W} - \bk{C} \tp \bk{G} \bk{C} \\&=
    \bk{W} - \bk{W} \bmat{ 0 \\ \bk{Y}\tp \b{A}\tp} \bk{G} \bmat{ 0 & \b{A} \bk{Y}} \bk{W} \\
    &= \bk{W} - \bk{W} \bmat{ 0 \\ \bk{Y}\tp \b{A}\tp} 
    \big(\b{I} + \bmat{ 0 & \bk{F} \b{A} \bk{Y}} \bk{W} \bmat{ 0 \\ \bk{Y}\tp \b{A}\tp} \big)^{-1} 
    \bmat{ 0 & \bk{F}\b{A} \bk{Y}} \bk{W} \\
    &= \left( \bk{W}^{-1} + \bmat{ 0 \\ \bk{Y}\tp \b{A}\tp} \bmat{ 0 & \bk{F} \b{A} \bk{Y}}  \right)^{-1},
\end{align*}
where the third equality is obtained by applying the SMW formula in
reverse. Since only the (2,2) block in the low-rank matrix 
of the third equality is nonzero, and since 
$ \bk{F} = \tau_k (\b{A} \b{A} \tp)^{-1} $, \jbr{note} that
$$
(\bk{W}^{-1})_{22} + \bk{Y} \tp \b{A} \tp \bk{F} \b{A} \bk{Y} 
=-(\tau_k(\tau_k \bk{D} + \bk{Y}\tp \bk{Y} - 
\bk{Y}\tp \b{A} \tp (\b{A} \b{A} \tp)^{-1} \b{A}  \bk{Y}) ),
$$
which corresponds to the $(2,2)$ block $ \bks{N} $ in \eqref{eq:redComps}. Because 
all other blocks are unaffected, it holds that $ \bk{W} - \bk{C} \tp \bk{G} \bk{C} = \bks{N} $.
Subsequently, by factoring $ \b{P} = I - \b{A} \tp (\b{A} \b{A} \tp )^{-1} \b{A} $ we deduce the compact representation \eqref{eq:redComps}. 
\end{proof}
%\qed
% $
%     % \big[
%     %     \begin{smallmatrix}
%     %             \b{0} \\   
%     %             \bk{Y} \tp \b{A} \tp 
%     %     \end{smallmatrix}
%     % \big]
%     \big[
%         \begin{smallmatrix}
%                 \b{0} &   
%                 \bk{Y} \tp \b{A} \tp 
%         \end{smallmatrix}
%     \big]\tp
%     \big[
%         \begin{smallmatrix}
%                 \b{0} &   
%                 \bk{Y} \tp \b{A} \tp 
%         \end{smallmatrix}
%     \big]
% $

\smallskip
 
 Note that
 $ \bk{S}\tp \bk{Z} = \bk{S}\tp \bk{Y} = \bk{L} + \bk{D} + \bar{\b{T}}_k $,
 with $ \bk{T} = \bk{D} + \bar{\b{T}}_k $, means that the RCR
 for $ \bks{V} $ is fully specified by storing $ \bk{S} $ and $ \bk{Z} $. An exception
 is the scalar $ \delta_k $, which may depend on the most recent $ \bk{y} $. Also
 when $ \sigma = 0 $, the representations
 \eqref{eq:redComp} and \eqref{eq:redComps} coincide. We apply $ \bks{V} $ to a vector $ \b{g} $ as
 \begin{equation*}
     h = \bmat{ \bk{S} \tp \\ \bk{Z} \tp } \b{g}, \qquad
     \bks{V} \b{g} = 
     \bmat{ \bk{S} & \bk{Z} } \bks{N}h + \frac{1}{\tau_k} \b{P} \b{g}.
 \end{equation*}

\subsection{\texorpdfstring{$\ell_2$}{L2}-norm search direction}
\label{sec:l2search}
To compute the $\ell_2$ TR minimizer we first set $ \sigma = 0 $
and $ \biidx{L2}{s}(0) = - \bk{V}(0) \bk{g} $. If $ \| \biidx{L2}{s}(0) \|_2 \le \Delta_k $,
the \jbbo{minimizer with the \jbrt{$\ell_2$-norm}} is given by $ \biidx{L2}{s}(0) $. Otherwise
($ \| \biidx{L2}{s}(0) \|_2 > \Delta_k  $) we define the so-called secular equation
\cite{ConGT00a} as
$$ \phi(\sigma) \equiv \frac{1}{\| \biidx{L2}{s}(\sigma)  \|_2} - \frac{1}{\Delta_k}. $$
To solve the secular equation we apply the 1D Newton iteration
\begin{equation*}
    \sigma_{j+1} = \sigma_j - \frac{\phi(\sigma_j)}{ \phi'(\sigma_j) },
\end{equation*}
where $ \phi'(\sigma_j) = 
- (\biidx{L2}{s}(\sigma_j)\tp \biidx{L2}{s}(\sigma_j)') /  \| \biidx{L2}{s}(\sigma_j) \|^3_2  $ and
$ \biidx{L2}{s}(\sigma_j)' = -\bk{V}(\sigma_j) \biidx{L2}{s}(\sigma_j) $
(with prime `` $'$ " denoting the derivative).
\jbr{Note that $ \biidx{L2}{s}(\sigma_j)' $ can be derived from the shifted system %in 
\eqref{eq:KKTshifted} by differentiation with respect to $\sigma$. Applying the
product rule in \eqref{eq:KKTshifted} and regarding the solutions as functions of
$\sigma$, i.e., $ \biidx{L2}{s}' \equiv \biidx{L2}{s}(\sigma)' $ and
$ \biidx{L2}{\lambda}' \equiv \biidx{L2}{\lambda}(\sigma)' $, one obtains the 
differentiated system
\begin{equation*}
    %\label{eq:KKTshifted}
    \bmat{\bk{B} + \sigma \b{I}   & \b{A} \tp
       \\ \b{A}                   & \b{0}}
    \bmat{\biidx{L2}{s}' \\ \biidx{L2}{\lambda}'}
    =
    \bmat{ -\biidx{L2}{s} \\ \b{0}}.
\end{equation*}
Since the system matrix is the same as in \eqref{eq:KKTshifted} (only the
right-hand side differs), $\biidx{L2}{s}(\sigma_j)'$ is fully determined by
$\bk{V}(\sigma_j)$ and $ \biidx{L2}{s}(\sigma_j)$.}
Starting from $ \sigma_0 = 0 $, we terminate the Newton iteration if 
$ | \phi(\sigma_{j+1}) | \le \varepsilon $ or an iteration limit is reached.
The search direction is then computed as 
$ \biidx{L2}{s}(\sigma_{j+1}) = - \bk{V}(\sigma_{j+1}) \bk{g} $. 

Our approach with the $ \ell_2 $ norm is summarized in Algorithm \ref{alg:LTRL2-LEC}.
\jbbo{This algorithm is based on storing and updating $ \bk{S}, \bk{Z} $,
and the small blocks of $ \bks{N} $ in \eqref{eq:redComps}. Suppose that
%at $ k=0 $,
$ \biidx{0}{s} $ and $ \biidx{0}{z} $ are obtained by \jbrt{an} 
initialization procedure (\jbrt{for instance, Init.~1 from section \ref{sec:numex}}). With $k=0$, the
initial matrices that define $ \bks{V} $ are given as
\begin{align}
\label{eq:initSZ}
\bk{S} = \bmat{ \bk{s} }, &\quad \bk{Z} = \bmat{ \bk{z} },
\\
    \label{eq:initVks}
    \bk{D} = \bmat{ \bk{s} \tp \bk{z} }, \quad
    \bk{T} = \bmat{ \bk{s} \tp \bk{z} }, &\quad
    \bk{Z}\tp \bk{Z} = \bmat{ \bk{z} \tp \bk{z} }, \quad
    \bk{L} = \bmat{ \b{0} }.
\end{align}
Once the iteration starts, we update
\begin{equation}
\label{eq:upSZ}
\bko{S} = \text{colUpdate}(\bk{S},\bk{s}), \quad 
\bko{Z} = \text{colUpdate}(\bk{Z},\bk{z}),
\end{equation}
%and
\begin{align}
   \label{eq:upVks}
   \bko{D} &= \text{prodUpdate}(\bk{D},\b{0},\b{0}, \bk{s}, \bk{z})\jbr{,}
\\ \bko{T} &= \text{prodUpdate}(\bk{T},\bk{S},\b{0}, \bk{s}, \bk{z})\jbr{,} \nonumber
\\ \bko{Z}\tp \bko{Z} &= \text{prodUpdate}(\bk{Z}\tp \bk{Z},\bk{Z},\bk{Z}, \bk{z}, \bk{z}), \text{ \jbr{and}} \nonumber
\\ \bko{L} &= \text{prodUpdate}(\bk{L},\b{0},\bk{Z}, \bk{s}, \b{0}) \nonumber.
\end{align}
% \begin{alignat}{2}
%     \label{eq:upVks}
%     \bko{D} &= \text{prodUpdate}(\bk{D},\b{0},\b{0}, \bk{s}, \bk{z}),
%      && \bko{T}= \text{prodUpdate}(\bk{T},\bk{S},\b{0}, \bk{s}, \bk{z}), \\
%     \bko{Z}\tp \bko{Z} &= \text{prodUpdate}(\bk{Z}\tp \bk{Z},\bk{Z},\bk{Z}, \bk{z}, \bk{z}),
%      && \bko{L}= \text{prodUpdate}(\bk{L},\b{0},\bk{Z}, \bk{s}, \b{0}) \nonumber.
% \end{alignat}
Note that we store and update matrices like $ \bk{Z}\tp \bk{Z} \in \mathbb{R}^{l \times l} $ instead of recomputing them. Because of the limited memory
technique (typically $ 3 \le l \le 7 $ \cite{ByrNS94}), such matrices 
are very small relative to large $n$.
Subsequently, $ \bks{N} \in \mathbb{R}^{2l \times 2l} $, defined by the
blocks in \eqref{eq:upVks}, remains very small compared to $n$.
}

%\clearpage

% Test algorithm
%\fbox{\parbox{0.90\linewidth}{
%\begin{mdframed}
%{
%\textbf{Algorithm 1}
%\hline
%\label{alg:alg1}

\begin{algorithm}
\caption{{\small LTRL2-SLEC} (Limited-Memory Trust-Region 2-norm for Sparse Linear Equality Constraints)} % $\ell_2$
\label{alg:LTRL2-LEC}
%\begin{mdframed}
\begin{algorithmic}[1]
%\hline
\ENSURE
%{ 
$  0 < c_1$, $0 < c_2, c_3,c_4,c_5,c_6 < 1 < c_7 $, $0 < \varepsilon_1, \varepsilon_2 $,
$ 0 < i_{\text{max}} $, $ k=0 $, $ \jbrt{0< l} $, $ \Delta_k = \| \bk{x} \|_2$, 
$\bk{g} = \nabla f(\bk{x})$, $\texttt{P} \in [0,1]$, 
$ \bk{g}^P = \text{compProj}(\b{A},\bk{g},\texttt{P}) $,
$ \bko{g}^P, \bk{s}, \bk{z}, \bk{y} $ (from initialization), % (e.g., Init.1, Sec. \ref{sec:numex})
$ \bk{S}, \bk{Z}, \bk{D}, \bk{T}, \bk{L}, \bk{Z}\tp \bk{Z} $ from 
\eqref{eq:initSZ} and \eqref{eq:initVks},
% Find: $ \bko{x} = \bk{x} - \alpha \bk{g}^P / \|\bk{g}^P \|_2 $ (backtrack line-search, $\alpha>0$),
% $ \bko{g} = \nabla f(\bko{x}) $, 
% $ \bko{g}^P = \text{compProj}(\b{A},\bko{g},\texttt{P}) $,
% $ \bk{s} = \bko{x} - \bk{x} $, $ \bk{z} = \bko{g}^P - \bk{g}^P $, $ \bk{y} = \bko{g} - \bk{g} $,
$\delta_k = \bk{s}\tp \bk{z} / \bk{y} \tp \bk{y} $, 
% $ \bk{S} = \bmat{ \bk{s} } $, $\bk{Z} = \bmat{ \bk{z} }$, $ \bk{D} = \bmat{ \bk{s} \tp \bk{z} } $,
% $ \bk{T} = \bmat{ \bk{s} \tp \bk{z} } $,
% $ \bk{L} = \bmat{ \b{0} } $,
% $ \bk{Z}\tp \bk{Z} = \bmat{ \bk{z} \tp \bk{z} } $,
$ \sigma = 0 $, $\tau_k = (1/\delta_k + \sigma)$, $ \theta_k = \tau_k(1-\delta_k \tau_k) $, $ \bk{N}(\sigma) $ from \eqref{eq:redComps}, $k = k + 1$
% $\bk{g} = \nabla f(\bk{x})$,  $\delta_k = 1/\| \b{\Phi}_k(0) \b{g}_k \|_2$, $\gamma_k = \delta^{-1}_k $, $\b{\Psi}_k = \b{A}^T$, $\b{M}_k = \delta_k(\b{A}\b{A}^T)^{-1}$, %$\b{N}_k = \gamma_k \b{A}^T\bk{M} $, 
%}
\WHILE{$ (\varepsilon_1 \le \| \b{g}^P_k \|_{\infty} )  $}
\STATE{$ h = -\bmat{ \bk{S} & \bk{Z} } \tp \bk{g} $} 
\STATE{$
        \bk{s} = \bmat{ \bk{S} & \bk{Z} } \bk{N}(0)h - \delta_k \b{g}^P_k
    $; $ \rho_k = 0 $\;} \COMMENT{Equality constrained step}
%    \nllabel{alg1:eqs} \tcc*[r]{Equality constrained step} 
	 %$ \rho_k = 0 $\;	 
\IF{$ \| {s}_k \|_2 \le \Delta_k $}
	 \STATE{$\rho_k = (f(\bk{x}) - f(\bk{x}+\bk{s}))/ (q(\b{0}) - q(\bk{s})) $} % \;	
\ENDIF
%\IF{ $ \rho_k \le c_1  $ }
%{
	\WHILE{ $   \rho_k \le c_1 $ }
	%{ 	  % c_1 < \rho_k	
		\STATE{$ \sigma = 0 $, $ i = 0 $; $\tau_k = (1/\delta_k + \sigma)$, $ \theta_k = \tau_k(1-\delta_k \tau_k) $}
		\STATE{$ h' = -\bmat{ \bk{S} & \bk{Z} } \tp \bk{s} $}
    		\STATE{$\bk{s}' = \bmat{ \bk{S} & \bk{Z} } \bks{N}h' - \delta_k \bk{s}$;
    		} % $ \phi(\sigma)$, $ \phi'(\sigma) $ 
% 	$\bk{s}'(\sigma) = 
% 	-( \delta_k \b{I} + \bk{ \Psi } \bk{ M }\bk{ \Psi}^T   ) \bk{s} $\;
	\WHILE{ $ \varepsilon_2 < |\phi(\sigma)| $ \AND $ i < i_{\max} $ } % \AND 
	%\COMMENT{Newton's method}
	%{ % (\tcc*[f]{Newton's method })
		\STATE{$ \sigma =  \sigma - \phi(\sigma) / \phi'(\sigma)  $}
		\STATE{$\tau_k = (1/\delta_k + \sigma)$, $ \theta_k = \tau_k(1-\delta_k \tau_k) $}
		\STATE{$ h = -\bmat{ \bk{S} & \bk{Z} } \tp \bk{g} $;
    				$
        					\bk{s} = \bmat{ \bk{S} & \bk{Z} } \bks{N}h - \frac{1}{\tau_k} \b{g}^P_k
    				$}
% 		$\bk{s} = \bk{s}(\sigma) = 
% 		-((\gamma_k + \sigma)^{-1} \b{I} + \bk{ \Psi } \bk{ M }(\sigma) \bk{ \Psi}^T   ) \b{g}_k $ \nllabel{alg1:ss}\;
% 		$\bk{s}' = \bk{s}'(\sigma) = 
% 		-( (\gamma_k + \sigma)^{-1} \b{I} + \bk{ \Psi } \bk{ M }(\sigma) \bk{ \Psi}^T   ) \bk{s} $ \nllabel{alg1:ssp}\;
       		\STATE{$ h' = -\bmat{ \bk{S} & \bk{Z} } \tp \bk{s} $;
        		$
            	\bk{s}' = \bmat{ \bk{S} & \bk{Z} } \bks{N}h' - \frac{1}{\tau_k} \bk{s}
        		$; } % $ \phi(\sigma)$, $ \phi'(\sigma) $ 
		\STATE{$i = i+1$} 	  % \leftarrow
	%}
	\ENDWHILE
	\COMMENT{Newton's method}
	\STATE{$\rho_k = 0$}
	\IF{$ 0 < (f(\bk{x}) - f(\bk{x}+\bk{s}))  $}
	%{
		\STATE{$\rho_k = (f(\bk{x}) - f(\bk{x}+\bk{s}))/ (q(\b{0}) - q(\bk{s})) $}
	%}
	\ENDIF		
	\IF{ $\rho_k \le c_2 $}
	%{
		\STATE{$ \Delta_k = \min(c_3 \| \bk{s} \|_2, c_4 \Delta_k)$} %\tcc*[r]{Reduce radius}
	%}
	\ENDIF
	%}
	\ENDWHILE
%}
%\ENDIF
\STATE{$\b{x}_{k+1} = \bk{x} + \bk{s}$}
\COMMENT{Accept step} %\tcc*[r]{Accept step}
\IF{$c_5\Delta_k \le \| \bk{s} \|_2$ \AND $ c_6 \le \rho_k $}
%{
	\STATE{$\Delta_k = c_7 \Delta_k$}
%}
\ENDIF
\STATE{$\bko{g} = \nabla f(\bko{x}) $, $ \bko{g}^P = \text{compProj}(\b{A},\bko{g},\texttt{P}) $, 
$ \bk{z} = \bko{g}^P - \bk{g}^P $,
$ \bk{y} = \bko{g} - \bk{g} $, 
$ \bko{S}, \bko{Z}, \bko{D}, \bko{T}, \bko{L}, \bko{Z}\tp \bko{Z} $ from 
\eqref{eq:upSZ} and \eqref{eq:upVks}
% $\bko{S} = \text{colUpdate}(\bk{S},\bk{s})$,
% $\bko{Z} = \text{colUpdate}(\bk{Z},\bk{z})$,
% $\bko{D} = \text{prodUpdate}(\bk{D},\b{0},\b{0}, \bk{s}, \bk{z})$,
% $\bko{T} = \text{prodUpdate}(\bk{T},\bk{S},\b{0}, \bk{s}, \bk{z})$,
% $\bko{L} = \text{prodUpdate}(\bk{L},\b{0},\bk{Z}, \bk{s}, \b{0})$,
% $\bko{S}\tp \bko{S} = \text{prodUpdate}(\bk{S}\tp \bk{S},\bk{S},\bk{S}, \bk{s}, \bk{s})$,
% $\bko{S}\tp \bko{Z} = \text{prodUpdate}(\bk{S}\tp \bk{Z},\bk{S},\bk{Z}, \bk{s}, \bk{z})$,
% $\bko{Z}\tp \bko{Z} = \text{prodUpdate}(\bk{Z}\tp \bk{Z},\bk{Z},\bk{Z}, \bk{z}, \bk{z})$\;
$ \delta_{k+1} = \bk{z}\tp \bk{s} / \bk{y}\tp \bk{y}  $,
$ \sigma = 0 $, $\tau_k = (1/\delta_k + \sigma)$, $ \theta_k = \tau_k(1-\delta_k \tau_k) $}
\STATE{Update $ \bks{N}$ from \eqref{eq:redComps}, $k = k+1$}
\ENDWHILE
% $ \bko{g} = \nabla f(\bko{x}) $ and $ \bk{y} = \bko{g} - \bk{g} $\;
% $ \delta_{k+1} = \bk{y}^T \bk{s} / \bk{s}^T \bk{s}  $ and $ \gamma_{k+1} = \delta^{-1}_{k+1}  $\;
% Update $ \bk{\Psi}, \bk{M} $ from \eqref{eq:compKi},  $ k = k+1$;
\end{algorithmic}
%\end{mdframed}
\end{algorithm}

\section{Eigendecomposition of \texorpdfstring{$ \bk{V} $}{Vk}}
\label{sec:eigen}

We describe how to exploit the structure of the RCR \eqref{eq:redComp}
to compute an implicit eigendecomposition of $ \bk{V} $, and how to
combine this with a shape-changing norm. The effect is that the trust-region subproblem solution
is given by an analytic formula.
Since the RCR is equivalent to representation
\eqref{eq:cmpVk_L1}, we can apply previous results. However, % apply the approach from \cite[Section 5]{BMP19}
using representation \eqref{eq:redComp} is computationally more efficient.
First, note that $ \bk{N} \in \mathbb{R}^{2l \times 2l} $ is a small symmetric
square matrix.
%(typically $l \in [3, 7]$ cf.~\cite{ByrNS94}).
Therefore, computing the 
nonzero eigenvalues and corresponding eigenvectors of the matrix
$
    \bmat{ \bk{S} & \bk{Z} } \bk{N} \bmat{ \bk{S} & \bk{Z} }\tp 
    = \biidx{2}{U} \biidx{2}{\Lambda} \biidx{2}{U} \tp
$ 
is inexpensive. In particular, we compute the thin QR factorization
$ \bmat{\bk{S} & \bk{Z}} =  \biidxh{2}{Q}\biidxh{2}{R}  $ and the small eigendecomposition
$ 
    \biidxh{2}{R} \bk{N} \biidxh{2}{R} \tp = \biidxh{2}{P} \biidx{2}{\Lambda} \biidxh{2}{P} \tp.
$
The small factorization is then
$$
    \bmat{ \bk{S} & \bk{Z} } \bk{N} \bmat{ \bk{S} & \bk{Z} }\tp =
    \biidxh{2}{Q} (\biidxh{2}{R} \bk{N} \biidxh{2}{R} \tp) \biidxh{2}{Q} \tp 
    = \biidxh{2}{Q} (\biidxh{2}{P} \biidx{2}{\Lambda} \biidxh{2}{P} \tp) \biidxh{2}{Q} \tp
    \equiv \biidx{2}{U} \biidx{2}{\Lambda} \biidx{2}{U} \tp,
$$
where the orthonormal matrix on the right-hand side is defined as 
$ \biidx{2}{U} \equiv  \biidxh{2}{Q} \biidxh{2}{P} $. %  and $\biidx{2}{\Lambda}$ is diagonal
Since $ \b{A} \tp (\b{A} \b{A} \tp)^{-1} \b{A} = \biidx{1}{Q} \biidx{1}{Q} \tp $
from \eqref{eq:QR}, we express $ \bk{V} $ as
\begin{equation*}
    \bk{V} = \delta_k \b{I} +
        \bmat{ \biidx{1}{Q} & \biidx{2}{U} }
        \bmat{ - \delta_k \b{I}_m   & \\
                                    & \biidx{2}{\Lambda}}
        \bmat{ \biidx{1}{Q} \tp \\ \biidx{2}{U} \tp },                            
\end{equation*}
where $ \biidx{1}{Q} \in \mathbb{R}^{n \times m} $ and
$ \biidx{2}{U} \in \mathbb{R}^{n \times 2l} $ are orthonormal,
while $ \biidx{2}{\Lambda} \in \mathbb{R}^{2l \times 2l} $ is diagonal.
Defining the orthogonal matrix 
$ 
    \b{U} \equiv  
    \bmat{\biidx{1}{Q} & \biidx{2}{U} & \biidx{3}{U}},
$
where $\biidx{3}{U} \in \mathbb{R}^{n \times n-(m+2l)}$ represents
the orthogonal complement of $ \bmat{\biidx{1}{Q} & \biidx{2}{U}} $,
we obtain the implicit eigen\-decomposition of $ \bk{V} $ as
\begin{equation}
    \label{eq:eigVk}
    \bk{V} =
    \bmat{ \biidx{1}{Q} & \biidx{2}{U} & \biidx{3}{U} }
    \bmat{ \b{0}_{m}    & & \\
                        & \delta_k \b{I}_{2l} + \biidx{2}{\Lambda}  & \\
                        &                                           & \delta_k \b{I}_{n-(m+2l)}}
    \bmat{ \biidx{1}{Q} \tp \\ \biidx{2}{U}\tp \\ \biidx{3}{U} \tp}
    \equiv \b{U} \b{\Lambda} \b{U} \tp.
\end{equation}

Note that %in factorization \eqref{eq:eigVk}
we do not 
explicitly form the potentially expensive \jbrt{to compute} orthonormal matrix $ \biidx{3}{U} $,
as only scaled projections $ \delta_k \biidx{3}{U} \biidx{3}{U} \tp  $ are needed.
We therefore refer to factorization \eqref{eq:eigVk} as being implicit. In particular,
from the identity $ \b{U} \b{U} \tp = \b{I} $, we obtain that 
$
    \biidx{3}{U} \biidx{3}{U} \tp = \b{I} - \biidx{1}{Q} \biidx{1}{Q} \tp - \biidx{2}{U} \biidx{2}{U} \tp = \b{P} - \biidx{2}{U} \biidx{2}{U} \tp.
$
\jbbo{Note here and above that $ \biidx{2}{U} $ is a thin rectangular matrix with only $2l$ columns.}
%, since its number of columns are the limited memory size $l$.}

\section{\texorpdfstring{\jbb{Shape-changing-norm trust-region constraint}}{Shape-changing-norm}}
%\subsection{\texorpdfstring{ \jbbo{Trust-region method with the} shape-changing norm}{Shape-changing-norm}} %  search direction
\label{sec:SCNorm}

To make use of the implicit eigensystem \eqref{eq:eigVk},
we apply the so-called shape-changing infinity norm introduced
in \cite{BurdakovLMTR16}:
\begin{equation*}
    \| \b{s} \|_{\b{U}} \equiv
    \text{max}\left\{ \Norm{\bmat{ \biidx{1}{Q} & \biidx{2}{U} }\tp \! \b{s}}_{\infty}, \Norm{\biidx{3}{U} \tp \b{s}}_2 \right\}.
\end{equation*}
With this norm, the trust-region subproblem has  \jbrt{a computationally efficient} solution that
can be \jbrt{obtained from}
\begin{equation*}
    \biidx{SC}{s} = 
    \underset{ \| \b{s} \|_{\b{U}} \le \Delta_k }{\text{ arg min }}
    q(\b{s})
    \quad \text{subject to} \quad \b{A} \b{s} = \b{0}.
    % \b{s} \tp \bk{g} + \frac{1}{2} \b{s} \tp \bk{B} \b{s}
\end{equation*}
Since the RCR is equivalent to
\eqref{eq:cmpVk_L1}, we invoke \cite[section 5.5]{BMP19} to obtain
 an \jbrt{direct} formula for the search direction:
\begin{equation*}
    \biidx{SC}{s} = \biidx{2}{U}(\biidx{2}{v} - \beta \biidx{2}{U}\tp \bk{g} )
    + \beta \b{P} \bk{g}, %\quad
    % (\biidx{2}{v})_i =
    % \begin{cases}
    %     \frac{-(\biidx{2}{U}\tp \bk{g})_i}{(\delta_k + (\biidx{2}{\Lambda})_{ii})^{-1}}
    %     & \text{ if } \left| \frac{(\biidx{2}{U}\tp \bk{g})_i}{(\delta_k + (\biidx{2}{\Lambda})_{ii})^{-1}}  \right| \le \Delta_k
    % \end{cases}
\end{equation*}
where with
$ \biidx{2}{U} \tp \bk{g} 
    = \biidxh{2}{P}\tp \biidxh{2}{R}^{-\top} \bmat{ \bk{S} & \bk{Z} }\tp \bk{g} \equiv \bk{u}  $, and
    $ \mu_i = (\delta_k + (\biidx{2}{\Lambda})_{ii})^{-1} $,
\begin{alignat}{2}
    \label{eq:SCanalytic1}
    &&(\biidx{2}{v})_i =&
    \begin{cases}
        \frac{-(\bk{u})_i}{\mu_i} % \biidx{2}{U}\tp \bk{g}
        & \text{ if } \left| \frac{(\bk{u})_i}{\mu_i}  \right| \le \Delta_k, \\
        \frac{-\Delta_k (\bk{u})_i}{|(\bk{u})_i|} % \biidx{2}{U}\tp \bk{g}
        & \text{ otherwise},
    \end{cases} \\
    \label{eq:SCanalytic2}
    &&\beta =& 
    \begin{cases}
        -\delta_k
        & \text{ if } \norm{\delta_k \biidx{3}{U} \tp \bk{g}}_2 \le \Delta_k,
          \\ \frac{-\Delta_k}{ \|\biidx{3}{U}\tp \bk{g} \|_2}
        & \text{ otherwise},
    \end{cases}
\end{alignat}
for $ 1 \le i \le 2l $. More details for the computation of $ \b{s}_{SC} $
are in Appendix C.
Note that the norm $ \| \biidx{3}{U}\tp \bk{g} \|_2 $ can be computed without explicitly
forming $ \biidx{3}{U} $, since 
$ 
    \| \biidx{3}{U}\tp \bk{g} \|^2_2 = \bk{g} \tp (\b{P} - \biidx{2}{U}\biidx{2}{U}\tp) \bk{g}
    = \| \b{P} \bk{g} \|_2^2 - \| \biidx{2}{U}\tp \bk{g} \|_2^2.
$
The trust-region algorithm using the RCR and the shape-changing norm
is summarized in Algorithm \ref{alg:LTRSC-LEC} below.
\jbbo{Like Algorithm \ref{alg:LTRL2-LEC}, this algorithm is based on %  \Red{1?}
storing and updating $ \bk{S}, \bk{Z} $ and the small blocks of 
$ \bk{N} $ in \eqref{eq:redComp}. Therefore, the initializations
 \eqref{eq:initSZ}--\eqref{eq:initVks} and updates 
 \eqref{eq:upSZ}--\eqref{eq:upVks} can be used. In addition,
 since in the thin QR factorization 
 $ \bmat{ \bk{S} & \bk{Z} } = \bh{Q}_2 \bh{R}_2 $ the triangular
 $ \bh{R}_2 $ is computed from a Cholesky factorization of 
 $ \bmat{ \bk{S} & \bk{Z} }\tp \bmat{ \bk{S} & \bk{Z} } $, we %also
 initialize the matrices
 \begin{equation}
 \label{eq:initVk}
 \bk{S}\tp \bk{S} = \bmat{ \bk{s}\tp \bk{s} }, \quad
 \bk{S}\tp \bk{Z} = \bmat{ \bk{s}\tp \bk{z} },
 %\Red{\quad \text{Should }k=1?}
 \end{equation}
 with corresponding updates
 \begin{align}
 \label{eq:upVk}
 \bko{S}\tp \bko{S} &= \text{prodUpdate}(\bk{S}\tp \bk{S},\bk{S},\bk{S}, \bk{s}, \bk{s}), \textnormal{ \jbr{and} } \\
 \bko{S}\tp \bko{Z} &= \text{prodUpdate}(\bk{S}\tp \bk{Z},\bk{S},\bk{Z}, \bk{s}, \bk{z}) \nonumber.
 \end{align}
 As before, with a small memory parameter $ l $, these matrices are very
 small compared to large $n$, and computations with them are inexpensive.
 }

% \begin{equation*}
%     (\biidx{2}{v})_i =
%     \begin{cases}
%         \frac{-(\biidx{2}{U}\tp \bk{g})_i}{\mu_i}
%         & \text{ if } \left| \frac{(\biidx{2}{U}\tp \bk{g})_i}{\mu_i}  \right| \le \Delta_k \\
%         \frac{-\Delta_k}{|(\biidx{2}{U}\tp \bk{g})_i|}
%         & \text{ otherwise },
%     \end{cases}
%     %\quad 
%     % \beta = 
%     % \begin{cases}
%     %     \delta_k
%     %     & \text{ if } \left| \delta_k \biidx{3}{U} \tp \bk{g}  \right| \le \Delta_k \\
%     %     \frac{-\Delta_k}{ \|\biidx{3}{U}\tp \bk{g} \|_2}
%     %     & \text{ otherwise }.
%     % \end{cases},
% \end{equation*}
% and
% \begin{equation*}
%     \beta = 
%     \begin{cases}
%         \delta_k
%         & \text{ if } \left| \delta_k \biidx{3}{U} \tp \bk{g}  \right| \le \Delta_k \\
%         \frac{-\Delta_k}{ \|\biidx{3}{U}\tp \bk{g} \|_2}
%         & \text{ otherwise }.
%     \end{cases}
% \end{equation*}

%$ \bmat{\bk{S} & \bk{Z}} =  \biidx{SZ}{Q} \biidx{SZ}{R}  $

%  \begin{equation*}
%      h = \bmat{ \bk{S} \tp \\ \bk{Z} \tp } \b{g}, \qquad
%      \bk{V} \b{g} = 
%      \bmat{ \bk{S} & \bk{Z} } \bk{N}h + \delta_k \b{P} \b{g}.
%  \end{equation*}

\begin{algorithm}%[p]
\caption{%{\small TR1x}  
{\small LTRSC-SLEC} (Limited-Memory Trust-Region Shape-Changing Norm for Sparse Linear Equality Constraints)}
\label{alg:LTRSC-LEC}
\begin{algorithmic}[1]
%\SetAlgoNoLine % No outside lines around blocks 
%\SetKwInput{init}{In} % Initialization
%\SetNlSty{}{}{} % Linenumbering style {<font>}{<txt before>}{<txt after>}
\ENSURE
{ 
$  0 < c_1$, $0 < c_2, c_3,c_4,c_5,c_6 < 1 < c_7 $, $0 < \varepsilon_1 $, $ \jbrt{0 < l} $, $ k=0 $, $ \Delta_k = \| \bk{x} \|_2$,  $\bk{g} = \nabla f(\bk{x})$, $\texttt{P} \in [0,1]$, 
$ \bk{g}^P = \text{compProj}(\b{A},\bk{g},\texttt{P}) $, 
$ \bko{g}^P, \bk{s}, \bk{z}, \bk{y} $ (from initialization), % (e.g., Init.1, Sec. \ref{sec:numex}),
$ \bk{S}, \bk{Z}, \bk{D}, \bk{T}, \bk{Z}\tp \bk{Z},
\bk{S}\tp \bk{S}, \bk{S}\tp \bk{Z}$ from 
\eqref{eq:initSZ}, \eqref{eq:initVks} and \eqref{eq:initVk},
% Find: $ \bko{x} = \bk{x} - \alpha \bk{g}^P / \|\bk{g}^P \|_2 $ (backtrack line-search, $\alpha>0$),
% $ \bko{g} = \nabla f(\bko{x}) $, 
% $ \bko{g}^P = \text{compProj}(\b{A},\bko{g},\texttt{P}) $,
% $ \bk{s} = \bko{x} - \bk{x} $, $ \bk{z} = \bko{g}^P - \bk{g}^P $, $ \bk{y} = \bko{g} - \bk{g} $,
$\delta_k = \bk{s}\tp \bk{z} / \bk{y} \tp \bk{y} $, 
% $ \bk{S} = \bmat{ \bk{s} } $, $\bk{Z} = \bmat{ \bk{z} }$, $ \bk{D} = \bmat{ \bk{s} \tp \bk{z} } $,
% $ \bk{T} = \bmat{ \bk{s} \tp \bk{z} } $,
% $ \bk{S}\tp \bk{S} = \bmat{ \bk{s} \tp \bk{s} } $,
% $ \bk{S}\tp \bk{Z} = \bmat{ \bk{s} \tp \bk{z} } $,
% $ \bk{Z}\tp \bk{Z} = \bmat{ \bk{z} \tp \bk{z} } $,
$ \bk{N} $ from \eqref{eq:redComp}, $k = k + 1$
}
%$\b{\Psi}_k = \b{A}^T$, $\b{M}_k = \delta_k(\b{A}\b{A}^T)^{-1}$} %$\b{N}_k = \gamma_k \b{A}^T\bk{M} $,
\WHILE{$ (\varepsilon_1 \le \| \b{g}^P_k \|_{\infty} )  $} 
%{
    \STATE{$ h = -\bmat{ \bk{S} & \bk{Z} } \tp \bk{g} $}
% 	\STATE{$
% 	\bk{s} = -\left( \delta_k \b{I} + \bk{ \Psi } \bk{ M } \bk{ \Psi}^T   \right) \bk{g}
% 	$}
    	\STATE{$
        \bk{s} = \bmat{ \bk{S} & \bk{Z} } \bk{N}h - \delta_k \b{g}^P_k; $ \bk{\rho} = 0 $;
    	$}
	\COMMENT{Equality constrained step}
	%\nllabel{alg2:eqs} \tcc*[r]{Equality constrained step} 	
\IF{$ \| {s}_k \|_2 \le \Delta_k $}
%{
	 \STATE{$\rho_k = (f(\bk{x}) - f(\bk{x}+\bk{s}))/ (q(\b{0}) - q(\bk{s})) $;
	 $ \| {s}_k \| = \| {s}_k \|_2 $}
%}
\ENDIF
\IF{ $ \rho_k \le c_1  $ }
%{ %\bmat{ \bk{S} & \bk{Z} }\tp \bmat{ \bk{S} & \bk{Z}}
	\STATE{$
	    \bh{R}_2\tp\bh{R}_2 =  
	        \bigg[
	            \begin{smallmatrix}
	                \bk{S}\tp \bk{S} & \bk{S}\tp \bk{Z} \\
	                \bk{Z}\tp \bk{S} & \bk{Z}\tp \bk{Z}
	            \end{smallmatrix}
	        \bigg]
	$}
	\COMMENT{Cholesky factorization}
	%\nllabel{alg2:R}  \tcc*[r]{Cholesky factorization}%{LDL$^\text{T}$ or Cholesky factorization}
	\STATE{$ \bh{P}_2 \b{\Lambda}_2 \bh{P}_2\tp  = \bh{R}_2 \bk{N} \bh{R}_2\tp $} % \nllabel{alg2:V} \tcc*[r]{Eigendecomposition}
	\COMMENT{Eigendecomposition}
	\STATE{$\bk{u}= \bh{P}_2\tp \bh{R}_2^{-\top} \bmat{ \bk{S} & \bk{Z} } \tp \bk{g} $}
	\STATE{$ \xi_k = 
		(  \| \bk{g}^P \|_2^2   - \| \bk{u} \|_2^2 )^{\frac{1}{2}}$}
	\WHILE{ $   \rho_k \le c_1 $ }
	%{
	%$ \Delta_k = d_{-} \Delta_k $\;
	\STATE{Set $ \b{v}_2 $ from \eqref{eq:SCanalytic1} using $ \bk{u} $, $ \b{\Lambda}_2 $}
	\STATE{Set $ \beta $ from \eqref{eq:SCanalytic2} using $ \xi_k = \| \biidx{3}{U}\tp \bk{g} \|_2 $}
	%Set $ \beta $ from \eqref{eq:subsolnperp}\; %with $ \| \bper{g} \|_2 $\;
	\STATE{$	\bk{s} = \bmat{ \bk{S} & \bk{Z} }
		\bh{R}_2^{-1} \bh{P}_2 (\b{v}_2 - \beta \bk{u}) + \beta\bk{g}^P$; % \nllabel{alg2:ss}
	\quad $\rho_k = 0$}
	\IF{$ 0 < (f(\bk{x}) - f(\bk{x}+\bk{s}))  $}
	%{
		\STATE{$\rho_k = (f(\bk{x}) - f(\bk{x}+\bk{s}))/ (q(\b{0}) - q(\bk{s})) $}
	%}
	\ENDIF
	%$\rho_k = (f(\bk{x}) - f(\bk{x}+\bk{s}))/ (Q(\b{0}) - Q(\bk{s})) $\;
			\IF{ $\rho_k \le c_2 $}
			%{
				\STATE{$ \Delta_k = \min(c_3 \| \bk{s} \|_{\b{U}}, c_4 \Delta_k)$} %\tcc*[r]{Reduce radius}
			%}
			\ENDIF
	%}
	\ENDWHILE
	\STATE{$ \| \bk{s} \| = \| \bk{s} \|_U $}
%}
\ENDIF
\STATE{$\b{x}_{k+1} = \bk{x} + \bk{s}$}% \tcc*[r]{Accept step}
\COMMENT{Accept step}
\IF{$c_5\Delta_k \le \| \bk{s} \|$ {\bf and} $ c_6 \le \rho_k $}
%{
	\STATE{$\Delta_k = c_7 \Delta_k$}
%}
\ENDIF
\STATE{$\bko{g} = \nabla f(\bko{x}) $, $ \bko{g}^P = \text{compProj}(\b{A},\bko{g},\texttt{P}) $, 
$ \bk{z} = \bko{g}^P - \bk{g}^P $,
$ \bk{y} = \bko{g} - \bk{g} $, 
$ \bko{S}, \bko{Z}, \bko{D}, \bko{T}, \bko{Z}\tp \bko{Z},
\bko{S}\tp \bko{S}, \bko{S}\tp \bko{Z}$ from 
\eqref{eq:upSZ}, \eqref{eq:upVks} and \eqref{eq:upVk};
% $\bko{S} = \text{colUpdate}(\bk{S},\bk{s})$,
% $\bko{Z} = \text{colUpdate}(\bk{Z},\bk{z})$,
% $\bko{D} = \text{prodUpdate}(\bk{D},\b{0},\b{0}, \bk{s}, \bk{z})$,
% $\bko{T} = \text{prodUpdate}(\bk{T},\bk{S},\b{0}, \bk{s}, \bk{z})$,
% $\bko{S}\tp \bko{S} = \text{prodUpdate}(\bk{S}\tp \bk{S},\bk{S},\bk{S}, \bk{s}, \bk{s})$,
% $\bko{S}\tp \bko{Z} = \text{prodUpdate}(\bk{S}\tp \bk{Z},\bk{S},\bk{Z}, \bk{s}, \bk{z})$,
% $\bko{Z}\tp \bko{Z} = \text{prodUpdate}(\bk{Z}\tp \bk{Z},\bk{Z},\bk{Z}, \bk{z}, \bk{z})$\;
$ \delta_{k+1} = \bk{z}\tp \bk{s} / \bk{y}\tp \bk{y}  $}
\STATE{Update $ \bk{N}$ from \eqref{eq:redComp}; $k = k+1$}
%}
%\ENDIF
\ENDWHILE
\end{algorithmic}
\end{algorithm}

\section{Numerical experiments}
\label{sec:numex}
The numerical experiments are carried out in MATLAB 2016a on a 
MacBook Pro @2.6 GHz Intel Core i7 with 32 GB of memory. For comparisons, we use the implementations of 
\rfm{Algorithms 1 and 2 from \cite{BMP19}, which we label \TR1 and \TR2.} %TR1 and TR2 
All codes are available in the public domain:
\begin{center}
\url{https://github.com/johannesbrust/LTR_LECx}
\end{center}
For \TR1, \TR2 we use the modified stopping criterion
$ \| \b{P} \bk{g} \|_{\infty} \le \epsilon $ in place of $\| \b{P} \bk{g} \|_2 / \text{max}(1, \b{x}_k ) \le \epsilon $ in order to compare consistently across solvers. Unless otherwise specified,
the default parameters of these two algorithms are used. %The names of our proposed algorithms
%are as follows:
\rfm{We use the following names for our proposed algorithms:}
\begin{center}
    \begin{tabular}{ll}
    \jbbo{\TR1H}: &Alg.~\ref{alg:LTRL2-LEC} with representation  \eqref{eq:redComps} and Householder QR
    \\ \jbbo{\TR1L}: &Alg.~\ref{alg:LTRL2-LEC} with representation  \eqref{eq:redComps} and preconditioned LSQR
    \\ \jbbo{\TR2H}: & Alg.~\ref{alg:LTRSC-LEC} with representation \eqref{eq:redComp} and Householder QR
    \\ \jbbo{\TR2L}: &Alg.~\ref{alg:LTRSC-LEC} with representation  \eqref{eq:redComp} and preconditioned LSQR
    \end{tabular}
\end{center}

Note that \TR1 and \TR2 were developed for low-dimensional linear equality constraints. In addition, we include {\small IPOPT } \cite{WaechterBiegler06} with an L-BFGS quasi-Newton matrix (we use a precompiled Mex file with {\small IPOPT } 3.12.12 that includes MUMPS and MA57 libraries). \jbr{We note that a commercial state-of-the-art quasi-Newton trust-region solver that uses a projected conjugate gradient solver is implemented in the {\small K{\footnotesize NITRO}-I{\footnotesize NTERIOR}/CG} \cite[Algorithm 3.2]{Byrd2006}.}
For \jbr{the freely available} {\small IPOPT } we specify the limited-memory BFGS option using the option
\texttt{hessian\_approximation=`limited memory'} with 
\texttt{tol=1e-5}. %, % and 
%\texttt{dual\_inf\_tol=1e-5} 
\jbr{(The parameter \texttt{tol} is used by {\small IPOPT } to ensure that the (scaled) projected gradient in the infinity norm and the constraint violation are below the specified threshold.
The default value is \texttt{tol=1e-8}.)}
% \jbr{
% (here \texttt{dual\_inf\_tol}} \jbr{ 
% ensures that the projected gradient satisfied $ \| \b{P} \bk{g} \|_{\infty} \le \texttt{1e-5} $ enabling comparisons, while \texttt{tol} ensures that the residual $ \| A\bk{x}-\b{b} \| $ is
% sufficiently small. IPOPT's default value is $ \texttt{tol=1e-8} $)}.
All other parameters in {\small IPOPT } are at their default values unless otherwise
specified. The parameters in \TR{1\{H,L\}} and \TR{2\{H,L\}} are set to $ c_1 $ (as machine epsilon), $ c_2 = 0.75 $, $ c_3 = 0.5 $, $ c_4 = 0.25 $, $ c_5 = 0.8 $, $ c_6 = 0.25 $, $ c_7 = 2 $, and  $i_{\text{max}}=10$.
The limited-memory parameter of all compared TR solvers is set to $l=5$ ({\small IPOPT }'s default is $l=6$). Because the proposed
methods are applicable to problems with a large number of constraints, problems with large dimensions such as $m \ge 10^4$, $n \ge 10^5$ are included. Throughout the experiments, $\b{A} \in \mathbb{R}^{m \times n}$ with $m < n$.

\jbrt{To initialize the algorithm, we distinguish two main cases.} If $ \biidx{0}{x} $ is not available, it 
is computed as the 
\jbbo{minimum-norm solution
$
    \biidx{0}{x} = \text{argmin}_{\b{x}} \norm{\b{x}}_2 \text{~s.t.~} \b{A}\b{x} = \b{b}.
$
(e.g., $\biidx{0}{x}=\b{A}\tp(\b{A}\b{A}\tp)^{-1}\b{b} $ when 
$\b{A}$ is full rank.)}
% least-squares solution
% $
%     %\biidx{0}{x} = \underset{\b{x}}{\text{ argmin }} \| \b{A}\b{x} - \b{b} \|_2.
%     \biidx{0}{x} = \text{argmin}_{\b{x}} \norm{\b{A}\b{x} - \b{b}}_2.
% $
%\Red{This is not unique.}
If $ \bh{x}_0 $ is provided but is infeasible, the initial
vector can be computed from \jbbo{
$
    \biidx{0}{p} = \text{argmin}_{\b{p}} \norm{\b{p}}_2 \text{~s.t.~}
    \b{A}\b{p} = \b{b} -\b{A} \bh{x}_0
$}
% $
%     \biidx{0}{p} = \text{argmin}_{\b{p}} \norm{\b{A}\b{p} + (\b{A} \bh{x}_0 - \b{b})}_2
% $
and $ \biidx{0}{x} = \bh{x}_0 + \biidx{0}{p}$. 
To compute the initial vectors 
$ \biidx{0}{s} = \biidx{1}{x} - \biidx{0}{x} $,
$ \biidx{0}{z} = \b{P}\biidx{1}{g} - \b{P}\biidx{0}{g} $,
and
$ \biidx{0}{y} = \biidx{1}{g} - \biidx{0}{g} $
we determine an initial $ \biidx{1}{x} $ value also.
Suppose that at $ k = 0 $, all of
$ \bk{x} $, $ \bk{g} = \nabla f(\bk{x}) $ and 
$ \bk{g}^{P} = \b{P} \bk{g} $ are known. An initialization for 
$ \bk{s}$, $\bk{z}$ and $ \bk{y} $ at $ k =0 $ is the following:
\begin{center}
   \begin{tabular}{ll}
   \multicolumn{2}{l}{Init.~1:}
\\ 1. & Backtracking line-search: $ \bko{x} = \bk{x} - \alpha \bk{g}^P / \|\bk{g}^P \|_2 $ (cf.\ \cite[Alg.\ 3.1]{NocW06})
\\ 2. & $ \bko{g} = \nabla f(\bko{x}) $, $ \bko{g}^P = \text{compProj}(\b{A},\bko{g},\texttt{P}) $
\\ 3. & $ \bk{s} = \bko{x} - \bk{x} $
\\ \textcolor{white}{3.} & $ \bk{z} = \bko{g}^P - \bk{g}^P $
\\ \textcolor{white}{3.} & $ \bk{y} = \bko{g} - \bk{g} $ %\\
        % \color{white}{1.} & \color{white}{Backtracking line-search: $ \bko{x} = \bk{x} - \alpha \bk{g}^P / \|\bk{g}^P \|_2 $ (cf. \cite[Alg. 3.1]{NocW06})}
    % 3.  & $ \bk{s} = \bko{x} - \bk{x} $, $ \bk{z} = \bko{g}^P - \bk{g}^P $, $ \bk{y} = \bko{g} - \bk{g} $
%     4.  & Update iterate $ \b{x} \leftarrow \b{x} + \b{s} $
%     Find: $ \bko{x} = \bk{x} - \alpha \bk{g}^P / \|\bk{g}^P \|_2 $ (backtrack line-search, $\alpha>0$),
% $ \bko{g} = \nabla f(\bko{x}) $, 
% $ \bko{g}^P = \text{compProj}(\b{A},\bko{g},\texttt{P}) $,
% $ \bk{s} = \bko{x} - \bk{x} $, $ \bk{z} = \bko{g}^P - \bk{g}^P $, $ \bk{y} = \bko{g} - \bk{g} $,
  \end{tabular}
\end{center}
Once $\biidx{0}{s}$, $\biidx{0}{z}$ and $\biidx{0}{y}$ have been initialized (with initial radius $ \Delta_0 = \| \biidx{0}{s} \|_2 $), all other updates are done automatically within the trust-region
strategy.%}
% and to IPOPT \cite{WaechterBiegler06} with IPOPT's quasi-Newton option. 

% 01/12/21, J.B., commenting out \clearpage
%\clearpage

% Description of performance profiles
% , II and 
The \jbr{outcomes from the subsequent Experiments I--III} are summarized in Figures  \ref{fig:EX_I}--\ref{fig:EX_III} as performance profiles (Dolan and Mor\'{e} \cite{DolanMore02}, extended in \cite{MahajanLeyfferKirches11} and often used
to compare the effectiveness of various solvers). \jbr{Detailed information for each
problem instance is in Tables \ref{tab:exp1}--\ref{tab:exp3}.}
Relative performances are displayed in terms of iterations and computation times. 
The performance metric $\rho_s(\tau)$ on $n_p$ test problems is given by
\begin{equation*}
		\rho_s(\tau) = \frac{\text{card}\left\{ p : \pi_{p,s} \le \tau \right\}}{n_p} \quad \text{and} \quad \pi_{p,s} = \frac{t_{p,s}}{ \underset{1\le i \le S,\ i \ne s}{\text{ min } t_{p,i}} },
\end{equation*} 
where $ t_{p,s}$ is the ``output'' (i.e., iterations or time) of
``solver'' $s$ on problem $p$, and $S$ denotes the total number of solvers for a given comparison. This metric measures
the proportion of how close a given solver is to the best result. Extended performance profiles are the same as the classical ones but include the part of 
the domain where $\tau \le 1$.
In the profiles we include a dashed vertical grey line to indicate $ \tau = 1 $.
 \jbbo{We note that although the iteration numbers are recorded 
differently for each solver, they correspond approximately to the number of KKT systems solved.}

\jbbo{Overall, we observe that the number of iterations used by the respective
solvers is relatively similar across different problems. However, the differences
in computation times are large. In particular, the RCR implementations use the least time in almost
all problem instances. This is possible because RCR
enables an efficient decoupling of computations with the constraint matrix
$ \b{A} $ and remaining small terms.}

\subsection{Experiment I}
\label{sec:numex1}
\begin{figure*}[t!]  % Figure 1
	\begin{minipage}{0.48\textwidth}
		\includegraphics[trim=0 0 20 0,clip,width=\textwidth]{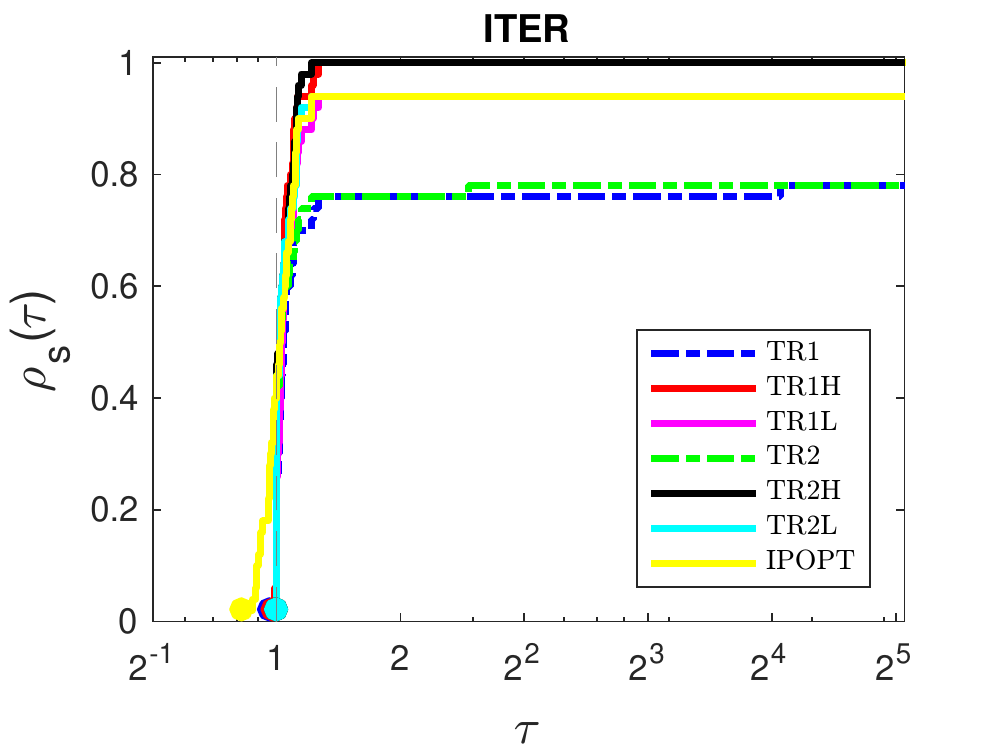}
	\end{minipage}
		\hfill
	\begin{minipage}{0.48\textwidth}
		\includegraphics[trim=0 0 20 0,clip,width=\textwidth]{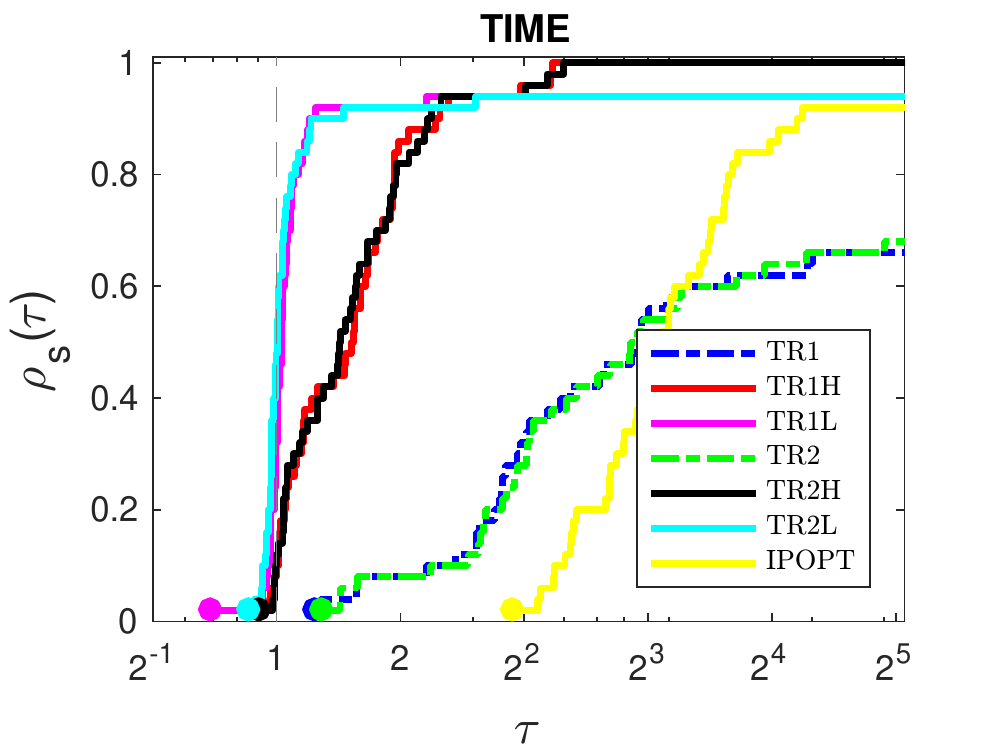}
	\end{minipage}
	\caption{\jbbo{Comparison of the 7 solvers from Experiment I using performance profiles  \cite{DolanMore02} on 50
	test problems from \cite{SuiteSparseMAT}.}
	 \jbbo{\TR{2H} and \TR{1H} converge on all problem instances (100\%).
	\TR{2L}, \TR{1L} and {\small IPOPT} converge on 47 problems (94\%).
	\TR{2} and \TR{1} are not applied to 9 large problems. 
	In the right plot, \TR{2L} and \TR{1L} are the fastest
	\jbr{(as seen from their curves being above others)}, while
	\TR{2H} and \TR{1H} are the most robust \jbr{(as seen from their curves ultimately
	reaching the top of the plot)}.
	Overall, \TR{2\{H,L\}} and \TR{1\{H,L\}} are faster than the other solvers.}}
		\label{fig:EX_I}       
\end{figure*}

This experiment uses problems with sparse and possibly low-rank $ \b{A} \in \mathbb{R}^{m \times n} $. The objective is the Rosenbrock function
\begin{equation*}
	f(\b{x}) = \sum_{i=1}^{n/2} (\b{x}_{2i}-\b{x}_{2i-1})^2 + (1-\b{x}_{2i-1})^2,
\end{equation*}
where $ n $ is an even integer. The matrices $ \b{A} \in \mathbb{R}^{m \times n} $ are obtained from the Suite\-Sparse Matrix Collection \cite{SuiteSparseMatrix}. Because \TR1 and \TR2
were not developed for problems with a large number of constraints, these solvers are only applied to problems for which $ m \le 2500 $. All other solvers were run on all test problems. Convergence of an algorithm is determined when two conditions are satisfied:
\begin{equation} \label{eq:converge}
    \norm{\b{P} \bk{g}}_{\infty} < 10^{-5} \quad \text{and} \quad
    \norm{\b{A} \bk{x} - \b{b}}_2 < 10^{-7}.
\end{equation}
 We summarize the outcomes in Figure\jbr{~\ref{fig:EX_I}} and Table \ref{tab:exp1}.

In this experiment we observe that our proposed algorithms (any of \TR{1\{H,L\}}, \TR{2\{H,L\}}) 
perform well in terms of computing time. Both ``H" versions of the proposed algorithms converged to the prescribed tolerances on all problems.
On the other hand, the ``L" versions are often the overall fastest, yet they did 
not converge on 3 problem instances (\texttt{beacxc}, \texttt{lp\_cre\_d}, \texttt{fit2d}).

\jbk{After rerunning the 3 problems for which {\small IPOPT } did not converge, we note that {\small IPOPT } did converge to its own (scaled) tolerances on one of these problems (\texttt{beacxc}), yet the computed solution did not satisfy 
\eqref{eq:converge}.
%the specified tests
%$\| Pg_k \|_{\infty} \le 10^{-5} $ \& $\| A\bk{x}-b \|_2 \le 10^{-7}$.
On the other two problems (\texttt{lp\_cre\_d}, \texttt{fit2d}), {\small IPOPT } returned a message such as info.status=$-2$, which is caused by an abort when the ``restoration phase" is called at an almost feasible point.}

\subsection{Experiment II}
\label{sec:EX2}
\begin{figure*}[t!]  % Figure 2
	\begin{minipage}{0.48\textwidth}
		\includegraphics[trim=0 0 20 0,clip,width=\textwidth]{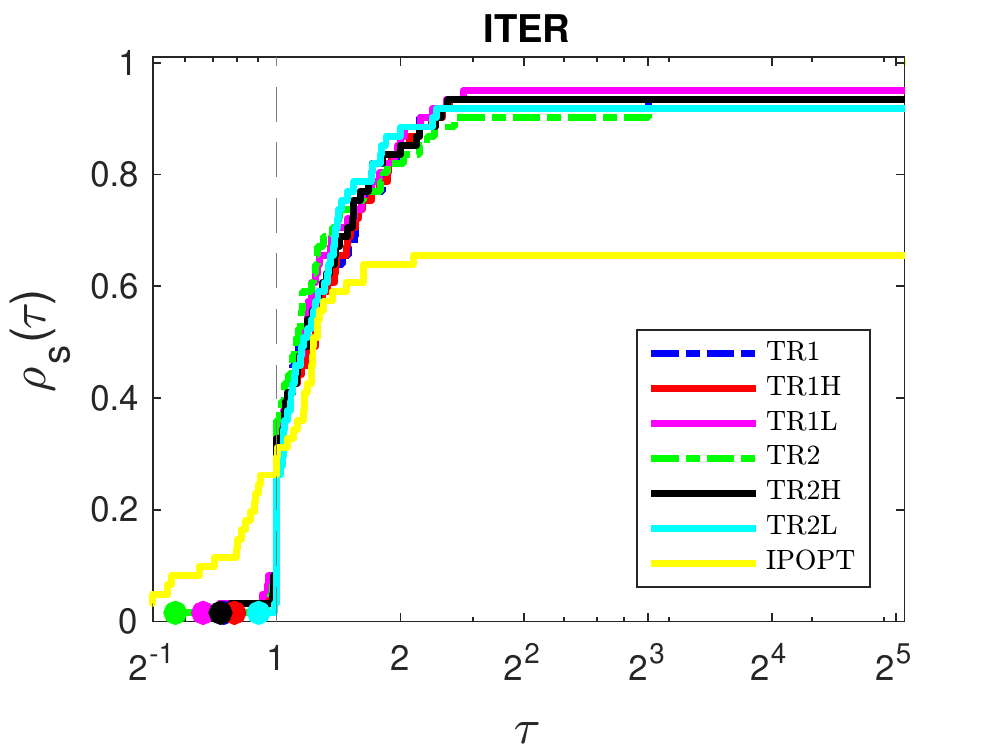}
	\end{minipage}
		\hfill
	\begin{minipage}{0.48\textwidth}
		\includegraphics[trim=0 0 20 0,clip,width=\textwidth]{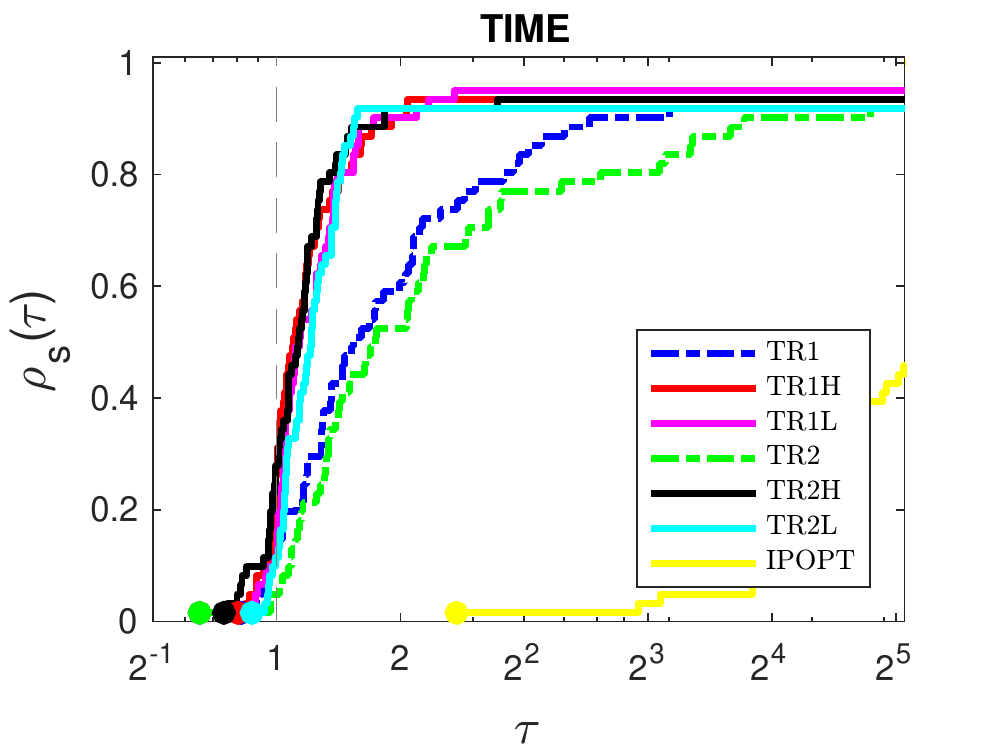}
	\end{minipage}
	\caption{\jbbo{Comparison of the 7 solvers from Experiment II using performance profiles %\cite{DolanMore02}
	on 62 test problems from \cite{GouOT03}.}
	 \jbbo{\TR{1L} converged on 58 problems. %, the largest number of problems. 
	All other solvers except {\small IPOPT} converged on 57 problems.
	In the left plot, the iteration numbers for \TR1, \TR{1\{H,L\}},
	\TR2 and \TR{2\{H,L\}} are similar, as seen by the tight clustering
	of the lines. However, the computational times of \TR1 and \TR2 are
	markedly higher than those of \TR{1\{H,L\}} and \TR{2\{H,L\}},
	as seen from the widening gap in the right plot.
	} }
		\label{fig:EX_II}       
\end{figure*}

%\jbb{[Section last edited before 11/19/20]}
In a second experiment, we compare the 7 solvers on large problems from the CUTEst collection
\cite{GouOT03}. The dimension $n$ is determined by the size of the corresponding CUTEst problem, % m=\text{ceil}(0.25n)
while we set $m$ to be about $25\%$ of $n$, i.e., $ \texttt{m=ceil(0.25n)} $.  The
matrices $\b{A}$ are formed as $ \texttt{A=sprand(m,n,0.1)} $, with $ \texttt{rng(090317)} $. Convergence is determined
by each algorithm internally. For \TR1, \TR{1H}, \TR{1L}, \TR2, \TR{2H}, \TR{2L} the conditions
$ \| \b{P} \bk{g} \|_{\infty} < 1 \times 10^{-5} $ and $ \| \b{A} \bk{x} - \b{b} \|_2 < 5 \times 10^{-8} $ are explicitly enforced, while for {\small IPOPT } we set $ \texttt{options\_ipopt.ipopt.tol=1e-5} $. We use the iteration limit of $100,000$ for all solvers. The limited-memory parameter is $l=5$
for all TR solvers and $l=6$ (default) for {\small IPOPT }. We summarize the outcomes in 
\jbr{Figure~\ref{fig:EX_II}} and Table \ref{tab:exp2}.
%\jbb{and include performance profiles in Appendix D}.
%\Red{and Figure~\ref{fig:EX_II}.}

% Figure 2 here

\subsection{Experiment III}
\label{sec:EX3}
\begin{figure*}[t!]  % Figure 3
	\begin{minipage}{0.48\textwidth}
		\includegraphics[trim=0 0 20 0,clip,width=\textwidth]{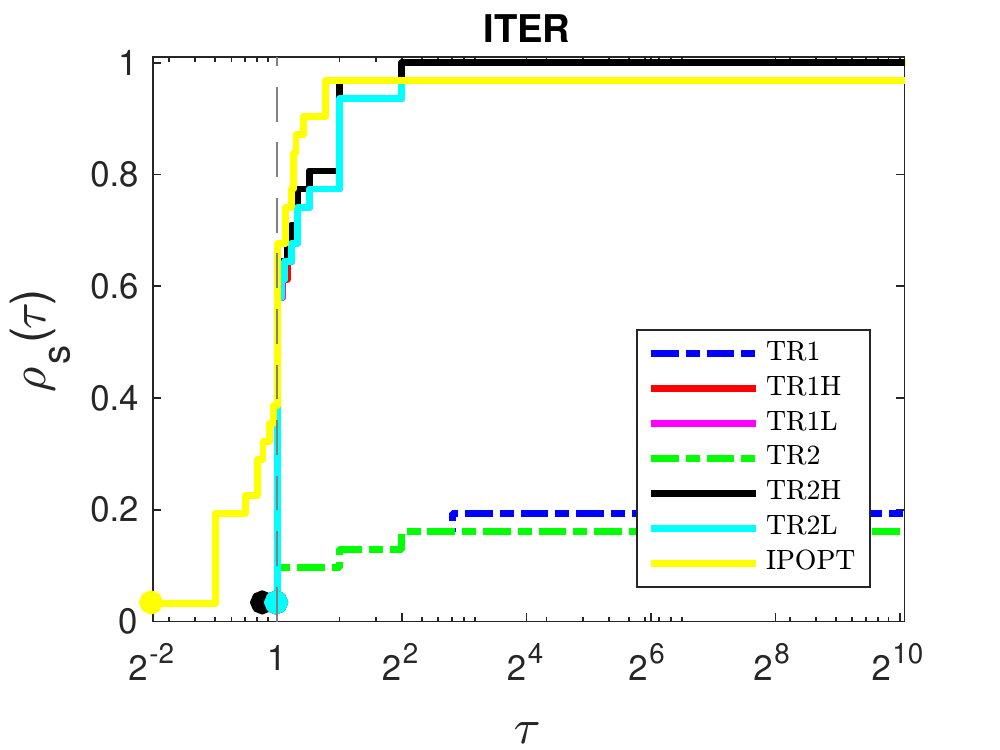}
	\end{minipage}
		\hfill
	\begin{minipage}{0.48\textwidth}
		\includegraphics[trim=0 0 20 0,clip,width=\textwidth]{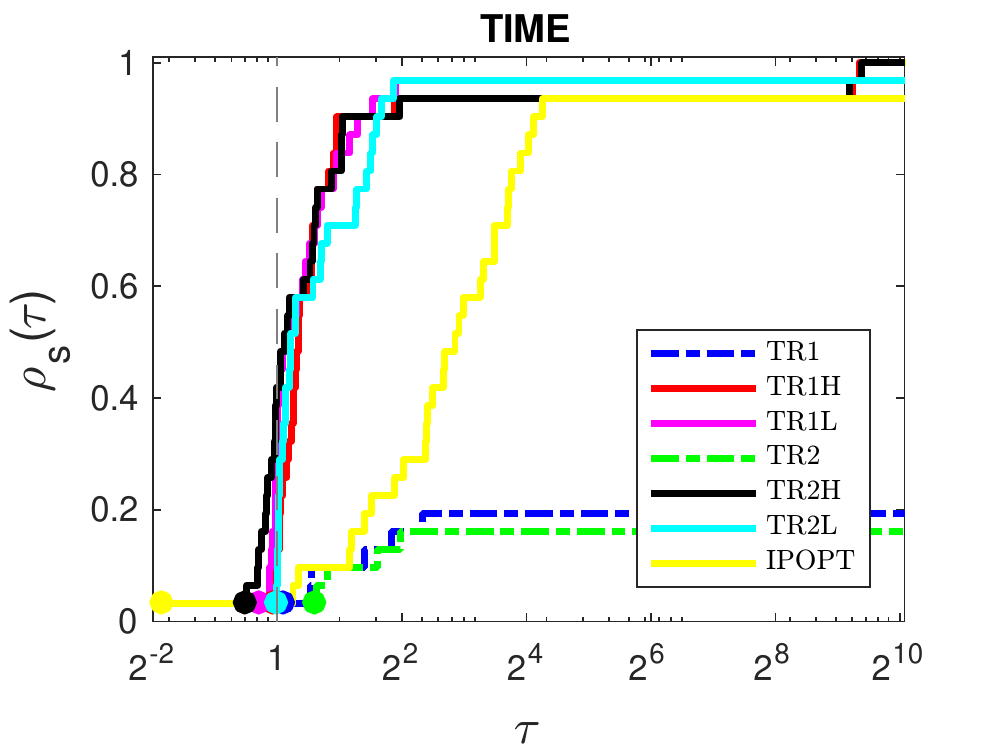}
	\end{minipage}
	\caption{\jbr{Comparison of the 7 solvers from Experiment III using performance profiles %\cite{DolanMore02}
	on 31 large linear equality constrained test problems from \cite{GouOT03}.
	\TR1 and \TR2 are applied to 6 problems
	(they are not practical on the remaining problems because of their size). \TR{2H}
	(also \TR{1H}) converged on all 31 instances. \TR{1L}, \TR{2L}, and {\small IPOPT } converged
	on 30 problems. In the {\small ITER} plot the number of iterations is relatively
	similar across the solvers that converged. In the {\small TIME } plot
	there is a gap between \TR1\{H,L\},\TR2\{H,L\} and {\small IPOPT}. \TR{2L} can
	have computational advantages, but appears slightly less robust than
	\TR{2H}, as seen from the final staircase in the {\small TIME } plot.}
% 	 \jbbo{\TR{1L} converged on 58 problems. %, the largest number of problems. 
% 	All other solvers except {\small IPOPT} converged on 57 problems.
% 	In the left plot, the iteration numbers for \TR1, \TR{1\{H,L\}},
% 	\TR2 and \TR{2\{H,L\}} are similar, as seen by the tight clustering
% 	of the lines. However, the computational times of \TR1 and \TR2 are
% 	markedly higher than those of \TR{1\{H,L\}} and \TR{2\{H,L\}},
% 	as seen from the widening gap in the right plot.}
	 }
		\label{fig:EX_III}       
\end{figure*}

%\jbb{[Section last edited before 11/19/20]}
\jbr{In a third experiment we compare the 7 solvers on 31 linear equality constrained problems
from CUTEst. Four of these problems (\texttt{AUG2D}, \texttt{AUG2DC}, \texttt{AUG3D}, \texttt{AUG3DC}) directly correspond
to the problem formulation \eqref{eq:main}. The remaining problems have additional
bound constraints, which are relaxed in this experiment. %In particular,
Problems 1--19 in Table \ref{tab:exp3} are convex and can immediately be attempted by 
the solvers (with bounds released). Problems 20--31 are not convex when the bounds
are relaxed, 
%so that
but adding the term $ \frac{\delta}{2} \| \b{x} \|_2^2 $ with $ \delta = 10 $
to the objective functions
%, which is numerically observed to 
produced finite solutions for these problems.
% In a third experiment, we compare the 7 solvers on large problems from the CUTEst collection
% \cite{GouOT03}. The dimension $n$ is determined by the size of the corresponding CUTEst problem, % m=\text{ceil}(0.25n)
% while we set $m$ to be about $25\%$ of $n$, i.e., $ \texttt{m=ceil(0.25n)} $.  The
% matrices $\b{A}$ are formed as $ \texttt{A=sprand(m,n,0.1)} $, with $ \texttt{rng(090317)} $. Convergence is determined
% by each algorithm internally. For \TR1, \TR{1H}, \TR{1L}, \TR2, \TR{2H}, \TR{2L} the conditions
% $ \| \b{P} \bk{g} \|_{\infty} < 1 \times 10^{-5} $ and $ \| \b{A} \bk{x} - \b{b} \|_2 < 5 \times 10^{-8} $ are explicitly enforced, while for IPOPT we set $ \texttt{options\_ipopt.ipopt.tol=1e-5} $ and
% $\texttt{options\_ipopt.dual\_inf\_tol=1e-5} $. We use the default limit on iterations, which is 
% $3000$ for IPOPT and $100,000$ for the other solvers. The limited-memory parameter is $l=5$
% for all solvers. We summarize the outcomes in Table \ref{tab:exp3}. 
As in the previous experiment, convergence is determined
by each algorithm internally. For \TR1, \TR{1H}, \TR{1L}, \TR2, \TR{2H}, \TR{2L} the conditions
$ \| \b{P} \bk{g} \|_{\infty} < 1 \times 10^{-5} $ and $ \| \b{A} \bk{x} - \b{b} \|_2 < 5 \times 10^{-8} $ are explicitly enforced, while for {\small IPOPT} we set $ \texttt{options\_ipopt.ipopt.tol=1e-5} $. We use the iteration limit of $100,000$ for all solvers. The limited-memory parameter is $l=5$
for all TR solvers and $l=6$ (default) for {\small IPOPT}. Since \TR1 and \TR2 are not 
designed for large $m$, they are applied to problems with $ m < 2500 $,
with the exception of 3 problems (\texttt{BLOWEYA}, \texttt{BLOWEYB}, \texttt{BLOWEYC})
that did not terminate within hours using \TR1 and \TR2. All other solvers are applied
to all problems. The results
%experiment outcomes
are in Figure \ref{fig:EX_III} and
Table \ref{tab:exp3}}. 

% Here Figure 3

\section{Conclusion}
\label{sec:conclusion}
%\jbb{[Section last edited before 11/19/20]}
For subproblem \eqref{eq:trsub},
this article develops the reduced compact representation~(RCR) of the (1,1)
block in the inverse KKT matrix, when the objective Hessian %of 2$^{\text{nd}}$ derivatives
is approximated by a compact quasi-Newton matrix. 
The representation is based on the \jbr{fact} that part of the solution
to the KKT system is unaffected when it is projected onto the nullspace
of the constraints. An advantage of the RCR is that it enables
\jbbo{a decoupling of solves with the constraint matrix and remaining small terms.
Moreover,}
a projected gradient can be used in two places: % economic updates because 
once as part of the matrix update, and second as part of the new step.
By effectively handling
orthogonal projections,
in combination with limited memory techniques, \jbbo{we can compute search directions efficiently}.
%we have developed computationally efficient search directions for large-scale problems. 
\jbk{We apply the orthogonal projections with a sparse
QR factorization or a preconditioned LSQR iteration}, \jbbo{including 
large and potentially rank-deficient constraints}.
% {This enables us to 
% handle large and potentially rank-deficient constraints effectively.}
The RCRs are implemented in two trust-region algorithms, one of which
exploits the underlying matrix structures in order to compute the search
direction by an analytic formula. 
\jbk{The other is based on an $\ell_2$ norm and uses the RCR 
within a 1D Newton iteration to determine
the optimal scalar shift.}
In numerical experiments on large problems,
our % proposed
implementations of the RCR yield
often significant improvements in the computation time,
\jbk{as a result of the advantageous structure of the proposed matrices}.
%to solve a problem.

\Red{Applications of problem \eqref{eq:main} often include bounds $\ell \le x \le u$. When second derivatives of the objective function are available, the problem is best handled by an interior method.  Otherwise, a barrier function could be added to the objective, and the methods here may sometimes be effective on a sequence of large equality-constrained subproblems.}

% \section{Notes}
% \begin{itemize}
% 	\item In the email from March, Johannes asks if Michael Saunders is interested in extending
% 	the methods from \cite{BMP19} (COAP 2019) into a new work. The initial starting point 
% 	for the new work was to use Section 8 of the manuscript LLNL-JRNL-75523
% 	(Section 8 was not published in COAP 2019 and is titled ``\emph{Extensions to Higher-Dimensional Linear Constraints}").
% 	\item The approaches in Section \ref{sec:motiv} of this document seem a bit more straight-forward than
% 	the ones from Section 8 in LLNL-JRNL-75523, with the same goal of extending
% 	(COAP 2019) into a new work. 
% 	\item This new work is much at the beginning, and it can jointly be developed.
% 	\item Questions: 
% 		\begin{itemize}
% 			 \item Do you think our proposed contributions could justify publication?
% 			 \item If not, what could be considered to justify publication?
% 			 \item How best to solve with $ \b{A} \b{A}\tp $?
% 		 \end{itemize}		 
% \end{itemize}

%\clearpage
%\newpage

\section*{Appendix A}
\label{sec:appA}
Here we describe a simplified expression for the matrix $ \bk{C} \tp \bk{G} \bk{C} $
from section~\ref{sec:redComp}. Recall that the L-BFGS inverse 
$ 
    \bik{B} = \delta_k I + \bk{J} \bk{W} \bk{J} \tp
$
is defined by
\begin{equation*}
    \bk{J}  = \bmat{ \bk{S} & \bk{Y} }, \quad
    \bk{W}  =
    \bmat{  \bk{T}^{-\top}(\bk{D} + \delta_k \bk{Y} \tp \bk{Y})\bik{T}  & -\delta_k \bk{T}^{-\top} \\
            -\delta_k\bik{T} & 0_{l \times l}}.
\end{equation*}
First, \jbr{note} that
\begin{equation*}
    \bk{C} \equiv \b{A} \bk{J} \bk{W} 
    = \bmat{ 0 & \b{A} \bk{Y} } \bk{W} 
    = \bmat{ -\delta_k \b{A} \bk{Y} \bik{T} & 0 }.
\end{equation*}
Second, \jbr{it holds} that 
\begin{equation*}
    \bik{G} \equiv \b{A} \bik{B} \b{A} \tp 
    = \delta_k \b{A} \b{A}\tp + \b{A} \bk{J} \bk{W} \bk{J} \tp \b{A} \tp
    = \delta_k \b{A} \b{A}\tp + \bk{C}
    \bmat{ 0 \\  (\b{A}\bk{Y}) \tp },
\end{equation*}
so that $ \bik{G} = \delta_k \b{A} \b{A} \tp $, because the last term in 
the above expression for $ \bik{G} $ vanishes. Multiplying $ \bk{C} \tp $, $ \bk{G} $ and
$ \bk{C} $ we see that
\begin{equation*}
    \bk{C} \tp \bk{G} \bk{C} = 
    \bmat{\delta_k 
     \bk{T}^{-\top} \bk{Y} \tp \b{A} \tp (\b{A} \b{A} \tp )^{-1} \b{A} \bk{Y} \bik{T} & 0_{l \times l}\\
     0_{l \times l} & 0_{l \times l} }.
\end{equation*}

\section*{Appendix B}
\label{sec:appB}
This appendix describes how we apply the functions from the SuiteSparse
library \cite{SuiteSparseMAT} in our implementations. We use
SuiteSparse version 5.8.1 from
\url{https://github.com/DrTimothyAldenDavis/SuiteSparse/releases}.

\subsection*{B.1: Householder QR projection}
\label{sec:appB1}
The Matlab commands to compute the projection $ P \bk{g} $ using a 
Householder QR factorization are listed in Table \ref{tab:HHQR}. %for completeness. %\newline

%\begin{center}
\begin{table}[ht]  % Table 3
    \caption{Matlab commands to use SparseSuite functions for computing projections $z=Py$
             using a Householder QR factorization.}
    \label{tab:HHQR}
    \qquad
    \begin{tabular}{l}
       \texttt{\% Options}
    \\ \texttt{opts.Q = `Householder';}
    \\ \texttt{opts.permutation = `vector';}
    \\
    \\ \texttt{\% QR factorization using SPQR}
    \\ \texttt{[Q,\~{},\~{},info] = spqr(A',opts);}
    \\ \texttt{rankA = info.rank\_A\_estimate;}
    \\
    \\ \texttt{\% Projection}
    \\ \texttt{ztmp = spqr\_qmult(Q,y,0);}  % g
    \\ \texttt{zrkA = zeros(rankA,1);}
    \\ \texttt{z = [zrkA;ztmp(rankA+1:end)];} % gp
    \\ \texttt{z = spqr\_qmult(Q,z,1);} % gp
    \end{tabular}
\end{table}
%\end{center}

\subsection*{B.2: Preconditioned LSQR projection}
\label{sec:appB2}
The Matlab commands to compute the projection $ P \bk{g} $ using
preconditioned LSQR \cite{PS82a} are listed in Table \ref{tab:LSQR}. %for completeness. %\newline

%\begin{center}
\begin{table}[ht]  % Table 4
	%\caption{Matlab commands for computing projections}
    \caption{Matlab commands for computing projections
             $z=Py$ \hbox{using} preconditioned LSQR
             (where  $P = I - A\tp (AA\tp)\inv A$).
             \Red{If $A$ has full row rank ($\texttt{rankA}=m$),
             LSQR should need only 1 iteration.}
             Notes: SPQR uses all of $\b{A}\tp$ in the QR factorization
             $\b{A}\tp \b{P}_{\text{msk}} = \b{Q}\b{R}$,
             where $\b{P}_{\text{msk}}$ is a column permutation of
             $\b{A}\tp$ and $\b{R}$ is upper trapezoidal. We store the 
             permutation in the vector $\texttt{maskA}$. If $\b{A}\tp$ does not have full row rank,
             we use the first $\texttt{rankA}$ columns
             of $\b{A}\tp\b{P}_{\text{msk}}$ 
             (the command $\texttt{A(maskA(1:rankA),:)'}$).
    	     If $A$ contains some relatively dense columns, we should partition
             $\b{A}\b{P}_{\text{prt}} = [\: \b{A}_{S} \: \b{A}_{D} \:]$ into sparse and dense columns,
             then use $\b{A}_{S}$ in place of $\b{A}$
             in the call to $\texttt{spqr}$.
             }
    \label{tab:LSQR}
    \qquad
    \begin{tabular}{l}
       \texttt{\% Options}
    \\ \texttt{opts.econ = 0;}
    \\ \texttt{opts.Q = `Householder';}
    \\ \texttt{opts.permutation = `vector';}
    \\ \texttt{tol = 1e-15;}
    \\ \texttt{maxit = m;}
    \\
    \\ \texttt{\% Preconditioner using a triangular}
    \\ \texttt{\% factor from SPQR}
    \\ \texttt{[\~{},R,maskA,info] = spqr(A',opts);}
    \\ \texttt{rankA = info.rank\_A\_estimate;}
    %\\ \texttt{maskA = maskA1(1:rankA); }
    \\
    \\ \texttt{\% Projection} % [x,\~{},\~{},\~{}]
    \\ \texttt{x = lsqr(A(maskA(1:rankA),:)',y,...} % g
    \\ \quad \quad \quad \quad \ \ \texttt{tol,maxit,R(1:rankA,1:rankA));} %[]
    \\ \texttt{z = y - A(maskA(\jbbo{1:rankA}),:)'*x(1:rankA,1);} % gp, g
    % \texttt{ztmp = spqr\_qmult(QA,g,0);} \\ 
    % \texttt{zrkA = zeros(rankA,1);} \\
    % \texttt{gp = [zrkA;ztmp(rankA+1:end)];} \\
    % \texttt{gp = spqr\_qmult(QA,gp,1);}
    \end{tabular}
\end{table}
%\end{center}

%\clearpage

\section*{Appendix C}
\label{sec:appC}
This appendix overviews the subproblem solution with the shape-changing norm.
Note that $ \b{U} = \bmat{ \b{Q}_1 & \b{U}_2 & \b{U}_3 } \in \mathbb{R}^{n \times n} $
(from section \ref{sec:eigen}) represents an orthogonal matrix, and that the quadratic function is
\begin{equation*}
	q(\b{s}) = \b{s} \tp \bk{g} + \frac{1}{2} \b{s} \tp \bk{B} \b{s} 
	= \b{s} \tp \b{U} \b{U} \tp \bk{g} + \frac{1}{2} \b{s} \tp \b{U} \b{U} \tp \bk{B} \b{U} \b{U} \tp \b{s}.
\end{equation*}
We introduce the change of variables 
$ \b{v} \tp =  \bmat{ \b{v}_1\tp & \b{v}_2\tp & \b{v}_3\tp } \equiv \b{s}\tp \b{U} $. Moreover, it
holds that 
$$
	\b{U} \tp \bk{B} \b{U} = 
		\bmat{
			\b{Q}_1 \tp \bk{B} \b{Q}_1	 &  \b{Q}_1 \tp \bk{B} \b{U}_2 & \b{Q}_1 \tp \bk{B} \b{U}_3 \\
			\b{U}_2\tp \bk{B} \b{Q}_1 & (\delta_k \b{I} + \b{\Lambda}_2)^{-1} & \\
			\b{U}_3\tp \bk{B} \b{Q}_1 & & \delta_k^{-1} \b{I}}
$$
% 	\bmat{
% 			\vline	 & \rule[.5ex]{2.5em}{0.4pt}\b{Q}_1 \tp \bk{B} \b{U}\rule[.5ex]{2.5em}{0.4pt} & \rule[.5ex]{1.5em}{0.4pt} \\
% 			\b{U}\tp \bk{B} \b{Q}_1 & (\delta_k \b{I} + \b{\Lambda}_2)^{-1} & \\
% 			\vline & & \delta_k^{-1} \b{I}
% 	}
%$
% $ 
% 	\b{U} \tp \bk{B} \b{U} = 
% 		\bmat{
% 			\b{Q}_1 \tp \bk{B} \b{Q}_1	 & \multicolumn{2}{c}{\rule[.5ex]{2.5em}{0.4pt}\b{Q}_1 \tp \bk{B} \b{U}\rule[.5ex]{2.5em}{0.4pt}} \\
% 			\b{U}\tp \bk{B} \b{Q}_1 & (\delta_k \b{I} + \b{\Lambda}_2)^{-1} & \\
% 			\vline & & \delta_k^{-1} \b{I}
% 	}
% % 	\bmat{
% % 			\vline	 & \rule[.5ex]{2.5em}{0.4pt}\b{Q}_1 \tp \bk{B} \b{U}\rule[.5ex]{2.5em}{0.4pt} & \rule[.5ex]{1.5em}{0.4pt} \\
% % 			\b{U}\tp \bk{B} \b{Q}_1 & (\delta_k \b{I} + \b{\Lambda}_2)^{-1} & \\
% % 			\vline & & \delta_k^{-1} \b{I}
% % 	}
% $
(cf.\ \cite[Lemma 2]{BMP19}), and that 
$$ 
\b{A} \b{U} \b{U} \tp \b{s} = \b{A} \b{U} \b{v} 
	= \bmat{\b{R} & \b{0} & \b{0}} \bmat{\b{v}_1 \\ \b{v}_2 \\ \b{v}_3}
    = Rv_1.
$$
% from \eqref{eq:trsub}
With the constraint $ \b{A} \b{s} = \b{0} = \b{A} \b{U} \b{v} $, this implies $ \b{v}_1 = \b{0} $
(for $\b{R}$ nonsingular).
Therefore, the trust-region subproblem defined by the
shape-changing norm decouples into a problem with $ \b{v}_2 $ and $ \b{v}_3 $
only (once $ \b{v}_1 = \b{0} $ is fixed):
\begin{align*}
    \underset{ 
        \tiny
        \begin{array}{c}
         \| \b{s} \|_U \le \Delta_k  \\
         \b{A} \b{s} = \b{0} 
    \end{array} }{ \text{ minimize } } q(\b{s})
     =
        \bigg\{&
        \underset{ 
         \| \b{v}_2 \|_{\infty} \le \Delta_k
        }{ \text{ minimize } } \biidx{2}{v} \tp \biidx{2}{U} \tp \bk{g} +
            \frac{1}{2} \biidx{2}{v} \tp (\delta_k \b{I} + \Lambda_2)^{-1} \biidx{2}{v} \\
         & + 
        \underset{ 
         \| \b{v}_3 \|_{2} \le \Delta_k  
        }{ \text{ minimize } } \biidx{3}{v} \tp \biidx{3}{U} \tp \bk{g} +
            \frac{\| \biidx{3}{v} \|^2_2}{2\delta_k} 
        \bigg\}.
\end{align*}
This reformulated subproblem can be solved analytically and the component-wise
solution of $ \biidx{2}{v} $ is in \eqref{eq:SCanalytic1}. The analytic solution
of $ \biidx{3}{v} $ is $ \biidx{3}{v} = \beta \biidx{3}{U} \tp \bk{g} $
with $ \beta $ from \eqref{eq:SCanalytic2}. Subsequently, $\b{s}$ is
obtained by transforming variables as 
$ \b{s} = \b{U} \b{v} = \biidx{2}{U} \biidx{2}{v} + \biidx{3}{U} \biidx{3}{v}  $.
The orthonormal matrix $ \biidx{2}{U} $ is computed as 
$ \biidx{2}{U} = \bmat{ \bk{S} & \bk{Z} } \bh{R}_2^{-1} \bh{P}_2 $,
and since $ \biidx{3}{U} \biidx{3}{U} \tp = \b{P} - \biidx{2}{U} \biidx{2}{U} \tp $, the optimal step with the shape-changing norm is as in \eqref{sec:SCNorm}:
\begin{equation*}
    \b{s}_{SC} = \biidx{2}{U}(\biidx{2}{v} - \beta \biidx{2}{U}\tp \bk{g})
    + \beta \b{P} \bk{g}.
\end{equation*}
With $ \bk{u} \equiv \biidx{2}{U} \tp \bk{g} $, the step is then
computed as in Algorithm \ref{alg:LTRSC-LEC} (line 15):
\begin{equation*}
    \b{s}_{SC} = \bmat{ \bk{S} & \bk{Z} } \bh{R}_2^{-1} \bh{P}_2(\biidx{2}{v} - \beta \bk{u})
    + \beta \b{P} \bk{g}.
\end{equation*}

\section*{Appendix D}
\label{sec:appD}
% This section to hold the tables

\subsection*{D.1: \jbr{Detailed Table for Experiment I}}
\label{sec:appD1}

\begin{table}[p]  % Table 1
 \captionsetup{width=1\linewidth}
\caption{Experiment I compares 7 solvers on problems from the SuiteSparse Matrix Collection \cite{SuiteSparseMatrix}. 
%``Dty" in column 3 is the density of a particular matrix $ \b{A} $ calculated as $ \text{Dty} = \frac{\text{nnz}(\b{A})}{m \cdot n} $.
Entries with $\texttt{N/A}^{*}$ denote problems to which \TR1 and \TR2 were not applied, because they are too large.
$\texttt{NC}^{\dagger}$ means the solver did not converge to tolerances.  \jbbo{\TR{2H} and \TR{1H} converged on all problem
instances. Overall, the computational times of \TR{2\{H,L\}} 
and \TR{1\{H,L\}} were lower by a 
significant factor compared to the times of \TR1, \TR2, and {\small IPOPT}.
The number of iterations for each solver is similar across all problems.}}
\label{tab:exp1}

% Adjusting table column separation
\setlength{\tabcolsep}{2pt} % Default value: 6pt

\scriptsize
\hbox to 1.00\textwidth{\hss % 1.25\textwidth
\begin{tabular}{|l|c|c| c c | c c | c c | c c | c c | c c | c c |} % tabular
		\hline
		\multirow{2}{*}{Problem} & \multirow{2}{*}{$m$/$n$} & \multirow{2}{*}{$\text{rank}(\b{A})$} & % /\text{Dty}
		\multicolumn{2}{c|}{\jbbo{TR2}}  & 
		\multicolumn{2}{c|}{\jbbo{TR2H}} &  		
		\multicolumn{2}{c|}{\jbbo{TR2L}} &  	
		\multicolumn{2}{c|}{\jbbo{TR1}} &  	
		\multicolumn{2}{c|}{\jbbo{TR1H}} &  	
		\multicolumn{2}{c|}{\jbbo{TR1L}} &  	
		\multicolumn{2}{c|}{IPOPT} \\  	
		\cline{4-17}
		& & %&		
		& It &  Sec
		& It &  Sec 
		& It &  Sec 
		& It &  Sec 
		& It &  Sec 
		& It &  Sec 
		& It &  Sec \\
		\hline
%		\endhead
%		\hline \multicolumn{17}{r}{\textit{Continued on next page}} \\
%		\endfoot
%		\endlastfoot
$\texttt{beacxc}$ & 497/506 & 449/0.2 &73 &0.52&25 &\emph{0.044}&25 &0.15&419 &3.8&25 &\textbf{0.041}&25 &0.15&$\texttt{NC}^{\dagger}$ & $\texttt{NC}$ \\ 
$\texttt{lp\_25fv47}$ & 821/1876 & 820/0.007 &60 &0.82&60 &0.21&60 &\textbf{0.14}&62 &0.85&62 &0.22&62 &\emph{0.14}&61 &0.73\\ 
$\texttt{lp\_agg2}$ & 516/758 & 516/0.01 &40 &0.21&40 &\emph{0.054}&40 &\textbf{0.052}&42 &0.21&42 &0.056&42 &0.055&41 &0.22\\ 
$\texttt{lp\_agg3}$ & 516/758 & 516/0.01 &39 &0.21&39 &\textbf{0.051}&39 &\emph{0.051}&39 &0.2&39 &0.052&39 &0.051&44 &0.24\\ 
$\texttt{lp\_bnl1}$ & 643/1586 & 642/0.005 &70 &0.57&70 &0.14&70 &\emph{0.079}&67 &0.6&67 &0.14&67 &\textbf{0.078}&62 &0.59\\ 
$\texttt{lp\_bnl2}$ & 2324/4486 & 2324/0.001 &69 &11&69 &0.62&69 &\emph{0.28}&69 &11&69 &0.52&69 &\textbf{0.27}&67 &2.2\\ 
$\texttt{lp\_cre\_a}$ & 3516/7248 & 3428/0.0007 &$\texttt{N/A}^{*}$ & $\texttt{N/A}$ &83 &0.65&83 &\textbf{0.37}&$\texttt{N/A}^{*}$ & $\texttt{N/A}$ &88 &0.71&88 &\emph{0.38}&87 &3.3\\ 
$\texttt{lp\_cre\_d}$ & 8926/73948 & 6476/0.0004 &$\texttt{N/A}^{*}$ & $\texttt{N/A}$ &556 &1.2e+02&510 &\textbf{24}&$\texttt{N/A}^{*}$ & $\texttt{N/A}$ &503 &1e+02&552 &\emph{25}&$\texttt{NC}^{\dagger}$ & $\texttt{NC}$ \\ 
$\texttt{lp\_czprob}$ & 929/3562 & 929/0.003 &17 &0.27&17 &0.059&17 &\emph{0.032}&17 &0.25&17 &0.049&17 &\textbf{0.028}&18 &0.34\\ 
$\texttt{lp\_d6cube}$ & 415/6184 & 404/0.01 &35 &0.22&35 &0.4&35 &\textbf{0.17}&36 &0.21&36 &0.44&36 &\emph{0.18}&38 &1.1\\ 
$\texttt{lp\_degen3}$ & 1503/2604 & 1503/0.006 &39 &2.2&39 &0.25&39 &\emph{0.25}&39 &2.3&39 &0.27&39 &\textbf{0.24}&40 &1.5\\ 
$\texttt{lp\_dfl001}$ & 6071/12230 & 6071/0.0005 &$\texttt{N/A}^{*}$ & $\texttt{N/A}$ &226 &\textbf{16}&231 &19&$\texttt{N/A}^{*}$ & $\texttt{N/A}$ &226 &\emph{16}&238 &20&207 &1.2e+02\\ 
$\texttt{lp\_etamacro}$ & 400/816 & 400/0.008 &78 &0.27&78 &0.12&78 &\textbf{0.085}&86 &0.29&86 &0.13&86 &\emph{0.09}&68 &0.44\\ 
$\texttt{lp\_fffff800}$ & 524/1028 & 524/0.01 &$\texttt{NC}^{\dagger}$ & $\texttt{NC}$ &61 &\textbf{0.095}&$\texttt{NC}^{\dagger}$ & $\texttt{NC}$ &$\texttt{NC}^{\dagger}$ & $\texttt{NC}$ &59 &\emph{0.097}&$\texttt{NC}^{\dagger}$ & $\texttt{NC}$ &57 &0.45\\ 
$\texttt{lp\_finnis}$ & 497/1064 & 497/0.005 &150 &0.72&151 &0.2&156 &\textbf{0.13}&159 &0.69&155 &0.2&155 &\emph{0.14}&167 &1.2\\ 
$\texttt{lp\_fit2d}$ & 25/10524 & 25/0.5 &266 &1.3&266 &\emph{0.9}&258 &\textbf{0.88}&247 &1.4&261 &0.91&279 &0.99&$\texttt{NC}^{\dagger}$ & $\texttt{NC}$ \\ 
$\texttt{lp\_ganges}$ & 1309/1706 & 1309/0.003 &41 &1.3&41 &0.11&41 &\textbf{0.067}&41 &1.4&41 &0.12&41 &\emph{0.073}&37 &0.43\\ 
$\texttt{lp\_gfrd\_pnc}$ & 616/1160 & 616/0.003 &$\texttt{NC}^{\dagger}$ & $\texttt{NC}$ &54 &0.054&54 &\emph{0.043}&$\texttt{NC}^{\dagger}$ & $\texttt{NC}$ &54 &0.052&54 &\textbf{0.042}&48 &0.36\\ 
$\texttt{lp\_greenbea}$ & 2392/5598 & 2389/0.002 &149 &47&149 &1.2&149 &\textbf{0.63}&157 &33&153 &1.3&150 &\emph{0.72}&181 &6.8\\ 
$\texttt{lp\_greenbeb}$ & 2392/5598 & 2389/0.002 &149 &45&149 &1.2&149 &\emph{0.65}&157 &31&153 &1.3&150 &\textbf{0.65}&181 &6.5\\ 
$\texttt{lp\_grow22}$ & 440/946 & 440/0.02 &79 &0.24&79 &0.079&79 &\emph{0.071}&79 &0.24&79 &0.08&79 &\textbf{0.069}&65 &0.36\\ 
$\texttt{lp\_ken\_07}$ & 2426/3602 & 2426/0.001 &34 &12&34 &0.091&34 &\textbf{0.067}&34 &7.4&34 &0.093&34 &\emph{0.07}&31 &0.85\\ 
$\texttt{lp\_maros}$ & 846/1966 & 846/0.006 &74 &0.87&74 &\emph{0.21}&$\texttt{NC}^{\dagger}$ & $\texttt{NC}$ &74 &0.86&74 &\textbf{0.21}&$\texttt{NC}^{\dagger}$ & $\texttt{NC}$ &71 &0.92\\ 
$\texttt{lp\_maros\_r7}$ & 3136/9408 & 3136/0.005 &$\texttt{N/A}^{*}$ & $\texttt{N/A}$ &57 &\emph{2.1}&57 &2.2&$\texttt{N/A}^{*}$ & $\texttt{N/A}$ &57 &\textbf{2.1}&57 &2.3&51 &25\\ 
$\texttt{lp\_modszk1}$ & 687/1620 & 686/0.003 &71 &0.51&71 &0.13&71 &\textbf{0.07}&71 &0.51&71 &0.14&71 &\emph{0.071}&70 &0.53\\ 
$\texttt{lp\_osa\_30}$ & 4350/104374 & 4350/0.001 &$\texttt{N/A}^{*}$ & $\texttt{N/A}$ &46 &8.5&46 &\textbf{1.9}&$\texttt{N/A}^{*}$ & $\texttt{N/A}$ &45 &8.8&45 &\emph{2}&43 &30\\ 
$\texttt{lp\_osa\_60}$ & 10280/243246 & 10280/0.0006 &$\texttt{N/A}^{*}$ & $\texttt{N/A}$ &47 &24&47 &\textbf{5.9}&$\texttt{N/A}^{*}$ & $\texttt{N/A}$ &44 &23&44 &\emph{5.9}&42 &1.1e+02\\ 
$\texttt{lp\_pds\_02}$ & 2953/7716 & 2953/0.0007 &$\texttt{N/A}^{*}$ & $\texttt{N/A}$ &25 &0.25&25 &\emph{0.14}&$\texttt{N/A}^{*}$ & $\texttt{N/A}$ &25 &0.24&25 &\textbf{0.099}&26 &1.3\\ 
$\texttt{lp\_pds\_10}$ & 16558/49932 & 16558/0.0001 &$\texttt{N/A}^{*}$ & $\texttt{N/A}$ &61 &13&61 &\emph{8.1}&$\texttt{N/A}^{*}$ & $\texttt{N/A}$ &60 &13&60 &\textbf{7.9}&59 &62\\ 
$\texttt{lp\_perold}$ & 625/1506 & 625/0.007 &58 &0.39&58 &0.16&58 &\textbf{0.084}&58 &0.38&58 &0.16&58 &\emph{0.087}&57 &0.59\\ 
$\texttt{lp\_pilot}$ & 1441/4860 & 1441/0.006 &105 &13&105 &1.3&105 &\textbf{0.83}&109 &6.2&109 &1.3&109 &\emph{0.97}&117 &5.8\\ 
$\texttt{lp\_pilot87}$ & 2030/6680 & 2030/0.006 &102 &17&102 &2.6&102 &\textbf{1.9}&104 &15&104 &2.7&104 &\emph{1.9}&110 &15\\ 
$\texttt{lp\_pilot\_we}$ & 722/2928 & 722/0.004 &73 &0.69&73 &0.24&73 &\emph{0.11}&73 &0.65&73 &0.21&73 &\textbf{0.11}&81 &1.4\\ 
$\texttt{lp\_pilotnov}$ & 975/2446 & 975/0.006 &77 &1.3&77 &\textbf{0.35}&$\texttt{NC}^{\dagger}$ & $\texttt{NC}$ &77 &1.3&77 &\emph{0.35}&$\texttt{NC}^{\dagger}$ & $\texttt{NC}$ &78 &1.3\\ 
$\texttt{lp\_qap12}$ & 3192/8856 & 3192/0.001 &$\texttt{N/A}^{*}$ & $\texttt{N/A}$ &27 &\emph{3.3}&27 &3.7&$\texttt{N/A}^{*}$ & $\texttt{N/A}$ &26 &\textbf{3.1}&26 &3.6&25 &1.4e+02\\ 
$\texttt{lp\_qap8}$ & 912/1632 & 912/0.005 &20 &0.42&20 &0.15&20 &\emph{0.11}&22 &0.31&22 &0.15&22 &\textbf{0.1}&21 &1.9\\ 
$\texttt{lp\_scfxm1}$ & 330/600 & 330/0.01 &45 &0.1&45 &0.042&45 &\emph{0.036}&44 &0.098&44 &0.041&44 &\textbf{0.036}&44 &0.19\\ 
$\texttt{lp\_scfxm2}$ & 660/1200 & 660/0.007 &52 &0.42&52 &0.079&52 &\textbf{0.056}&57 &0.43&57 &0.094&57 &\emph{0.06}&55 &0.46\\ 
$\texttt{lp\_scfxm3}$ & 990/1800 & 990/0.005 &45 &0.8&45 &0.096&45 &\textbf{0.061}&45 &0.76&45 &0.097&45 &\emph{0.062}&48 &0.54\\ 
$\texttt{lp\_scsd1}$ & 77/760 & 77/0.04 &74 &0.035&74 &0.035&74 &0.044&74 &\textbf{0.033}&74 &\emph{0.034}&74 &0.043&$\texttt{NC}^{\dagger}$ & $\texttt{NC}$ \\ 
$\texttt{lp\_scsd6}$ & 147/1350 & 147/0.02 &84 &0.077&84 &\textbf{0.054}&84 &0.06&92 &0.084&92 &\emph{0.06}&92 &0.065&75 &0.34\\ 
$\texttt{lp\_scsd8}$ & 397/2750 & 397/0.008 &66 &0.21&66 &0.07&66 &\textbf{0.066}&65 &0.21&65 &0.069&65 &\emph{0.066}&66 &0.54\\ 
$\texttt{lp\_sctap1}$ & 300/660 & 300/0.009 &107 &0.25&107 &0.088&107 &\emph{0.075}&102 &0.24&102 &0.081&102 &\textbf{0.07}&100 &0.45\\ 
$\texttt{lp\_sctap2}$ & 1090/2500 & 1090/0.003 &145 &5.5&146 &0.43&146 &\textbf{0.18}&145 &3.6&143 &0.46&146 &\emph{0.19}&157 &2.3\\ 
$\texttt{lp\_sctap3}$ & 1480/3340 & 1480/0.002 &204 &27&205 &0.75&201 &\emph{0.39}&199 &13&197 &0.74&202 &\textbf{0.39}&220 &4.3\\ 
$\texttt{lp\_ship04l}$ & 402/2166 & 360/0.007 &84 &0.25&84 &0.12&84 &\textbf{0.077}&84 &0.27&84 &0.12&84 &\emph{0.08}&92 &0.84\\ 
$\texttt{lp\_ship04s}$ & 402/1506 & 360/0.007 &74 &0.17&74 &0.063&74 &\textbf{0.053}&74 &0.16&74 &0.065&74 &\emph{0.055}&71 &0.48\\ 
$\texttt{lp\_stair}$ & 356/614 & 356/0.02 &47 &0.11&47 &0.047&47 &\emph{0.046}&47 &0.11&47 &0.047&47 &\textbf{0.045}&47 &0.23\\ 
$\texttt{lp\_standata}$ & 359/1274 & 359/0.007 &78 &0.22&78 &0.072&78 &\emph{0.058}&79 &0.21&79 &0.067&79 &\textbf{0.057}&80 &0.65\\ 
$\texttt{lp\_standmps}$ & 467/1274 & 467/0.007 &52 &0.21&52 &0.06&52 &\textbf{0.042}&52 &0.21&52 &0.065&52 &\emph{0.043}&58 &0.48\\

 \hline
 \end{tabular}
 \hss}
\end{table} % Table 1

In this experiment the degree of difficulty in solving a problem depends largely
on handling $ \b{A} $, because the structure of the
objective function is the same for all instances.
We observe that our proposed algorithms (any of \TR{1\{H,L\}}, \TR{2\{H,L\}}) always use less
computation time (often significantly), except for two problem instances.
On problem $ \texttt{lp\_d6cube} $, \TR2 used less time than \TR{2H}, as did \TR1 over
\TR{1H}. However, the ``L" versions were fastest overall on this problem.
On problem $ \texttt{lp\_scsd1} $, \TR1 used the least time. In these
two problems the number of constraints is not large, and one can expect
that \TR1, \TR2 do comparatively well. %on such problems.
However, for all other 48 
problems the new methods used the least time. We observe that
both ``H" versions converged to the prescribed tolerances on all problems.
On the other hand, the ``L" versions are often the overall fastest, yet they did 
not converge on 3 problem instances (\texttt{beacxc}, \texttt{lp\_cre\_d}, \texttt{fit2d}). 

%\jbk{After rerunning the 3 problems for which IPOPT's column reports ``\texttt{NC}'', we find that IPOPT did converge to its own (scaled) tolerances on one of these problems (\texttt{beacxc}), yet the computed solution did not satisfy 
% \eqref{eq:converge}
% %the specified tests
% %$\| Pg_k \|_{\infty} \le 10^{-5} $ \& $\| A\bk{x}-b \|_2 \le 10^{-7}$.
% On the other two problems (\texttt{lp\_cre\_d}, \texttt{fit2d}), IPOPT returned a message such as info.status=$-2$, which is caused by an abort when the ``restoration phase" is called at an almost feasible point.}

\subsection*{D.2: \jbr{Detailed Table for Experiment II}}
\label{sec:appD2}

% Table 2 here

\begin{table}[p]  % Table 2
 \captionsetup{width=1\linewidth}
\caption{Experiment II compares 7 solvers on 61 large problems from the CUTEst collection \cite{GouOT03}. 
$\texttt{NC}^{\dagger}$ means the solver did not converge to tolerances. $\texttt{MX}^{\dagger}$ means the iteration limit was reached.
 \jbbo{\TR{1L} converged on 58 problems, the largest number of problems
amongst the solvers. \TR{2H} was faster than \TR2 on 51 problems,
and \TR{2L} was faster than \TR2 on 46 problems (the differences are often significant). \TR{1H} was faster than \TR{1} on 49 problems and
\TR{1L} was faster than \TR{1} on 41 problems (often significantly).
All of \TR{1\{H,L\}} and \TR{2\{H,L\}} were faster than
{\small IPOPT}.}}
\label{tab:exp2}

% Adjusting table column separation
\setlength{\tabcolsep}{1.5pt} % Default value: 6pt

\scriptsize
\hbox to 1.00\textwidth{\hss % 1.25\textwidth
\begin{tabular}{|l | c | c c | c c | c c | c c | c c | c c | c c|} % tabular
		\hline
		\multirow{2}{*}{Problem} & \multirow{2}{*}{$m$/$n$} & %\multirow{2}{*}{$\text{rank}(\b{A})$/\text{Dty}} & 
% 		\multicolumn{2}{c|}{TR1}  & 
% 		\multicolumn{2}{c|}{TR1x1} &  		
% 		\multicolumn{2}{c|}{TR1x2} &  	
% 		\multicolumn{2}{c|}{TR2} &  	
% 		\multicolumn{2}{c|}{TR2x1} &  	
% 		\multicolumn{2}{c|}{TR2x2} &
        \multicolumn{2}{c|}{\jbbo{TR2}}  & 
		\multicolumn{2}{c|}{\jbbo{TR2H}} &  		
		\multicolumn{2}{c|}{\jbbo{TR2L}} &  	
		\multicolumn{2}{c|}{\jbbo{TR1}} &  	
		\multicolumn{2}{c|}{\jbbo{TR1H}} &  	
		\multicolumn{2}{c|}{\jbbo{TR1L}} & 
		\multicolumn{2}{c|}{IPOPT} \\  	
		\cline{3-16}
		& %& %&		
		& It &  Sec
		& It &  Sec 
		& It &  Sec 
		& It &  Sec 
		& It &  Sec 
		& It &  Sec 
		& It &  Sec \\
		\hline
%		\endhead
%		\hline \multicolumn{17}{r}{\textit{Continued on next page}} \\
%		\endfoot
%		\endlastfoot
$\texttt{ARWHEAD}$ & 1250/5000 &343 &1.7e+02&349 &19&372 &19&264 &72&304 &\textbf{16}&315 &\emph{16}&$\texttt{NC}^{\dagger}$ & $\texttt{NC}$ \\ 
$\texttt{BDQRTIC}$ & 1250/5000 &181 &50&174 &\textbf{8.1}&187 &9.9&174 &31&186 &8.9&160 &\emph{8.4}&78 &1.2e+02\\ 
$\texttt{BOX}$ & 2500/10000 &240 &1.5e+03&280 &63&281 &79&218 &2.1e+02&258 &\textbf{54}&208 &\emph{58}&$\texttt{NC}^{\dagger}$ & $\texttt{NC}$ \\ 
$\texttt{BROYDN7D}$ & 1250/5000 &355 &20&370 &\emph{18}&367 &18&355 &20&370 &\textbf{17}&381 &19&432 &6.5e+02\\ 
$\texttt{BRYBND}$ & 1250/5000 &897 &1.5e+02&883 &\textbf{45}&1273 &64&1396 &1.2e+02&1177 &\emph{60}&1421 &70&1027 &1.7e+03\\ 
$\texttt{COSINE}$ & 2500/10000 &$\texttt{NC}^{\dagger}$ & $\texttt{NC}$ &5028 &1e+03&4527 &1.2e+03&4755 &2e+03&7318 &1.6e+03&3292 &\emph{910}&$\texttt{NC}^{\dagger}$ & $\texttt{NC}$ \\ 
$\texttt{CRAGGLVY}$ & 1250/5000 &373 &63&371 &\textbf{18}&369 &\emph{19}&400 &45&390 &20&397 &21&205 &3.4e+02\\ 
$\texttt{CURLY10}$ & 2500/10000 &1563 &7.2e+02&2498 &5.3e+02&1496 &\emph{429}&1512 &4.5e+02&1549 &\textbf{347}&1759 &4.9e+02&1775 &3e+04\\ 
$\texttt{CURLY20}$ & 2500/10000 &1951 &9.5e+02&2015 &\textbf{455}&1993 &\emph{552}&3149 &9.5e+02&4110 &8.7e+02&3836 &1.1e+03&$\texttt{NC}^{\dagger}$ & $\texttt{NC}$ \\ 
$\texttt{CURLY30}$ & 2500/10000 &4457 &2.8e+03&4210 &\emph{952}&3669 &1e+03&2744 &\textbf{783}&6940 &1.6e+03&6145 &1.7e+03&$\texttt{NC}^{\dagger}$ & $\texttt{NC}$ \\ 
$\texttt{DIXMAANA}$ & 750/3000 &10 &0.53&10 &\textbf{0.43}&10 &0.51&10 &0.5&10 &0.47&10 &\emph{0.46}&13 &8.3\\ 
$\texttt{DIXMAANB}$ & 750/3000 &9 &0.55&9 &0.5&9 &\emph{0.5}&9 &0.59&9 &\textbf{0.5}&9 &0.55&11 &8.1\\ 
$\texttt{DIXMAANC}$ & 750/3000 &12 &0.73&12 &0.67&12 &0.72&12 &\emph{0.65}&12 &\textbf{0.63}&12 &0.69&14 &10\\ 
$\texttt{DIXMAAND}$ & 750/3000 &23 &1.7&23 &1.1&23 &1.2&22 &\emph{0.93}&22 &\textbf{0.83}&22 &1&27 &16\\ 
$\texttt{DIXMAANE}$ & 750/3000 &35 &1.1&35 &1&35 &1.1&35 &\textbf{0.83}&35 &\emph{0.88}&35 &1.1&41 &18\\ 
$\texttt{DIXMAANF}$ & 750/3000 &183 &5.2&194 &\textbf{3.9}&194 &5.7&194 &6.6&195 &\emph{4.9}&203 &6.7&297 &1.3e+02\\ 
$\texttt{DIXMAANG}$ & 750/3000 &434 &19&397 &\textbf{8.3}&439 &12&435 &13&408 &\emph{9.8}&404 &11&$\texttt{NC}^{\dagger}$ & $\texttt{NC}$ \\ 
$\texttt{DIXMAANH}$ & 750/3000 &433 &14&470 &11&454 &13&459 &\emph{11}&421 &\textbf{9.3}&443 &12&422 &1.8e+02\\ 
$\texttt{DIXMAANI}$ & 750/3000 &82 &2&82 &\emph{1.8}&82 &2.4&82 &\textbf{1.6}&82 &1.8&82 &2.5&103 &46\\ 
$\texttt{DIXMAANJ}$ & 750/3000 &1054 &41&1506 &35&1023 &\emph{27}&1415 &42&1490 &34&944 &\textbf{24}&$\texttt{NC}^{\dagger}$ & $\texttt{NC}$ \\ 
$\texttt{DIXMAANK}$ & 750/3000 &2971 &1e+02&3026 &65&3082 &71&2831 &80&2870 &\textbf{61}&2691 &\emph{62}&$\texttt{NC}^{\dagger}$ & $\texttt{NC}$ \\ 
$\texttt{DIXMAANL}$ & 750/3000 &1461 &\textbf{38}&3198 &69&2609 &60&2690 &66&2728 &\emph{58}&2597 &59&$\texttt{NC}^{\dagger}$ & $\texttt{NC}$ \\ 
$\texttt{DIXON3DQ}$ & 2500/10000 &51 &17&51 &\emph{12}&51 &17&51 &17&51 &\textbf{12}&51 &16&56 &6.7e+02\\ 
$\texttt{DQDRTIC}$ & 1250/5000 &13 &1.7&7 &0.85&7 &0.77&13 &1.5&7 &\textbf{0.75}&7 &\emph{0.75}&7 &13\\ 
$\texttt{DQRTIC}$ & 1250/5000 &63 &\emph{4.6}&107 &6.7&107 &7&63 &\textbf{4.5}&107 &5.7&107 &6.1&93 &1.5e+02\\ 
$\texttt{EDENSCH}$ & 500/2000 &32 &\emph{0.33}&32 &0.4&32 &0.38&32 &\textbf{0.32}&32 &0.39&32 &0.36&34 &5\\ 
$\texttt{EG2}$ & 250/1000 &423 &2.2&504 &\emph{1.3}&439 &\textbf{1.3}&514 &5.2&624 &2&502 &1.9&908 &23\\ 
$\texttt{ENGVAL1}$ & 1250/5000 &31 &2.7&31 &\textbf{1.8}&31 &2&31 &2.6&31 &\emph{1.9}&31 &2&38 &61\\ 
$\texttt{EXTROSNB}$ & 250/1000 &148 &\emph{0.44}&148 &0.45&148 &0.49&145 &0.53&145 &0.46&145 &\textbf{0.39}&129 &3\\ 
$\texttt{FLETCHCR}$ & 250/1000 &150 &\emph{0.4}&150 &0.46&150 &0.51&150 &\textbf{0.37}&150 &0.42&150 &0.41&137 &3.1\\ 
$\texttt{FMINSRF2}$ & 1407/5625 &122 &10&122 &\emph{9.3}&122 &10&122 &10&122 &\textbf{7.8}&122 &9.6&167 &4e+02\\ 
$\texttt{FREUROTH}$ & 1250/5000 &287 &1e+02&247 &\textbf{12}&235 &\emph{13}&274 &37&255 &13&234 &13&202 &3.2e+02\\ 
$\texttt{GENHUMPS}$ & 1250/5000 &2215 &1.2e+02&1762 &99&1829 &\textbf{93}&2215 &1.3e+02&1762 &98&1829 &\emph{95}&$\texttt{NC}^{\dagger}$ & $\texttt{NC}$ \\ 
$\texttt{LIARWHD}$ & 1250/5000 &3854 &1.6e+03&3998 &4.4e+02&2726 &\emph{196}&2638 &1.2e+03&2408 &2.6e+02&1591 &\textbf{128}&$\texttt{NC}^{\dagger}$ & $\texttt{NC}$ \\ 
$\texttt{MOREBV}$ & 1250/5000 &151 &23&151 &22&151 &20&151 &19&151 &\emph{16}&151 &\textbf{16}&$\texttt{NC}^{\dagger}$ & $\texttt{NC}$ \\ 
$\texttt{MSQRTALS}$ & 256/1024 &$\texttt{MX}^{\dagger}$ & $\texttt{MX}$ &$\texttt{MX}^{\dagger}$ & $\texttt{MX}$ &$\texttt{MX}^{\dagger}$ & $\texttt{MX}$ &$\texttt{MX}^{\dagger}$ & $\texttt{MX}$ &78461 &6.6e+02&99724 &\emph{620}&$\texttt{NC}^{\dagger}$ & $\texttt{NC}$ \\ 
$\texttt{MSQRTBLS}$ & 256/1024 &$\texttt{MX}^{\dagger}$ & $\texttt{MX}$ &$\texttt{MX}^{\dagger}$ & $\texttt{MX}$ &$\texttt{MX}^{\dagger}$ & $\texttt{MX}$ &$\texttt{MX}^{\dagger}$ & $\texttt{MX}$ &$\texttt{MX}^{\dagger}$ & $\texttt{MX}$ &$\texttt{MX}^{\dagger}$ & $\texttt{MX}$ &$\texttt{NC}^{\dagger}$ & $\texttt{NC}$ \\ 
$\texttt{NCB20}$ & 1253/5010 &345 &47&348 &18&349 &18&314 &33&317 &\textbf{16}&307 &\emph{16}&252 &3.9e+02\\ 
$\texttt{NONCVXU2}$ & 1250/5000 &185 &20&185 &\textbf{9}&185 &9.5&186 &14&187 &\emph{9.2}&186 &9.4&120 &1.9e+02\\ 
$\texttt{NONCVXUN}$ & 1250/5000 &282 &33&283 &\textbf{14}&282 &\emph{14}&360 &31&354 &17&370 &19&199 &3.1e+02\\ 
$\texttt{NONDIA}$ & 1250/5000 &1612 &6.9e+02&1600 &88&1734 &\emph{88}&2764 &7.3e+02&1407 &\textbf{78}&1907 &98&$\texttt{NC}^{\dagger}$ & $\texttt{NC}$ \\ 
$\texttt{NONDQUAR}$ & 1250/5000 &897 &4.3e+02&865 &47&811 &\textbf{42}&816 &2.1e+02&876 &47&857 &\emph{44}&332 &8.1e+02\\ 
$\texttt{PENALTY1}$ & 250/1000 &8 &0.051&2 &0.018&2 &\emph{0.017}&8 &0.056&2 &0.019&2 &\textbf{0.016}&1 &0.043\\ 
$\texttt{POWELLSG}$ & 1250/5000 &88 &6.1&88 &\textbf{4.4}&88 &4.6&88 &5.6&88 &\emph{4.4}&88 &4.6&99 &1.5e+02\\ 
$\texttt{POWER}$ & 2500/10000 &51 &\textbf{17}&$\texttt{MX}^{\dagger}$ & $\texttt{MX}$ &$\texttt{MX}^{\dagger}$ & $\texttt{MX}$ &51 &\emph{17}&$\texttt{MX}^{\dagger}$ & $\texttt{MX}$ &$\texttt{MX}^{\dagger}$ & $\texttt{MX}$ &62 &6.9e+02\\ 
$\texttt{QUARTC}$ & 1250/5000 &70 &\emph{4.9}&104 &5.3&104 &5.4&70 &\textbf{4.5}&104 &5.1&104 &5.6&89 &1.4e+02\\ 
$\texttt{SCHMVETT}$ & 1250/5000 &$\texttt{MX}^{\dagger}$ & $\texttt{MX}$ &70882 &3.9e+03&$\texttt{MX}^{\dagger}$ & $\texttt{MX}$ &$\texttt{NC}^{\dagger}$ & $\texttt{NC}$ &$\texttt{MX}^{\dagger}$ & $\texttt{MX}$ &96572 &5.1e+03&$\texttt{NC}^{\dagger}$ & $\texttt{NC}$ \\ 
$\texttt{SINQUAD}$ & 1250/5000 &236 &56&282 &15&214 &\textbf{11}&247 &32&216 &\emph{12}&277 &14&116 &1.8e+02\\ 
$\texttt{SPARSQUR}$ & 2500/10000 &35 &13&43 &\emph{10}&43 &14&35 &13&43 &\textbf{9.9}&43 &14&31 &3.5e+02\\ 
$\texttt{SPMSRTLS}$ & 1250/4999 &2222 &2.7e+02&1791 &\textbf{95}&2377 &1.2e+02&2792 &2e+02&2475 &1.3e+02&1834 &\emph{98}&$\texttt{NC}^{\dagger}$ & $\texttt{NC}$ \\ 
$\texttt{SROSENBR}$ & 1250/5000 &5561 &4.1e+02&8211 &4.3e+02&4814 &\textbf{235}&6400 &4.3e+02&6747 &3.6e+02&5280 &\emph{270}&$\texttt{NC}^{\dagger}$ & $\texttt{NC}$ \\ 
$\texttt{TOINTGSS}$ & 1250/5000 &39 &3.1&39 &\textbf{2.2}&39 &2.3&39 &3&39 &2.3&39 &\emph{2.3}&49 &76\\ 
$\texttt{TQUARTIC}$ & 1250/5000 &2069 &8.7e+02&1155 &\textbf{64}&1508 &\emph{78}&1494 &3.7e+02&1867 &1e+02&1871 &98&$\texttt{NC}^{\dagger}$ & $\texttt{NC}$ \\ 
$\texttt{TRIDIA}$ & 1250/5000 &147 &9&82 &\textbf{4.2}&82 &4.3&147 &9.1&82 &\emph{4.2}&82 &4.4&66 &1e+02\\ 
$\texttt{WOODS}$ & 1000/4000 &1192 &45&1157 &38&1077 &\textbf{27}&1236 &44&1167 &37&1132 &\emph{29}&971 &1.3e+03\\ 
$\texttt{SPARSINE}$ & 1250/5000 &1504 &1.7e+02&1476 &79&1464 &\textbf{74}&2188 &1.6e+02&1407 &\emph{74}&3999 &2e+02&2294 &5.6e+03\\ 
$\texttt{TESTQUAD}$ & 1250/5000 &10988 &\textbf{623}&14186 &7.3e+02&13357 &6.5e+02&10988 &\emph{643}&14186 &7.3e+02&13357 &6.6e+02&$\texttt{NC}^{\dagger}$ & $\texttt{NC}$ \\ 
$\texttt{JIMACK}$ & 888/3549 &$\texttt{NC}^{\dagger}$ & $\texttt{NC}$ &$\texttt{NC}^{\dagger}$ & $\texttt{NC}$ &$\texttt{NC}^{\dagger}$ & $\texttt{NC}$ &$\texttt{NC}^{\dagger}$ & $\texttt{NC}$ &$\texttt{NC}^{\dagger}$ & $\texttt{NC}$ &$\texttt{NC}^{\dagger}$ & $\texttt{NC}$ &$\texttt{NC}^{\dagger}$ & $\texttt{NC}$ \\ 
$\texttt{NCB20B}$ & 1250/5000 &57 &4.1&56 &3.2&56 &3.2&57 &4.2&56 &\textbf{3.1}&56 &\emph{3.2}&47 &73\\ 
$\texttt{EIGENALS}$ & 638/2550 &202 &\emph{3.2}&204 &3.7&203 &4.1&202 &\textbf{3}&204 &3.6&203 &4&161 &43\\ 
$\texttt{EIGENBLS}$ & 638/2550 &28 &0.59&28 &0.65&28 &0.6&28 &\textbf{0.51}&28 &\emph{0.52}&28 &0.62&28 &7.7\\ 

 \hline
 \end{tabular}
 \hss}
\end{table} % Table 2

In Experiment II, the objective functions for each problem are defined by
a large CUTEst problem, whereas the corresponding $ \b{A} $ matrices are not meant
to be overly challenging. We observe that the proposed algorithms 
(the ones including ``\{H,L\}") improve the computation times on the majority of 
problems. For the 10 instances in which \TR2 used less time than \TR2H, the 
differences are relatively small. An exception is $ \texttt{DIXMAANL} $, where 
the difference amounts to 31s. However, for the other 51 problems, \TR2H resulted
in often significant improvements in computation time. For instance, in
$ \texttt{LIARWHD} $ this difference amounts to 1182s (more than 19 minutes).
These observations carry over when comparing \TR1 with \TR{1H}. The ``L" versions
exhibit similar outcomes as the ``H" ones, with occasional increases in computation
times. Overall, \TR{1L} converged to the specified tolerances on the largest number of problems. \jbk{The problems reported as ``\texttt{NC}" in {\small IPOPT }'s column correspond to status flags other than ``0, 1, 2" $\equiv$ ``solved, solved to acceptable level, infeasible problem detected".}

\subsection*{D.3: \jbr{Detailed Table for Experiment III}}
\label{sec:appD3}

% Table 3 here

\jbr{In Experiment III, \TR{2H} and \TR{1H} converged on all 31 problems,
while all other solvers (besides \TR1 and TR2) converged on all problems except one:
\texttt{CVXQP2}. \TR{2H} was the fastest on 10 problems (the best outcome among the solvers),
while \TR{1L} was the fastest on 9 problems (the second best outcome). Problems 
\texttt{A0ESDNDL} and \texttt{A0ESINDL} appear noteworthy:
%. In particular, these two problems each 
they contain dense columns (satisfying the condition $ \textnormal{nnz}(A_{:,j}) \big / m > 0.1 $).
%On these problems the 
Sparse QR factorization %uses long computing times 
is expensive because
of fill-in. However, the iterative method LSQR (with the preconditioning technique
from section \ref{sec:LSQR}) can overcome these difficulties.}
%\jbb{and include performance profiles in Appendix D}.
%\Red{and Figure~\ref{fig:EX_II}.}

\begin{table}[p]  % Table 3
 \captionsetup{width=1\linewidth}
\caption{\jbr{Experiment III compares 7 solvers on 31 linear equality constrained problems from the CUTEst collection \cite{GouOT03}. 
$\texttt{NC}^{\dagger}$ means the solver did not converge to tolerances.
$\texttt{N/A}$ means that \TR1 and \TR2 were not applied because the problem size rendered them not practical. \TR{2H} and \TR{1H} converged on all 31 problems.
\TR{2L}, \TR{1L}, and {\small IPOPT} converged on 30 problems (the exception is 
\texttt{CVXQP2}). The fastest and second fastest solvers for each problem
are highlighted in bold and italic fonts, respectively. Overall,
\TR{2H} was fastest on 12 problems (the best outcome on this experiment), while
\TR{1L} was fastest on 11 problems (the second best outcome). Problems
\texttt{A0ESDNDL} and \texttt{A0ESINDL} contain dense columns in $A$, and the sparse QR
factorization takes additional time as seen from the entries of \TR{2H} and \TR{1H}.
However, preconditioned {\small LSQR} can overcome
this difficulty, as observed in the entries for \TR{2L} and \TR{1L} for these problem
instances.
%$\texttt{MX}^{\dagger}$ means the iteration limit was reached.
%  \TR{1L} converged on 58 problems, the largest number of problems
% amongst the solvers. \TR{2H} was faster than \TR2 on 51 problems,
% and \TR{2L} was faster than \TR2 on 46 problems (the differences are often significant). \TR{1H} was faster than \TR{1} on 49 problems and
% \TR{1L} was faster than \TR{1} on 41 problems (often significantly).
% All of \TR{1\{H,L\}} and \TR{2\{H,L\}} were faster than
% {\small IPOPT}.}
}}
\label{tab:exp3}

% Adjusting table column separation
\setlength{\tabcolsep}{2pt} % Default value: 6pt

\scriptsize
\hbox to 1.00\textwidth{\hss % 1.25\textwidth
\begin{tabular}{|l | c | c c | c c | c c | c c | c c | c c | c c|} % tabular
		\hline
		\multirow{2}{*}{Problem} & \multirow{2}{*}{$m$/$n$} & %\multirow{2}{*}{$\text{rank}(\b{A})$/\text{Dty}} & 
% 		\multicolumn{2}{c|}{TR1}  & 
% 		\multicolumn{2}{c|}{TR1x1} &  		
% 		\multicolumn{2}{c|}{TR1x2} &  	
% 		\multicolumn{2}{c|}{TR2} &  	
% 		\multicolumn{2}{c|}{TR2x1} &  	
% 		\multicolumn{2}{c|}{TR2x2} &
        \multicolumn{2}{c|}{\jbbo{TR2}}  & 
		\multicolumn{2}{c|}{\jbbo{TR2H}} &  		
		\multicolumn{2}{c|}{\jbbo{TR2L}} &  	
		\multicolumn{2}{c|}{\jbbo{TR1}} &  	
		\multicolumn{2}{c|}{\jbbo{TR1H}} &  	
		\multicolumn{2}{c|}{\jbbo{TR1L}} & 
		\multicolumn{2}{c|}{IPOPT} \\  	
		\cline{3-16}
		& %& %&		
		& It &  Sec
		& It &  Sec 
		& It &  Sec 
		& It &  Sec 
		& It &  Sec 
		& It &  Sec 
		& It &  Sec \\
		\hline
%		\endhead
%		\hline \multicolumn{17}{r}{\textit{Continued on next page}} \\
%		\endfoot
%		\endlastfoot
$\texttt{AUG2D}$ & 10000/20200 &$\texttt{N/A}^{*}$ & $\texttt{N/A}$ &7 &0.26&7 &\emph{0.15}&$\texttt{N/A}^{*}$ & $\texttt{N/A}$ &7 &0.24&7 &\textbf{0.13}&12 &1.4\\ 
$\texttt{AUG2DC}$ & 10000/20200 &$\texttt{N/A}^{*}$ & $\texttt{N/A}$ &2 &0.11&2 &\emph{0.067}&$\texttt{N/A}^{*}$ & $\texttt{N/A}$ &2 &0.1&2 &\textbf{0.067}&1 &0.15\\ 
$\texttt{AUG2DCQP}$ & 10000/20200 &$\texttt{N/A}^{*}$ & $\texttt{N/A}$ &2 &0.11&2 &\emph{0.072}&$\texttt{N/A}^{*}$ & $\texttt{N/A}$ &2 &0.11&2 &\textbf{0.07}&1 &0.16\\ 
$\texttt{AUG2DQP}$ & 10000/20200 &$\texttt{N/A}^{*}$ & $\texttt{N/A}$ &7 &0.23&7 &\textbf{0.13}&$\texttt{N/A}^{*}$ & $\texttt{N/A}$ &7 &0.24&7 &\emph{0.13}&12 &1.4\\ 
$\texttt{AUG3D}$ & 8000/27543 &$\texttt{N/A}^{*}$ & $\texttt{N/A}$ &10 &0.68&10 &\emph{0.52}&$\texttt{N/A}^{*}$ & $\texttt{N/A}$ &10 &0.6&10 &\textbf{0.51}&11 &2.6\\ 
$\texttt{AUG3DC}$ & 8000/27543 &$\texttt{N/A}^{*}$ & $\texttt{N/A}$ &2 &0.3&2 &\emph{0.28}&$\texttt{N/A}^{*}$ & $\texttt{N/A}$ &2 &0.3&2 &\textbf{0.26}&1 &0.31\\ 
$\texttt{AUG3DCQP}$ & 8000/27543 &$\texttt{N/A}^{*}$ & $\texttt{N/A}$ &2 &0.3&2 &\emph{0.27}&$\texttt{N/A}^{*}$ & $\texttt{N/A}$ &2 &0.33&2 &\textbf{0.26}&1 &0.33\\ 
$\texttt{AUG3DQP}$ & 8000/27543 &$\texttt{N/A}^{*}$ & $\texttt{N/A}$ &10 &0.74&10 &\emph{0.55}&$\texttt{N/A}^{*}$ & $\texttt{N/A}$ &10 &0.64&10 &\textbf{0.5}&11 &2.6\\ 
$\texttt{CVXQP1}$ & 5000/10000 &$\texttt{N/A}^{*}$ & $\texttt{N/A}$ &827 &7.8&805 &\textbf{3.8}&$\texttt{N/A}^{*}$ & $\texttt{N/A}$ &827 &7.3&805 &\emph{3.8}&740 &51\\ 
$\texttt{CVXQP2}$ & 2500/10000 &$\texttt{N/A}^{*}$ & $\texttt{N/A}$ &39596 &1.5e+02&$\texttt{NC}^{\dagger}$ & $\texttt{NC}$ &$\texttt{N/A}^{*}$ & $\texttt{N/A}$ &47572 &1.8e+02&$\texttt{NC}^{\dagger}$ & $\texttt{NC}$ &$\texttt{NC}^{\dagger}$ & $\texttt{NC}$ \\ 
$\texttt{CVXQP3}$ & 7500/10000 &$\texttt{N/A}^{*}$ & $\texttt{N/A}$ &169 &2.8&169 &\emph{1.4}&$\texttt{N/A}^{*}$ & $\texttt{N/A}$ &169 &2.4&169 &\textbf{1.4}&118 &8.9\\ 
$\texttt{STCQP1}$ & 4095/8193 &$\texttt{N/A}^{*}$ & $\texttt{N/A}$ &88 &\textbf{0.15}&88 &0.42&$\texttt{N/A}^{*}$ & $\texttt{N/A}$ &88 &\emph{0.18}&88 &0.36&75 &6.8e+02\\ 
$\texttt{STCQP2}$ & 4095/8193 &$\texttt{N/A}^{*}$ & $\texttt{N/A}$ &142 &\textbf{0.25}&142 &0.8&$\texttt{N/A}^{*}$ & $\texttt{N/A}$ &144 &\emph{0.28}&144 &0.72&136 &4.8\\ 
$\texttt{DTOC1L}$ & 3996/5998 &$\texttt{N/A}^{*}$ & $\texttt{N/A}$ &13 &\textbf{0.073}&13 &0.13&$\texttt{N/A}^{*}$ & $\texttt{N/A}$ &13 &\emph{0.075}&13 &0.14&16 &0.41\\ 
$\texttt{DTOC3}$ & 2998/4499 &$\texttt{N/A}^{*}$ & $\texttt{N/A}$ &5 &\textbf{0.025}&5 &0.059&$\texttt{N/A}^{*}$ & $\texttt{N/A}$ &5 &\emph{0.03}&5 &0.033&4 &0.09\\ 
$\texttt{PORTSQP}$ & 1/100000 &2 &0.09&2 &0.064&2 &0.067&2 &\emph{0.062}&2 &\textbf{0.059}&2 &0.062&1 &0.42\\ 
$\texttt{HUES-MOD}$ & 2/5000 &1 &0.0028&1 &\textbf{0.0018}&1 &0.0027&1 &0.0026&1 &\emph{0.0018}&1 &0.0026&1 &0.024\\ 
$\texttt{HUESTIS}$ & 2/5000 &2 &0.0073&2 &\textbf{0.0042}&2 &0.011&2 &0.0061&2 &\emph{0.0047}&2 &0.0094&2 &0.072\\ 
$\texttt{A0ESDNDL}$ & 15002/45006 &$\texttt{N/A}^{*}$ & $\texttt{N/A}$ &5 &69&5 &\emph{0.13}&$\texttt{N/A}^{*}$ & $\texttt{N/A}$ &5 &71&5 &\textbf{0.12}&6 &1.8\\ 
$\texttt{A0ESINDL}$ & 15002/45006 &$\texttt{N/A}^{*}$ & $\texttt{N/A}$ &5 &73&5 &\emph{0.12}&$\texttt{N/A}^{*}$ & $\texttt{N/A}$ &5 &70&5 &\textbf{0.11}&6 &1.8\\ 
$\texttt{PORTSNQP}$ & 2/100000 &$\texttt{NC}^{\dagger}$ & $\texttt{NC}$ &2 &\textbf{0.092}&2 &0.11&14 &0.47&2 &\emph{0.095}&2 &0.1&2 &0.88\\ 
$\texttt{BLOWEYA}$ & 2002/4002 &$\texttt{N/A}^{*}$ & $\texttt{N/A}$ &2 &\textbf{0.011}&2 &0.031&$\texttt{N/A}^{*}$ & $\texttt{N/A}$ &2 &\emph{0.015}&2 &0.021&2 &0.082\\ 
$\texttt{BLOWEYB}$ & 2002/4002 &$\texttt{N/A}^{*}$ & $\texttt{N/A}$ &2 &\textbf{0.015}&2 &0.019&$\texttt{N/A}^{*}$ & $\texttt{N/A}$ &2 &\emph{0.016}&2 &0.019&2 &0.082\\ 
$\texttt{BLOWEYC}$ & 2002/4002 &$\texttt{N/A}^{*}$ & $\texttt{N/A}$ &2 &\textbf{0.015}&2 &0.017&$\texttt{N/A}^{*}$ & $\texttt{N/A}$ &2 &\emph{0.015}&2 &0.021&2 &0.15\\ 
$\texttt{CONT5-QP}$ & 40200/40601 &$\texttt{N/A}^{*}$ & $\texttt{N/A}$ &2 &\emph{0.51}&2 &0.79&$\texttt{N/A}^{*}$ & $\texttt{N/A}$ &2 &\textbf{0.49}&2 &0.8&2 &1.3\\ 
$\texttt{DTOC1L}$ & 3996/5998 &$\texttt{N/A}^{*}$ & $\texttt{N/A}$ &5 &\textbf{0.03}&5 &0.09&$\texttt{N/A}^{*}$ & $\texttt{N/A}$ &5 &\emph{0.043}&5 &0.045&4 &0.12\\ 
$\texttt{FERRISDC}$ & 210/2200 &2 &0.084&2 &0.083&2 &0.077&2 &\emph{0.076}&2 &0.078&2 &0.079&0 &\textbf{0.021}\\ 
$\texttt{GOULDQP2}$ & 9999/19999 &$\texttt{N/A}^{*}$ & $\texttt{N/A}$ &2 &0.038&2 &\textbf{0.025}&$\texttt{N/A}^{*}$ & $\texttt{N/A}$ &2 &0.038&2 &\emph{0.026}&2 &0.2\\ 
$\texttt{GOULDQP3}$ & 9999/19999 &$\texttt{N/A}^{*}$ & $\texttt{N/A}$ &6 &0.076&6 &\emph{0.054}&$\texttt{N/A}^{*}$ & $\texttt{N/A}$ &6 &0.077&6 &\textbf{0.053}&7 &0.69\\ 
$\texttt{LINCONT}$ & 419/1257 &5 &0.058&5 &\emph{0.02}&5 &0.031&5 &0.05&5 &\textbf{0.019}&5 &0.03&5 &0.055\\ 
$\texttt{SOSQP2}$ & 2501/5000 &$\texttt{N/A}^{*}$ & $\texttt{N/A}$ &3 &\textbf{0.017}&3 &0.04&$\texttt{N/A}^{*}$ & $\texttt{N/A}$ &3 &0.022&3 &\emph{0.019}&4 &0.11\\ 
 \hline
 \end{tabular}
 \hss}
\end{table} % Table 3

% \jbr{In Experiment II, the objective functions for each problem are defined by
% a large CUTEst problem, whereas the corresponding $ \b{A} $ matrices are not meant
% to be overly challenging. We observe that the proposed algorithms 
% (the ones including ``\{H,L\}") improve the computation times on the majority of 
% problems. For the 10 instances in which \TR2 used less time than \TR2H, the 
% differences are relatively small. An exception is $ \texttt{DIXMAANL} $, where 
% the difference amounts to 31s. However, for the other 51 problems, \TR2H resulted
% in often significant improvements in computation time. For instance, in
% $ \texttt{LIARWHD} $ this difference amounts to 1182s (more than 19 minutes).
% These observations carry over when comparing \TR1 with \TR{1H}. The ``L" versions
% exhibit similar outcomes as the ``H" ones, with occasional increases in computation
% times. Overall, \TR{1L} converged to the specified tolerances on the largest number of problems. The problems reported as ``\texttt{NC}" in IPOPT's column correspond to status flags other than ``0, 1, 2" $\equiv$ ``solved, solved to acceptable level, infeasible problem detected".}

\clearpage

\section*{Acknowledgments}
We would like to acknowledge the valuable discussions initiated by Ariadna Cairo Baza and
spurred by the 9th ICIAM conference at the Universidad de Valencia.
\jbb{R. Marcia’s research was partially supported by NSF Grant IIS 1741490.}
\jbrt{We thank two referees for
%detailed comments that helped improve the presentation of the paper.}
their extremely detailed and helpful comments.}
%\pagebreak

\bibliographystyle{siamplain}
\bibliography{myrefs}

\medskip

\begin{center}
\fbox{%
\parbox{3in}{\scriptsize
The submitted manuscript has been created by UChicago Argonne, LLC, Operator of Argonne 
National Laboratory (``Argonne''). Argonne, a U.S. Department of Energy Office of Science 
laboratory, is operated under Contract No. DE-AC02-06CH11357. The U.S. Government retains 
for itself, and others acting on its behalf, a paid-up nonexclusive, irrevocable worldwide 
license in said article to reproduce, prepare derivative works, distribute copies to the 
public, and perform publicly and display publicly, by or on behalf of the Government. 
The Department of Energy will provide public access to these results of federally 
sponsored research in accordance with the DOE Public Access Plan.
\url{http://energy.gov/downloads/doe-public-accessplan}
}}
\end{center}

\end{document}